\documentclass[12pt]{article}
\usepackage{amsmath}
\usepackage{amssymb}
\usepackage{amsthm}
\usepackage[totalwidth=17cm, totalheight=25cm]{geometry}
\usepackage[all]{xy}

\newtheoremstyle{boldplain}
  {9pt}
  {9pt}
  {\itshape}
  {}
  {\bfseries}
  {.}
  {.5em}
  {\thmname{#1}\thmnumber{ #2}\thmnote{ (#3)}}%

\newtheoremstyle{bolddefinition}
  {9pt}
  {9pt}
  {}
  {}
  {\bfseries}
  {.}
  {.5em}
  {\thmname{#1}\thmnumber{ #2}\thmnote{ (#3)}}%

\theoremstyle{boldplain}

\newtheorem{thm}[equation]{Theorem}
\newtheorem{prop}[equation]{Proposition}

\newtheorem{ass}[equation]{Assumption}
\newtheorem{lem}[equation]{Lemma}

\newtheorem{cor}[equation]{Corollary}

\theoremstyle{bolddefinition}
\newtheorem{defn}[equation]{Definition}

\newtheorem{rem}[equation]{Remark}

\newfont{\bigbf}{cmbx10 scaled\magstep1}
\normalbaselines
\numberwithin{equation}{section}

\def\no{\noindent}

\def\R{{\mathbb R}}
\def\H{{\mathbb H}}

\def\N{{\mathbb N}}

\def\Z{{\mathbb Z}}

\def\al{\alpha}
\def\ga{\gamma}
\def\Ga{\Gamma}
\def\de{\delta}
\def\De{\Delta}
\def\eps{\epsilon}
\def\la{\lambda}

\def\si{\sigma}
\def\Si{\Sigma}
\def\om{\omega}

\def\acts{\curvearrowright}
\def\D{\partial}

\def\half{\frac{1}{2}}
\def\third{\frac{1}{3}}
\def\quart{\frac{1}{4}}
\def\3quart{\frac{3}{4}}
\def\fifth{\frac{1}{5}}
\def\sixth{\frac{1}{6}}
\def\seventh{\frac{1}{7}}
\def\tenth{\frac{1}{10}}
\def\eleventh{\frac{1}{11}}
\def\hunth{\frac{1}{100}}
\def\lra{\longrightarrow}

\def\pihalf{\frac{\pi}{2}}
\def\pithird{\frac{\pi}{3}}
\def\piquart{\frac{\pi}{4}}

\def\pitenth{\frac{\pi}{10}}
\def\pihunth{\frac{\pi}{100}}

\def\ol{\overline}

\def\diam{\mathop{\hbox{diam}}}

\def\rad{\mathop{\hbox{rad}}}
\def\cangle{\tilde\angle}
\def\tangle{\angle_{Tits}}
\def\tits{\mathop{\partial_{Tits}}}
\def\vol{\mathop{\hbox{vol}}}

\def\BI{\begin{itemize}}
\def\EI{\end{itemize}}

\title{On the topology of locally volume collapsed Riemannian 3-orbifolds}

\author{Daniel Faessler}

\date{January 19, 2011}

\begin{document}

\maketitle
\begin{abstract}
\noindent
We study the geometry and topology
of Riemannian 3-orbifolds
which are locally volume collapsed
with respect to a curvature scale.
We show that a
sufficiently collapsed closed 3-orbifold
without bad 2-suborbifolds
either admits a metric of nonnegative sectional curvature
or satisfies Thurston's Geometrization Conjecture.
We also prove a version with boundary.
Kleiner and Lott independently proved similar results
\cite{KL_preprint}.
\end{abstract}

\tableofcontents

\section{Introduction}

We study the geometry and topology
of Riemannian 3-orbifolds
which are locally volume collapsed
with respect to a curvature scale.
Such orbifolds are expected
to occur as the {\em thin} part
of solutions of the 
orbifold version of Perelman's Ricci flow with surgery
(as constructed on manifolds in \cite{Perelman})
after sufficiently long time.
Our main result (Theorem \ref{thm:mainclosed})
concerns the topology
of locally collapsed 3-orbifolds.
We show that a
sufficiently collapsed closed 3-orbifold
without bad 2-suborbifolds
either admits a metric of nonnegative sectional curvature
or satisfies Thurston's Geometrization Conjecture,
i.e.\ has a connected sum decomposition
(by spherical surgeries)
into components which in turn
admit a (toric) JSJ-decomosition
into {\em geometric} components.
Closed 3-orbifolds with nonnegative sectional curvature
are also expected to be geometric
by an orbifold version of
Hamilton's corresponding result
for 3-manifolds \cite{Hamilton1}.

In order to avoid the use of
an orbifold version of
Perelman's Stability Theorem
we prove our result
under an additional regularity assumption.
We require uniform control on the
derivatives of the curvature tensor
up to some (sufficiently large) finite order.

We also prove a version with boundary
of our result (Theorem \ref{thm:mainbdry}).
We expect that the assumptions are sufficiently general
to apply to the thin part
of a solution of the orbifold Ricci flow with surgery
if the thick-thin decomposition is nontrivial
(i.e.\ in a situation of partial collapse).

Corresponding results for 
collapsed orientable 3-manifolds
have been stated without proof in \cite{Perelman}
and proved in
\cite{SY05}, \cite{MT2}, \cite{KL_coll} and \cite{BBBMP}.
After writing this paper,
we learned that Kleiner and Lott
independently proved
results similar to our main result,
cf.\ \cite[Prop.\ 9.7]{KL_preprint}.
Their method is an extension
of their work \cite{KL_coll} in the manifold case
to the orbifold case,
whereas our approach 
is closer to an extension
of the approach in \cite{MT2}.

\medskip

The paper is organized as follows:
In section 2, we first review basic facts on orbifolds
in low dimensions.
We then discuss decompositions of 3-orbifolds
along spherical and toric 2-suborbifolds and prove
that {\em graph} orbifolds in the sense of Waldhausen
(cf. section \ref{sec:decos})
satisfy Thurston's Geometrization Conjecture
(Corollary \ref{cor:graphgeom}).

In the third section,
we discuss a coarse stratification
of roughly 2-dimensional Alexandrov spaces.
More precisely,
we use a conical approximation argument
to show that the points in such a space
which do not admit 1-strainers of a certain length and quality
accumulate in isolated regions.
Outside these regions,
the Alexandrov space is 1-strained
which allows us to perform
a (coarse) dimension reduction
by considering cross sections to these strainers.
We further distinguish points
according to whether they lie in coarse necks,
edges or the interior of the Alexandrov space
and study their geometric properties.
These considerations are similar in spirit
to considerations in \cite{MT2} and \cite{KL_coll}.

In section 4,
we restrict our attention
to closed volume collapsed 3-orbifolds.
We consider them as
Alexandrov spaces which are roughly of dimension $\le 2$
and promote their coarse stratification
to a certain decomposition into 3-suborbifolds.
To determine the local topology
of the components in this decomposition,
we use a variation
(and extension to additional situations)
of the blow-up arguments in \cite{SY}.
We derive a graph decomposition of
the collapsed 3-orbifolds.
Combined with the results of section 2,
the main result follows. 

The fifth section
contains a generalization of Theorem \ref{thm:mainclosed}
to compact orbifolds with boundary
where we require that neighbourhoods of the boundary
are close to pieces of hyperbolic cusps
(cf. Theorem \ref{thm:mainbdry}).

\medskip

{\em Acknowledgments}.
The author would like to thank
Professor Bernhard Leeb
for his guidance and support
and Richard Bamler for
helpful discussions.
He is grateful for funding by
Studienstiftung des Deutschen Volkes.

\section{Decompositions of 3-orbifolds along 2-suborbifolds}

\subsection{Orbifolds}
We refer to \cite[sec.\ 2]{BMP} for a more detailed discussion of orbifolds.

\subsubsection{Smooth orbifolds}

Roughly speaking, an orbifold is a space which looks locally like the orbit space of a linear action by a finite (orthogonal) group. Important examples are the orbit spaces of properly discontinuous group actions on manifolds, in particular of finite group actions. These orbifolds are called {\em good}, respectively, {\em very good}. To exclude exotic local phenomena, e.g.\ when working in the topological category, one has to require the group actions to be locally linearizable; this is automatic in the smooth category.

More formally, 
an $n$-dimensional {\em smooth orbifold} $O$ 
is a metrizable topological space 
together with a maximal atlas of $n$-dimensional orbifold charts satisfying certain compatibility conditions. An {\em orbifold chart} $(U,\tilde{U}, \Ga_U, \pi_U)$ consists of an open subset $U \subseteq O$, a smooth $n$-manifold $\tilde{U}$, a finite subgroup $\Ga_U \subset \textrm{Diff}(\tilde{U})$, and a continuous map $\pi_U: \tilde{U} \to U$ inducing a homeomorphism $\tilde{U}/\Ga_U \buildrel \cong \over \to U$.

Any two charts $(U_i,\tilde{U}_i, \Ga_{U_i},\pi_{U_i})$, $i=1,2$, must be compatible in the following sense: If $\tilde{x}_i \in \tilde{U}_i$ are points with $\pi_{U_1}(\tilde{x}_1)=\pi_{U_2}(\tilde{x}_2)$, then there exists a diffeomorphism $\tilde{\phi}:\tilde{V}_1 \to \tilde{V}_2$ of open neighborhoods $\tilde{V}_i$ of the $\tilde{x}_i$ with $\pi_{U_2} \circ \tilde{\phi} = \pi_{U_1}$. Finally, the charts must {\em cover} $M$.

We will denote by $|O|$ the topological space underlying the orbifold $O$.

More generally, we define smooth $n$-dimensional orbifolds $O$ {\em with boundary} by allowing the chart domains $\tilde U$ to be smooth $n$-manifolds with boundary.

The {\em boundary} $\D O$, respectively, the {\em interior} of $O$ consist of those points whose preimages in the chart domains are boundary, respectively, interior points. Since the local coordinate changes $\tilde\phi$ are smooth, they preserve boundaries. Consequently, $\D O$ is a closed subset and inherits a structure as a smooth $(n-1)$-orbifold without boundary. The boundary has a collar, i.e.\ $\D O$ has an open neighborhood in $O$ diffeomorphic to the product orbifold $\D O\times[0,1)$ where $[0,1)$ is to be understood as a 1-manifold with boundary. 

Note that the local coordinate change $\tilde\phi$ in the above definition must be equivariant with respect to an isomorphism $Stab_{\Ga_{U_1}}(\tilde{x}_1)\to Stab_{\Ga_{U_2}}(\tilde{x}_2)$, i.e.\ the local actions of the stabilizers $Stab_{\Ga_{U_i}}(\tilde{x}_i)$ near $\tilde{x}_i$ are conjugate. This is due to the fact that smooth orbit equivalences 
between effective finite smooth group actions on connected manifolds are conjugacies.

Given a point $x\in O$, we can choose a chart $(U,\tilde{U}, \Ga_U, \pi_U)$ with $x\in U$ and a preimage $\tilde{x} \in \pi_U^{-1}(x)$. By our previous remark, the germ of the action of $\Ga_x:=Stab_{\Ga_U}(\tilde x)$ near $\tilde x$ is independent of these choices, and so is its linearisation $\Ga_x\acts T_{\tilde x}\tilde U$. We may hence regard $\Ga_x$ as a subgroup of $O(n)$ 
well-defined up to conjugacy; it is called the {\em local group} of $O$ at $x$. 
We will call the quotient $T_{\tilde{x}}\tilde{U}/\Gamma_x$ 
the orbifold {\em tangent space} $T_{x}O$ of $O$ at $x$.
The orbifold $O$ is called {\em locally orientable} at the point $x$ if $\Ga_x\subset SO(n)$.

A point $x \in O$ is said to be {\em regular} if its local group $\Gamma_x$ is trivial, and {\em singular} otherwise. The subset $O^{sing}$ of singular points is the {\em singular locus}. (An orbifold is a manifold if and only if all its points are regular, i.e.\ if $O^{sing}=\emptyset$.) We call the conjugacy class of the local group $\Ga_x\subset O(n)$ the (singular) {\em type} of the point $x$. One observes that, in a chart around $x$, the fixed point set of $\Ga_x$ in $\tilde U$ is a submanifold and its connected component through $\tilde x$ projects to points in $O$ of the same singular type as $x$. Hence the equivalence classes of points of the same type inherit structures as smooth manifolds. Together they form a natural {\em stratification} of the orbifold $O$. We refer to the union of the $d$-dimensional strata as the singular {\em $d$-stratum} $O^{(d)}$. Note that $\overline{O^{(d)}} \setminus O^{(d)} \subset \bigcup_{k < d} O^{(k)}$. If $O$ has boundary, then the singular strata are manifolds with boundary and one has that $\partial O^{(d)} = (\partial O)^{(d-1)}$.

Clearly, the top-dimensional stratum $O^{(n)}$ consists precisely of the regular points. The singular $(n-1)$-stratum $O^{(n-1)}$ consists of the points with local group $\cong \mathbb{Z}_2$ generated by a hyperplane reflection. Its closure $\D_{refl}O:=\ol{O^{(n-1)}}$ is usually referred to as {\em reflector boundary} or {\em silvered boundary}, even though it is not contained in the boundary, $\D_{refl}O \cap \D O=\D_{refl}\D O$. It consists of the points whose local group contains a hyperplane reflection. We call $O^{(n-1)}$ the regular part of the reflector boundary.

Note that the underlying topological space $|O|$ of the orbifold $O$ 
contains only partial information on the orbifold structure,
e.g.\ on the singular stratification.
For instance, $|O|$ can be a topological manifold, 
although $O^{sing}\neq\emptyset$. 

We define an $m$-dimensional smooth {\em suborbifold} of a smooth orbifold 
to be a subset whose preimages in local charts are 
smooth $m$-dimensional submanifolds, cf.\ \cite[sec.\ B]{BS_split}. 
(This is more restrictive than other definitions used in the literature, 
cf.\ e.g.\ \cite[2.1.3]{BMP}). 
Analogously, a subset of an orbifold with boundary is called 
a {\em proper} suborbifold if its preimages in local charts 
are proper (smooth) submanifolds. 

A proper codimension-one suborbifold $\Si^{d-1} \subset O^d$ 
is called {\em two-sided} 
if it has a product neighborhood of the form $\Si\times(-1,1)$. 
It is called {\em locally two-sided} at a point $x$, 
if it has such a product neighborhood locally near $x$. 
If $\Si$ is not (globally) two-sided 
then it has a tubular neighborhood of the form 
$(\Si'\times(-1,1))/\Z_2$ 
where $\Z_2$ reflects on $(-1,1)$
and acts by a (possibly trivial) involution on $\Si'$.
(To verify this,
take e.g.\ $\Si'$ as the boundary of a tubular neighborhood (thickening)
of $\Si$.)

When speaking of a codim-zero suborbifold $O'\subset O$ 
we suppose that the components of $\D O'$ 
are either components of $\D O$ 
or disjoint from $\D O$, i.e.\ two-sided suborbifolds of $int(O)$. 

A continuous map $f:O\to O'$ of smooth orbifolds is called {\em smooth} (in the orbifold sense) if it lifts locally to an equivariant smooth (in the manifold sense) map of charts, i.e.\ if the following holds: For any point $x\in O$ exist charts $(U,\tilde{U},\Ga_U,\pi_U)$ around $x$ and $(U',\tilde U',\Ga_{U'},\pi_{U'})$ around $x'=f(x)$, and a smooth local lift $\tilde f_U:\tilde U\to\tilde U'$ of $f$, $f\circ\pi_U=\pi_{U'}\circ\tilde f_U$, which is equivariant with respect to a group homomorphism $\rho: \Gamma_U \to \Gamma_{U'}$. (Note that in general the equivariance of lifts is not automatic; there may exist non-equivariant lifts.)


A smooth map is a {\em submersion} ({\em immersion}) if it lifts locally to an equivariant submersion (immersion) of charts. 
It is a {\em local diffeomorphism} 
if it can be locally inverted everywhere by a smooth map, equivalently, if it lifts locally to diffeomorphisms of charts which are equivariant with respect to isomorphisms of local groups.
A global {\em diffeomorphism} is a homeomorphism which is also a local diffeomorphism.

A smooth map $p:O'\to O$ between orbifolds is a {\em covering} if every point $x\in O$ has an open neighborhood $U$ such that for every connected component $U'$ of $p^{-1}(U)$ exists a chart $(U',\tilde U',\Ga_{U'},\pi_{U'})$ for $O'$ and a possibly larger finite group $\Ga_U$, $\Ga_{U'}\subseteq\Ga_U\subset Diff(\tilde U')$, such that $(U,\tilde U',\Ga_U,p\circ\pi_{U'})$ is a chart for $O$. Coverings induce injective homomorphisms of the local groups which are well-defined up to postcomposition with inner automorphisms.

An orbifold is called {\em (very) good} if it is (finitely) covered by a manifold; an orbifold which is not covered by a manifold is called {\em bad}.

Let $F$ be an orbifold without boundary. Following \cite[sec.\ 2.4]{BMP} we define an {\em orbifold fiber bundle} or {\em orbifold fibration} with generic fiber $F$ as a submersion $p:O\to B$ of orbifolds, possibly with boundary, with the following property: For every point $x \in B$ exists a chart $\phi: \tilde{U} \to U$ around $x$, a smooth operation $\Gamma_x\acts F$ and a submersion $\si: \tilde{U} \times F \to O$ inducing a diffeomorphism between $(\tilde{U} \times F)/\Gamma_x$ (where we divide out the diagonal action) and $p^{-1}(U)$ such that $p \circ \si = \phi \circ \pi_{\tilde U}$. In the case with boundary we require that $p^{-1}(\D B)=\D O$;
then $p$ restricts over the boundary to the orbifold fiber bundle $p_{|\D O}:\D O\to\D B$. 
Note that orbifold coverings are 
(the same as) 
orbifold fiber bundles with 0-dimensional fiber. 

We say that a compact orbifold {\em fibers}
if it is the total space of an orbifold fibration 
whose base and fiber have strictly positive dimension 
and whose generic fiber is a closed orbifold.

An orbifold is called 
{\em spherical (discal, toric, solid toric)}
if it is diffeomorphic to the quotient of a 
round sphere $S^n$
(a closed unit disc $D^n$, a flat 2-torus, 
the compact 3-dimensional solid torus ($D^2\times S^1$))
by a finite isometric group action. 

\subsubsection{Riemannian and geometric orbifolds}
\label{sec:riemorb}

A {\em Riemannian orbifold} can be defined as a smooth orbifold together with compatible Riemannian metrics on the local uniformizations $\tilde{U}$, i.e.\ such that the operations $\Gamma_U\acts\tilde U$ and the local coordinate changes $\tilde{\phi}$ are isometric.

Notions from Riemannian geometry like lengths of curves, path metric, geodesics, exponential map directly generalize to orbifolds via local charts. The natural stratification of $O$ is totally geodesic. 

The smooth orbifold structure underlying a Riemannian orbifold 
is encoded in its metric structure. 
The singular points in the orbifold sense are also geometric singularities, 
and a Riemannian orbifold can be defined more directly as a metric space 
which is locally isometric to the quotient of a Riemannian manifold 
by a finite group of isometries.

\medskip

A {\em geometric structure} on a smooth 3-orbifold 
is a Riemannian metric 
which is modelled 
on one of the eight 3-dimensional Thurston geometries 
$S^3$, $\R^3$, $\H^3$,
$S^2 \times \R$, $\H^2\times\R$, 
$Nil$, $\widetilde{PSL(2,\R)}$ or $Solv$. 
This means that the Riemannian metrics on the local uniformizations
are locally isometric to the respective model space,
equivalently, 
the orbifold Riemannian metric is everywhere locally isometric 
to a quotient of the model space
by a finite group of isometries, 
cf.\ sec.~\ref{sec:riemorb}.

A compact 3-orbifold is called {\em geometric}
if its interior admits a complete geometric structure.
For a closed geometric 3-orbifold the model geometry is unique.

Geometric orbifolds are very good,
i.e.\ finitely covered by manifolds. 


\subsection{Low-dimensional orbifolds}

\subsubsection{1-orbifolds}

There are two connected closed 1-orbifolds, 
namely the circle and the {\em mirrored interval} $\bar I$. 
The latter has as underlying topological space $|\bar I|$ 
the compact interval $I$. 
The boundary points of $|\bar I|$ 
are reflector boundary points of $\bar I$, 
$\D_{refl}\bar I=\bar I^{(0)}=\D|\bar I|$, 
but no boundary points; 
their local groups are $O(1)\cong\Z_2$.

\subsubsection{2-orbifolds}
\label{sec:2orb}

We will denote by $D^2$ the {\em closed} 2-disc
and by $\ol D^2$ the closed 2-disc with {\em reflector boundary}. 

Let $O^2$ be a 2-orbifold, possibly with boundary. A singular point with local group $\cong\Z_p$, $p\geq2$, acting by rotations is called a {\em cone point} of order $p$, respectively in the Riemannian case, with {\em cone angle} $\frac{2\pi}{p}$. It is an isolated singular point in the interior of $O$ and has a neighborhood diffeomorphic to the disc $D^2(p):=D^2/\Z_p$ with cone point of order $p$.

A singular point whose local group $\Ga_x\subset O(2)$ is a dihedral group $D_q$, $q\geq2$, is called a {\em corner} vertex of order $q$, respectively in the Riemannian case, with angle $\frac{\pi}{q}$.\footnote{As a subgroup of $O(2)$, the dihedral group $D_q$ is defined as the isometry group of a regular $q$-gon. It is generated by the reflections at two lines through the origin with angle $\frac{\pi}{q}$. As an abstract group, it has the presentation $\langle s_1,s_2|s_1^2=s_2^2=(s_1s_2)^q=1\rangle$.} It is an interior point and has a neighborhood diffeomorphic to the sector $V^2(q):=D^2/D_q$. 
Note that $\D_{refl}V^2(q)=\D|V^2(q)|$.

A singular point with local group $\cong D_1\cong\Z_2$ acting by a reflection (on the disc or half-disc) is a regular reflection boundary point and has a neighborhood diffeomorphic to $V^2(1):=D^2/D_1$. It may be a boundary point, namely one of the points in $\D O\cap\D_{refl}O$.

The singular locus $O^{sing}=O^{(1)}\cup O^{(0)}$ consists of the reflector boundary $\D_{refl}O=\ol{O^{(1)}}$, which contains the set $\ol{O^{(1)}}-O^{(1)}=\ol{O^{(1)}}\cap O^{(0)}$ of corners, and of the set $O^{(0)}-\ol{O^{(1)}}$ of cone points. We call a connected component of $O^{(1)}$ a {\em reflector edge}. In a corner vertex, locally two reflector edges meet.

Sometimes we will also admit cone points and corner vertices of order 1 
which are nothing else than regular interior, respectively, 
regular reflector boundary points. 

A connected component of $\D|O|$ 
can be a connected component of $\D O$ or of $\D_{refl}O$,
or it can be a chain of consecutive boundary arcs and reflector edges.
In the latter case, 
any two of the boundary arcs are disjoint,
but there may be sequences of consecutive reflector edges 
meeting at corner vertices. 

We will use the following notation. 
For a 2-manifold $\Si$ with boundary, we denote by 
$\bar\Si(p_1, \hdots, p_k; q_1, \hdots, q_l)$ a 2-orbifold (without boundary) with underlying space $\Si$, reflector boundary $\D\Si$,
$k$ cone points of orders $p_i$ located in the interior of $\Si$ and $l$ corner vertices of orders $q_j$ lying on $\D\Si$. If there are no cone points ($k=0$) we write $\Si(; q_1, \hdots, q_l)$; if there are no corner vertices ($l=0$) we will sometimes write briefly $\Si(p_1, \hdots, p_k;)=:\Si(p_1, \hdots, p_k)$, and also $\Si(;)=:\Si$. In most instances where we will apply this notation the diffeomorphism type of the 2-orbifold is uniquely determined. Note in particular that $\ol D^2(p;)=\ol D^2(p)$ is a closed $2$-orbifold with a reflector boundary circle.

The {\em Euler characteristic} of a 2-orbifold $O = \Sigma(p_1, \hdots, p_k; q_1, \hdots, q_l)$ is given by
\begin{equation*} 
\chi(O) = \chi(\Sigma) - \sum_{i=1}^k (1 - \frac{1}{p_i}) - \frac{1}{2} \sum_{j=1}^l(1 - \frac{1}{q_j}).
\end{equation*}
We recall the classification of connected closed spherical 2-orbifolds, that is, of quotients of the 2-sphere: $S^2$, $\R P^2$, $\ol D^2$, $S^2(p,p)$, $\R P^2(p)$, $\ol D^2(p)$, $\ol D^2(;p,p)$, $S^2(2,2,p)$, $S^2(2,3,3)$, $S^2(2,3,4)$, $S^2(2,3,5)$, $\ol D^2(;2,2,p)$, $\ol D^2(;2,3,3)$, $\ol D^2(;2,3,4)$, $\ol D^2(;2,3,5)$, $\ol D^2(2;p)$ and $\ol D^2(3;2)$ with $p\geq2$.

The list of connected closed flat 2-orbifolds, i.e. of quotients of the 2-torus is: $T^2$, $K^2$, $\textrm{Ann}^2$, $\textrm{M\"ob}^2$, $S^2(2,3,6)$, $S^2(2,4,4)$, $S^2(3,3,3)$, $S^2(2,2,2,2)$, $\mathbb{R}P^2(2,2)$, $\ol D^2(;2,3,6)$, $\ol D^2(;2,4,4)$, $\ol D^2(;3,3,3)$, $\ol D^2(;2,2,2,2)$, $\ol D^2(4;2)$, $\ol D^2(3;3)$, $\ol D^2(2;2,2)$ and $\ol D^2(2,2)$. 

All connected bad 2-orbifolds are closed
and can be obtained by gluing two non-diffeomorphic discal 2-orbifolds 
along their boundaries. 
They are diffeomorphic to 
$S^2(p)$, $S^2(p,q)$, $\ol D^2(;p)$ or $\ol D^2(;p,q)$ with $2\leq p<q$. 

\begin{prop}
A connected closed $2$-orbifold admits a Riemannian metric with nonnegative sectional curvature if and only if it has Euler characteristic $\ge 0$ if and only if it is spherical, flat or bad.
\end{prop}

\subsubsection{3-orbifolds}

Let $O^3$ be a 3-orbifold, possibly with boundary. We first discuss the local structure of the singular locus. 

If $x\in O^{(1)}$, then the local group $\Ga_x\subset O(3)$ fixes a line, i.e.\ $\Ga_x\cong\Z_p$ or $D_p$ with $p\geq2$ and $S^2/\Ga_x\cong S^2(p,p)$ or $\ol D^2(;p,p)$. We call the connected component of $O^{(1)}$ containing $x$ a {\em singular edge}, respectively, {\em reflector edge} (or circle) of order $p$. In the reflector case $\Ga_x\cong D_p$, locally two reflector faces, that is, components of $O^{(2)}$ meet at the edge. 
The boundary points of the singular and reflector edges (i.e.\ of their underlying 1-manifolds) are the cone points and corner vertices of the 2-orbifold $\D O$. 

If $x\in O^{(0)}$, then we call $x$ a {\em singular vertex}. In this case, $\Ga_x\subset O(3)$ has no nontrivial fixed vector and $S^2/\Ga_x$ is a spherical 2-orbifold with diameter $<\pi$, i.e.\ isometric to $\R P^2$, $\R P^2(p)$, $\ol D^2(p)$, $S^2(2,2,p)$, $S^2(2,3,3)$, $S^2(2,3,4)$, $S^2(2,3,5)$, $\ol D^2(;2,2,p)$, $\ol D^2(;2,3,3)$, $\ol D^2(;2,3,4)$, $\ol D^2(;2,3,5)$, $\ol D^2(2;p)$ or $\ol D^2(3;2)$ with $p\geq2$. We have $x\in\D_{refl}O$ if and only if $\D_{refl}(S^2/\Ga_x)\neq\emptyset$. The cone points and corner vertices of $S^2/\Ga_x$ correspond to singular edges emanating from $x$.

\subsubsection{3-orbifolds with nonnegative sectional curvature}

Furthermore, we recall the classification of noncompact $3$-orbifolds (without boundary) admitting complete Riemannian metrics of nonnegative sectional curvature. It follows from an orbifold version of the Soul Theorem which states that a complete noncompact Riemannian orbifold with sectional curvature $\ge 0$ contains a totally convex and totally geodesic closed suborbifold, a so-called {\em soul}, and the orbifold is diffeomorphic to the normal bundle of the soul.

\begin{prop}
\label{prop:nonnegsec3orbi}
Every connected complete noncompact Riemannian $3$-orbifold with sectional curvature $\ge 0$ is diffeomorphic to one of the following:
\begin{enumerate}
\item Quotients $\mathbb{R}^3/\Gamma$ for finite subgroups $\Gamma\subset\textrm{O}(3)$. (The soul is a point.)
\item Bundles over $S^1$ with fiber $\mathbb{R}^2/\Gamma'$ for finite subgroups $\Gamma'\subset\textrm{O}(2)$. (The soul is a circle.)
\item A $3$-orbifold arising from $D^2(p)\times[-1,1]$, $p\geq1$, by gluing each of the boundary components $D^2(p)\times\{\pm1\}$ either to itself via a half-rotation or reflection or by making it a reflector boundary component. There are six such orbifolds. (The soul is a mirrored interval.)
\item A $3$-orbifold arising from $V^2(p)\times[-1,1]$, $p\geq1$, by gluing each of the boundary components $V^2(p)\times\{\pm1\}$ either to itself via the reflection at its bisector or by making it a reflector boundary component. There are three such orbifolds. (The soul is a mirrored interval.)
\item Products $\Si^2 \times \R$ for closed 2-orbifolds $\Si^2$ with Euler characteristic $\ge 0$. (The soul is 2-dimensional.)
\item A $3$-orbifold arising from $\Sigma^2 \times [0,1)$, where $\Sigma^2$ is a closed 2-orbifold with Euler characteristic $\ge 0$,
by gluing $\Sigma^2 \times \{0\}$ to itself by a non-trivial involution. (The soul is $2$-dimensional.)
\end{enumerate}
\end{prop}

We observe for future reference
that all orbifolds occuring in the proposition
are 3-discal if and only if the soul is a point,
and solid toric if and only if the soul is 1-dimensional.

Moreover, we note an alternative constructions
for noncompact complete 3-orbifolds with $\sec \ge 0$
and soul a mirrored interval:
Such orbifolds can also be obtained
by starting with two quotients of the 3-ball,
each with boundary $\R P^2(p)$, $S^2(2,2,p)$ or $\bar{D}^2(p)$
for some fixed $p \ge 1$,
and glueing them together
along a closed pointed disc $D^2(p)$
contained in both boundaries.
The (interior of the) resulting 3-orbifolds
are then diffeomorphic to those arising from
$D^2(p) \times [-1,1]$ as in {\em 3.)}.

Similarly, we could start with two quotients of the 3-ball,
each with boundary $\bar{D}^2(;2,2,p)$ or $\bar{D}^2(2;p)$
for some fixed $p \ge 1$,
and then identify two sectors $V^2(p)$.
In this way, we obtain precisely the orbifolds
arising from $V^2(p) \times [-1,1]$
as in {\em 4.)}.

\subsection{Fibrations and decompositions}

\subsubsection{Fibered 3-orbifolds}

An {\em orbifold Seifert fibration} is an orbifold fibration $p:O^3\to B^2$ 
with 3-dimensional total space $O$, 
2-dimensional base $B$ 
and 1-dimensional closed connected generic fiber $F$,
i.e.\ $F$ is the circle $S^1$ or the mirrored interval $\ol I$. 
A {\em Seifert orbifold} is a 3-orbifold admitting a Seifert fibration. 

Every fiber has a neighborhood 
which is fiber preserving diffeomorphic 
to a solid toric orbifold equipped with a canonical Seifert fibration. 
More precisely, 
suppose that $x\in B$ is a point in the base
and let $\Ga_x$ be its local group. 
Then the fiber $p^{-1}(x)$ has a saturated neighborhood 
of the form $(D^2\times F)/\Gamma_x$ 
with the natural fibration $(D^2\times F)/\Gamma_x\to D^2/\Gamma_x$. 
The action $\Gamma_x \acts D^2$ is effective, 
whereas the action $\Gamma_x \acts F$ is in general not. 
Fibers in the boundary have similar model neighborhoods.
A classification of Seifert orbifolds,
locally and globally,
has been given in \cite{BS_seif}. 

Seifert fibrations of solid toric 3-orbifolds 
as well as 1-dimensional fibrations of their toric boudaries 
are in general not unique. 
The next result describes 
which fibrations of the boundary extend to Seifert fibrations. 
Let $V\cong(D^2\times S^1)/\Ga$ 
be a solid toric 3-orbifold.
We call a 1-dimensional fibration of $\D V$ 
{\em horizontal}
if it is isotopic to the fibration 
$\D V\to S^1/\Ga$. 
\begin{lem}
\label{fibrationsoftori}
A 1-dimensional fibration 
of the boundary $\D V$ of a solid toric orbifold $V$ 
extends to a Seifert fibration of $V$ 
if and only if it is not horizontal. 
\end{lem}
\proof
This is a consequence of the fact 
that 1-dimensional fibrations of closed flat 2-orbifolds 
can be isotoped to be geodesic. 
\qed

All compact Seifert 3-orbifolds without bad 2-suborbifolds are geometric, 
see \cite[ch.\ 3]{Thurston_book} and \cite[2.4]{BMP}. 
More precisely, 
a connected closed Seifert orbifold without bad 2-suborbifolds
admits a geometric structure 
modelled on a unique Thurston geometry 
different from hyperbolic and solvgeometry
(i.e.\ on $S^2 \times \R^1$, $S^3$, 
$\R^3$, $Nil$, $\mathbb{H}^2 \times \R^1$, or
$\widetilde{PSL(2,\mathbb{R})}$). 
A solid toric 3-orbifold admits, 
depending on its topological type,
geometric structures modelled on 
some or all of the six contractible model geometries.

A non-solid toric connected compact Seifert orbifold with nonempty boundary 
contains no bad 2-suborbifolds and 
admits an $\H^2\times\R$- or $\R^3$-structure.
(If it admits an $\R^3$-structure then also an $\H^2\times\R$-structure.) 
In fact, 
it can be geometrized in a stronger sense;
namely, it admits a Riemannian metric with totally geodesic boundary 
locally modelled on either $\H^2\times\R$ or $\R^3$.

\medskip

A {\em toric fibration} of a 3-orbifold 
is an orbifold fibration whose generic fiber is a toric 2-orbifold.
(Fibrations with 2-dimensional fibers of other topological types 
will play no role in this text.) 

Connected 3-orbifolds admitting toric fibrations are geometric
with one of the three model geometries 
$\widetilde{PSL(2,\R)}$, $Nil$ or $\R^3$. 
Those with nonempty boundary 
admit euclidean metrics with totally geodesic boundary 
and complete $\R^3$-structures on their interior. 
They are diffeomorphic to 
$T\times[-1,1]$ or $(T\times[-1,1])/\Z_2$
with a toric 2-orbifold $T$ and, in the latter case,
with $\Z_2$ acting by a reflection on $[-1,1]$. 
Unlike in the manifold case,
they are not always Seifert. 
This is due to the fact that, 
whereas a 2-torus or Klein bottle admits 
(infinitely many, respectively, two) circle fibrations,
not all toric 2-orbifolds admit 1-dimensional orbifold fibrations.
(Compare the discussion in \cite{Dunbar}.)

\subsubsection{2-suborbifolds in 3-orbifolds}

A 3-orbifold $O$ is called {\em irreducible} 
if it does not contain any bad 2-suborbifold 
and if every two-sided spherical 2-suborbifold 
bounds a discal 3-suborbifold. 

It is called {\em (topologically) atoroidal}
if every incompressible two-sided toric 2-suborbifold $\Si\subset O$ 
is {\em boundary parallel}, 
i.e.\ bounds a collar neighborhood $\cong\Si\times[0,1]$ 
of a boundary component $\cong\Si$. 

Let $O$ be a 3-orbifold
and let $\Si\subset O$ be a proper 
$2$-suborbifold.
A {\em compressing discal 2-suborbifold} or {\em compression disc} for $\Si$ 
is a discal $2$-suborbifold $D\subset O$ 
which intersects $\Si$ transversally in $\partial D = D \cap \Si$ 
such that $\D D$ does not bound a discal 2-suborbifold in $\Si$.
(If $D$ is one-sided, 
we understand this to mean that 
splitting the connected component of $\Si$ containing $\D D$ along $\D D$
does not yield a discal 2-orbifold.
Anyway, 
one-sided compression discs can be replaced by two-sided ones
by passing to the boundary of a tubular neighborhood.)
Note that a spherical 2-suborbifold has no compression discs
because every closed 1-suborbifold of a spherical 2-orbifold 
bounds a discal 2-suborbifold. 

A {\em compression} of $\Si$ is either 
a discal 3-suborbifold whose boundary is a component of $\Si$ 
or a compression disc for $\Si$.
If $\Si$ admits a compression 
then it is called {\em compressible},
and otherwise {\em incompressible}. 

Thus a $3$-orbifold is irreducible
if it contains no bad $2$-suborbifolds 
and if all two-sided spherical 2-suborbifolds are compressible. 

The notion of incompressibility 
is particularly useful in the irreducible case
because then the position e.g.\ of closed 2-suborbifolds $\Si$
relative to incompressible 2-suborbifolds $\Si_{inc}$
can be simplified by isotopies.
Namely, it can be acheived 
that $\Si_{inc}$ divides $\Si$ into non-discal components. 

Discal 3-orbifolds are irreducible. 
This is formulated but not proved in \cite[Thm.\ 3.1]{BMP}.
A proof can be found in \cite[2.4]{act3mfs}.

More generally,
every closed 2-suborbifold of a discal 
3-orbifold 
is compressible.
For nonspherical suborbifolds this follows from the 
Equivariant Loop Theorem \cite{MeeksYau}, 
cf.\ \cite[Thm.\ 3.6]{BMP}.
We will use it only for toric 2-suborbifolds. 

\subsubsection{Decompositions of 3-orbifolds along 2-suborbifolds}
\label{sec:decos}

We suppose in the following that $O$ is a compact 3-orbifold. 

Let ${\mathcal F}$ be a finite family 
of disjoint two-sided closed 2-suborbifolds $\Si_j\subset int(O)$. 
The operation of removing from $O$ 
an open tubular neighborhood of $\cup_j\Si_j$ 
is called {\em splitting} $O$ 
along the $\Si_j$. 
We call 
the splitting spherical (toric, incompressible) 
if all $\Si_j$ are spherical (toric, incompressible). 
We will refer to the $\Si_j$ as {\em splitting 2-suborbifolds}. 

A {\em connected sum decomposition} of $O$ 
or a {\em surgery} on $O$ 
is performed by 
first splitting $O$ along a family of spherical 2-suborbifolds 
and then filling discal 3-orbifolds 
into the additional spherical boundary components 
created by the splitting. 
Conversely, 
$O$ is called a {\em connected sum}
of the 3-orbifolds resulting from this decomposition. 
Note that we allow connected sums of connected orbifolds (components) 
with themselves. 

The following result reduces the study 
of compact 3-orbifolds without bad 2-suborbifolds 
to the study of irreducible ones. 
It is due to Kneser \cite{Kneser} in the manifold case, 
see \cite[3.3]{BMP} for a proof in the case of orientable orbifolds.
The argument given there also extends to the nonorientable case.
\begin{thm}[Spherical decomposition]
A compact $3$-orbifold without bad $2$-suborbifolds 
can be decomposed by surgery 
into finitely many irreducible compact 3-orbifolds.
\end{thm}

Let us now consider toric splittings.
\begin{lem}
\label{lem:irrincpiec}
Suppose that $O$ is split along a toric family ${\mathcal T}$ 
into compact pieces $O_i$. 
Then $O$ is irreducible and ${\mathcal T}$ is incompressible (in $O$) 
if and only if all pieces $O_i$ are irreducible 
and for each piece $O_i$ 
the portion $\D O_i-\D O$ of its boundary corresponding to ${\mathcal T}$ 
is incompressible (in $O_i$). 
Moreover,
if in this situation all boundary components of the $O_i$ are incompressible,
then $O$ has incompressible boundary. 
\end{lem}
\proof
The standard proof in the manifold case carries over.
The ``only if'' direction uses the fact 
that toric 2-suborbifolds of discal 3-orbifolds are always compressible.
\qed

There is a canonical splitting of irreducible compact 3-orbifolds 
along incompressible toric suborbifolds. 
It is due to Jaco, Shalen and Johannson in the manifold case
and has be extended to orbifolds by Bonahon and Siebenmann \cite{BS_split},
see also \cite[3.3 and 3.15]{BMP}. 
\begin{thm}[JSJ-splitting]
An irreducible compact 3-orbifold 
admits an incompressible toric splitting 
into components
each of which is atoroidal or Seifert fibered (or both).
A minimal such splitting is unique up to isotopy.
\end{thm}
We will also consider a class of toric splittings
with weaker properties.
Following Waldhausen's definition \cite{Waldhausen} in the manifold case, 
we define a {\em graph splitting} 
of a compact 3-orbifold with toric boundary 
to be a (not necessarily incompressible) toric splitting into pieces 
which admit orbifold fibrations with 1- or 2-dimensional closed fibers.
Moreover, 
the 2-dimensional fibers are required to be toric. 
We will refer to the pieces with 2-dimensional fibrations 
as pieces {\em with toric fibrations}.
A 3-orbifold admitting a graph splitting is called a {\em graph orbifold}. 
Briefly,
it is a 3-orbifold which can be ``cut up into fibered pieces''.

Connected compact 3-orbifolds with toric fibrations 
and nonempty boundaries are diffeomorphic to 
$T\times[-1,1]$ or $(T\times[-1,1])/\Z_2$
with a toric 2-orbifold $T$ and, in the latter case,
with $\Z_2$ acting by a reflection on $[-1,1]$. 
Unlike in the manifold case,
they are not always Seifert. 
This is due to the fact that, 
whereas a 2-torus or Klein bottle admits 
(infinitely many, respectively, two) circle fibrations,
not all toric 2-orbifolds admit 1-dimensional orbifold fibrations,
as already discussed above.
Hence, in the orbifold case
a graph splitting may comprise non-Seifert pieces.

Seifert orbifolds with discal base orbifold are solid toric.
All other connected Seifert orbifolds with nonempty boundary 
have base orbifolds of Euler characteristic $\chi\leq0$ 
and admit nonpositively curved Riemannian metrics 
with totally geodesic boundary. 
(These metrics can be modelled on $\H^2\times\R$ or $\R^3$,
depending on whether $\chi<0$ or $\chi=0$.) 
The 3-orbifolds 
with toric fibrations and nonempty boundary 
admit flat metrics with totally geodesic boundary. 
The existence of these geometric structures 
on the non-solid toric pieces of a graph splitting 
implies that 
they are irreducible and have incompressible boundaries.\footnote{
For the irreducibility one uses the fact 
that discal 3-orbifolds are irreducible.
Namely, 
consider an embedded 2-sphere $S\subset\tilde P$ 
in the universal cover of such a piece $P$
preserved by a finite group $\Ga_S$ of isometries. 
With the Hadamard-Cartan Theorem it follows that $S$ is contained
in a $\Ga_S$-invariant closed ball on which the $\Ga_S$-action is standard.
Hence the action on the ball bounded by $S$ is also standard.}
Moreover,
the pieces 
with toric fibrations and nonempty boundaries 
are atoroidal. 

If in a non-trivial graph splitting of $O$ occur no solid toric pieces,
then all pieces are irreducible atoroidal or Seifert orbifolds 
with nonempty incompressible toric boundaries.
Hence $O$ is irreducible with incompressible toric boundary 
and the splitting is incompressible,
compare Lemma~\ref{lem:irrincpiec}. 
In particular,
a minimal incompressible graph splitting 
of an irreducible compact connected 3-orbifold 
with incompressible toric boundary 
is {\em canonical} up to isotopy
because it coincides with the JSJ-splitting,
unless the orbifold admits a toric fibration over a closed 1-orbifold
(in which case it is geometric). 
Indeed, 
suppose that a nontrivial minimal incompressible graph splitting 
were not minimal as a splitting into atoroidal and Seifert components.
Then for some splitting toric 2-suborbifold $T$
the union of the (one or two) components adjacent to it 
cannot be Seifert and must therefore be atoroidal.
The definition of atoroidality then implies 
that one of these components is $\cong T\times[0,1]$,
contradicting the minimality of the graph splitting. 

A {\em geometric splitting} 
of an irreducible compact connected 3-orbifold $O$ 
is an incompressible toric splitting 
into geometric pieces,
i.e.\ into irreducible compact 3-orbifolds 
whose interiors admit complete geometric structures. 
We refer to the components of the splitting as {\em geometric pieces}. 
Note that if $O$ itself is not geometric,
then the pieces have nonempty boundaries 
and admit geometric structures modelled on 
$\H^3$, $\H^2\times\R$ or $\R^3$.
In particular,
a nontrivial incompressible graph splitting of $O$ 
is a geometric splitting into pieces 
admitting $\H^2\times\R$- or $\R^3$-structures. 

A compact 3-orbifold 
is said to be 
{\em decomposable into geometric pieces}
or to {\em satisfy Thurston's Geometrization Conjecture}
if it can be decomposed by surgery 
into irreducible compact connected 3-orbifolds
which are geometric or admit a geometric splitting, 
cf.\ \cite[3.7]{BMP}.


\subsection{From graph splittings to geometric decompositions}

Graph splittings of compact 3-orbifolds have fairly weak properties 
and are in particular far from being unique.
In this section we show that 
a graph splitting can be improved to a geometric decomposition. 
The manifold case of this discussion
is due to Waldhausen \cite{Waldhausen}.

\begin{thm}
\label{thm:grgeo}
Suppose that $O$ is a compact connected 3-orbifold with toric boundary 
which admits a graph splitting.

If $O$ is irreducible,
then it is either solid toric 
or has incompressible toric boundary.
In the latter case,
it is geometric 
with model geometry different from $S^2\times\R$ and $H^3$,
or it admits an incompressible graph splitting
(and hence a JSJ-splitting without hyperbolic components).

If $O$ is not irreducible,
it can be decomposed by surgery into irreducible orbifolds of this kind. 

%
\end{thm}
\proof
Consider a graph splitting of $O$ 
along a toric family ${\mathcal T}$. 

If ${\mathcal T}=\emptyset$, 
then $O$ is solid toric,
or it is closed and admits an $S^2\times\R$-structure,
or it is irreducible with (possibly empty) incompressible boundary 
and admits a geometric structure modelled on 
one of the six geometries different from $S^2\times\R$ and $H^3$.
If $O$ admits an $S^2\times\R$-structure,
it can be decomposed by surgery along a spherical cross section 
into one or two spherical 3-orbifolds. 

If ${\mathcal T}\neq\emptyset$,
then the pieces of the splitting have nonempty boundary. 
If one of the pieces $P$ is of the form $T\times[0,1]$
with a toric 2-orbifold $T$,
we may reduce ${\mathcal T}$ by erasing one of the components of $\D P$, 
unless $P$ is the only piece.
In the latter case, $O$ fibers over the circle and is geometric
(with model geometry $\R^3$, $Nil$ or $Solv$). 

If no solid toric piece occurs in the graph splitting
(and if ${\mathcal T}\neq\emptyset$), 
then $O$ is irreducible with incompressible boundary 
and the graph splitting is incompressible,
cf.\ Lemma \ref{lem:irrincpiec}.

If there is a solid toric piece
and if adjacent to it there is another solid toric piece 
or a one-ended piece with toric fibration,
i.e.\ a piece diffeomorphic to $(T\times[-1,1])/\Z_2$ with $T$ toric 
and $\Z_2$ reflecting on $[-1,1]$,
then $O$ is closed and geometric 
with model geometry $S^3$ or $S^2\times\R$. 

Since in all other cases we are done or can reduce the splitting 
by removing a component from ${\mathcal T}$,
we assume that at least one solid toric piece 
$V_0\cong(D^2\times S^1)/\Ga$ occurs in the graph splitting 
and that adjacent to $V_0$ there is a non-solid toric Seifert piece $S$. 
We may further assume that all pieces with toric fibrations are one-ended,
i.e.\ diffeomorphic to $(T\times[-1,1])/\Z_2$ with $T$ toric 
and $\Z_2$ reflecting on $[-1,1]$.

We denote by $q:V_0\cong(D^2\times S^1)/\Ga\to S^1/\Ga$ 
the fibration of $V_0$ by discal cross sections.
Let $K\triangleleft\Ga$ be the kernel of the action $\Ga\acts S^1$. 
Via the action $K\acts D^2$ 
we may regard $K$ as a subgroup $K\subset O(2)$. 
The generic discal cross section of $V_0$ is $\cong D^2/K$. 
Let $p:S\to B$ denote the Seifert fibration of $S$. 
The base $B$ is a non-discal 2-orbifold with nonempty boundary.
Let $T_0=\D V_0=S\cap V_0\in{\mathcal T}$ 
denote the toric 2-suborbifold separating $V_0$ and $S$,
and $\D_0 B=p(T_0)$ the boundary component of $B$ 
corresponding to $T_0$. 
It is either a circle or an arc connecting two points of 
$\D_{refl}B\cap\D B$. 

If the fibrations $p|_{T_0}$ and $q|_{T_0}$ 
of $T_0$ by closed 1-orbifolds 
are not isotopic,
then the Seifert fibration $p$ can be extended over $V_0$,
compare Lemma~\ref{fibrationsoftori},
i.e.\ $S\cup V_0$ is Seifert 
and we reduce the graph splitting 
by removing the component $T_0$ from ${\mathcal T}$. 
(The Seifert piece ``swallows'' the adjacent solid toric piece.)

Otherwise, 
if $p|_{T_0}$ and $q|_{T_0}$ are isotopic,
we may assume that they agree, 
i.e.\ that the discal cross sections of $V_0$ fill in Seifert fibers. 
We then have the identification $\D_0 B\cong S^1/\Ga$.
In this situation we find the following class 
of two-sided spherical 2-suborbifolds 
adapted to the graph structure. 
Let $\al\subset B-B^{sing}$ be a properly embedded arc 
with endpoints in $\D_0B-\D_{refl}B$. 
It yields the spherical 2-suborbifold 
$\Si_{\al}\subset S\cup V_0$ 
obtained by taking $p^{-1}(\al)\subset S$ 
and attaching to it the pair of discal 2-suborbifolds 
$q^{-1}(\D\al)\subset V_0$. 
Hence $\Si_{\al}\cong S^2/K$ 
where we extend the action of $K\subset O(2)$ to $\R^3$ 
using the canonical embedding $O(2)\subset O(3)$.
If $\D_{refl}B\neq\emptyset$
we get another similar class of spherical 2-suborbifolds 
by taking embedded arcs $\hat\al$ 
connecting a regular boundary point on $\D_0 B$ 
to an interior point of a reflector edge. 
In this case we have $\Si_{\hat\al}\cong S^2/\hat K$
where $\hat K\subset O(2)$ is an index two extension 
of $K$ such that the elements in $\hat K-K$ 
switch the poles $(0,0,\pm1)$ of $S^2$. 

We note that the discal 3-orbifold $D^3/K$ has a decomposition 
into a Seifert and a solid toric piece 
analogous to the decomposition of $S\cup V_0$. 
Indeed, take an $O(2)$-invariant decomposition of $D^3$ 
into a tubular neighborhood of the equator $S^1\times\{0\}$ 
and the complement of this neighborhood. 
By dividing out $K$, one sees that $D^3/K$ is obtained 
by attaching the cylinder $(D^2/K)\times[-1,1]$ 
to a product fibration with fiber $S^1/K$ 
and base 
a bigon (topological disc) such that the fibrations of the boundaries match. 
Thus we may split $O$ along $\Si_{\al}$,
fill in copies of $D^3/K$ 
into the two spherical boundary components resulting from the splitting
to obtain a (possibly disconnected) 3-orbifold $O_{\al}$,
and extend the graph splitting to all of $O_{\al}$ 
by attaching copies of $(\D D^2/K)\times[-1,1]$ 
to the pieces of $T_0$. 
The effect on the base of the Seifert piece is that 
$B$ is split along $\al$
and two bigons are attached along the copies of $\al$. 
The discal 3-orbifold $D^3/\hat K$ is obtained analogously 
by attaching the cylinder quotient $(D^2\times[-1,1])/\hat K$, 
where $\hat K/K$ acts on $[-1,1]$ by a reflection, 
to a Seifert manifold with base orbifold a triangular disc $\De$ 
whose boundary consists of two boundary arcs and one reflector boundary arc. 
When performing the connected sum decomposition (surgery) 
of $O$ along $\Si_{\hat\al}$
and extending the graph splitting over $O_{\hat\al}$, 
the effect on the base of the Seifert piece is 
that $B$ is again split along $\hat\al$, 
but this time copies of $\De$ are attached along the copies of $\hat\al$. 

If ${\mathcal A}$ is a finite system of disjoint properly embedded arcs
$\al_i$ and $\hat\al_j$ in $B$ as above,
we denote the result of performing simultaneous surgeries
along the system of spherical 2-suborbifolds 
$\Si_{\al_i}$ and $\Si_{\hat\al_j}$
by $O_{\mathcal A}$ 
and equip it with an induced graph splitting 
along a toric family ${\mathcal T}_{\mathcal A}$
as explained. 
The components of ${\mathcal T}$ different from $T_0$ 
correspond to components of ${\mathcal T}_{\mathcal A}$, 
whereas $T_0$ may split up into several components. 
Furthermore, 
we denote by $S_{\mathcal A}\subset O_{\mathcal A}$ 
the Seifert suborbifold corresponding to $S$ 
and by $B_{\mathcal A}$ the base of its Seifert fibration. 

Now we choose the system of arcs ${\mathcal A}$ 
so that they split $B$ into pieces as simple as possible. 
After making a suitable choice, 
every connected component $B'$ of $B_{\mathcal A}$
is diffeomorphic to one of the compact 2-orbifolds 
in the following list: 
\BI
\item $Ann$ 
or the quadrangle $Q$ bounded by two boundary and two reflector edges 
occuring in alternating order
(a quotient of $Ann$ by an involution);
\item $\ol D^2(p)$ or $\ol V^2(p)$ with $p\geq1$. 
\EI
\begin{lem}
\label{lem:simplebase}
Suppose that the compact 3-orbifold $O'$ has a decomposition 
$O'=S'\cup V'$ along a toric 2-suborbifold $T'=S'\cap V'$
into a Seifert piece $S'$ and a solid toric piece $V'$ 
such that the discal cross sections of $V'$ fill in Seifert fibers of $S'$. 
If $B'$ belongs to the above list,
then $O'$ is either spherical or solid toric. 
\end{lem}
\proof
If $B'$ splits as the product of the compact interval 
and a connected closed 1-orbifold,
i.e.\ if $B'$ is the annulus or the quadrangle, 
then $O'=S'\cup V'\cong V'$ is solid toric. 

If $B'$ is discal,
then $S'$ is also solid toric
and has the form 
$S'\cong(D^2\times F')/\Ga'$
with faithful action $\Ga'\acts D^2$ 
and generic Seifert fiber $F'$.
Let $\De'$ denote the discal cross section of $V'$.
Using an identification $\D\De'\cong F'$,
we form the closed 3-orbifold $\hat O'$ 
by gluing $D^2\times F'$ and $\D D^2\times\De'$  
canonically along their boundaries. 
We extend the $\Ga'$-action from $D^2\times F'$ to $\hat O'$
by choosing an extension of the $\Ga'$-action on $F'$ to $\De'$. 
There exists a $\Ga'$-invariant spherical structure on $\hat O'$.
Hence $O'\cong\hat O'/\Ga'$ is spherical. 
\qed

\medskip\no
{\em Proof of Theorem~\ref{thm:grgeo} continued.}
As a consequence of Lemma~\ref{lem:simplebase},
the splitting of $O_{\mathcal A}$
along the subfamily of ${\mathcal T}_{\mathcal A}$
consisting of those toric 2-suborbifolds, 
which correspond to the components of ${\mathcal T}$ different from $T_0$, 
is still a graph splitting. 
(The pieces resulting from splitting $V_0$ swallow
the corresponding adjacent pieces of $S$.) 

Our discussion yields so far:
$O$ is irreducible and satisfies the conclusion of the theorem,
or $O$ is closed and admits an $S^2\times\R$-structure,
or $O$ can be decomposed by surgery into components
which admit graph splittings along strictly fewer toric 2-suborbifolds.

By repeating this process finitely many times,
it follows that $O$ can be decomposed by surgery into irreducible orbifolds 
satisfying the conclusion of the theorem. 
In particular, 
the assertion holds if $O$ is not irreducible.
If $O$ is irreducible,
then $O$ is diffeomorphic to one of the components
arising from the surgery,
and the assertion holds as well.
\qed

\begin{cor}
\label{cor:graphgeom}
A compact connected 3-orbifold with toric boundary 
which admits a graph splitting
is decomposable into geometric pieces
(and hence satisfies Thurston's Geometrization Conjecture).
\end{cor}

\begin{rem}
Let $O$ be as in Theorem~\ref{thm:grgeo}. 
Then the boundary of $O$ is incompressible 
if and only if 
no solid toric components occur in the surgery decomposition of $O$.
Indeed, 
suppose that $V$ is a solid toric component. 
Then there exists a finite family 
of disjoint embedded discal 3-suborbifolds $B_i\subset V$
such that $V-\cup_iB_i$ embeds into $O$.
There exists a compression disc for $\D V$ in $V$ avoiding the $B_i$, 
and hence a compression disc for $\D V\subset\D O$ in $O$.
Conversely,
suppose that $\D O$ is compressible and consider a compression disc $\De$. 
Using the property that $O$ contains no bad 2-suborbifolds,
we can make $\De$ step by step disjoint 
from the family of spherical 2-suborbifolds 
along which the surgery is performed 
until $\De$ is contained in an irreducible component. 
This component must be solid toric. 
\end{rem}

\section{Coarse stratification of roughly $\leq2$-dimensional 
Alexandrov spaces}


\subsection{Preliminaries}
For background on spaces with curvature bounded below
we refer to the basic text \cite{BGP}
and to the material in \cite[ch.\ 10]{BBI}. 

By a {\em segment} we mean more precisely a 
distance minimizing geodesic segment.
Given two points $x$ and $y$, then $xy$ denotes one of the possibly several 
segments connecting these points. 

\subsubsection{Alexandrov balls}
All arguments from Alexandrov geometry used in this paper will be local 
and accordingly work in an appropriate class of local Alexandrov spaces. 
\begin{defn}[Alexandrov ball]
\label{def:alexb} 
An {\em Alexandrov ball} of curvature $\geq\kappa$ 
is a local Alexandrov space with curvature $\geq\kappa$ 
of the form $X=B(x,\rho)$, $\rho>0$,  
with the additional properties that 
the closed balls $\ol B(x,r)$ for $r\in(0,\rho)$ are metrically complete,  
and that for any two points $y,z \in X$ with 
$d(y,z) + d(x,y) + d(x,z) < 2\rho$ 
exists a segment $yz$ joining $y$ and $z$.
\end{defn}
The first property can be viewed as metrical completeness
``up to radius $<\rho$''
and the second is a global form of the length space condition. 
We call $x$ the {\em center} of the Alexandrov ball $B(x,\rho)$ 
and the minimal number $r\in[0,\rho]$ such that 
$\ol B(x,r)=B(x,\rho)$ its {\em radius}
(with respect to $x$). 
An Alexandrov space may be regarded 
as an Alexandrov ball with infinite radius. 

For any pair of points in $B(x,\frac{\rho}{2})$
or, more generally, in a ball $B(y,r)$ 
with $d(x,y)+2r\leq\rho$
exists a segment connecting them. 
For any triple of vertices in $B(x,\frac{\rho}{3})$
or, more generally, in a ball $B(y,r)$ 
with $d(x,y)+3r\leq\rho$
geodesic triangles exist 
and they satisfy triangle comparison
due to the version of Toponogov's Theorem for Alexandrov spaces, 
cf.\ \cite[\S 3]{BGP}. 

\subsubsection{Strainers and cross sections}
\label{sec:strainalex}

Let $X=B(x,R)$ be an Alexandrov ball with curvature $\geq-1$. 
All points occuring in our discussion below 
are supposed to lie in $B(x,\frac{R}{3})$. 

For a (small) constant $\theta>0$,
a $\theta$-straight $n$-{\em strainer} of length $l$ ($>l$)
in a point $x\in X$ 
consists of $n$ pairs of points $a_i,b_i$ at distance $l$ ($>l$)
from $x$ 
such that 
$\cangle_p(a_i,b_i)\geq\pi-\theta$,
$\cangle_p(a_i,a_j)\geq\pihalf-\theta$ and 
$\cangle_p(b_i,b_j)\geq\pihalf-\theta$ 
for all $i\neq j$, and 
$\cangle_p(a_i,b_j)\geq\pihalf-\theta$ for all $i,j$.
(All comparison angles are taken in the hyperbolic plane.) 
We call the strainer $<\theta$-straight
if it is $\theta'$-straight for some $\theta'<\theta$. 
Compare the definition of {\em burst points} in \cite[\S 5.2]{BGP}
and the definition of strainers \cite[\S 10.8.2]{BBI}.
We say that a strainer is {\em equilateral}
if all points $a_i,b_i$
have the same distance from $x$.

Similarly,
we define a $\theta$-straight 
$n\half$-{\em strainer} of length $l$ ($>l$) in $x$ 
as such an $n$-strainer together with an additional point $a_{n+1}$
at distance $l$ ($>l$) from $x$
such that 
$\cangle_p(a_{n+1},a_i)\geq\pihalf-\theta$ and 
$\cangle_p(a_{n+1},b_i)\geq\pihalf-\theta$ 
for all $i\leq n$. 
We define an {\em infinitesimal} strainer
as a configuration of directions in $\Si_xX$ 
satisfying analogous inequalities. 

Due to the monotonicity of comparison angles, 
the existence of a strainer of length $l$ at $x$ 
implies the existence of strainers at $x$ 
of the same type and straightness 
with any length $l'<l$. 

For an $n$-strainer 
$(a_1,b_1,\dots,a_n,b_n)$ in $x$ 
we put 
\begin{equation*}
f_i:=f_{a_i,b_i}:=
\half(d(a_i,\cdot)-d(b_i,\cdot)) - \half(d(a_i,x)-d(b_i,x))
\end{equation*}
and $f:=f_{a_1,b_1,\dots,a_n,b_n}:=(f_1,\dots,f_n)$, 
normalized to vanish in $x$. 
The functions $f_i$ are 1-Lipschitz. 
We call the level sets $f^{-1}(t)$ the 
{\em cross sections} of the strainer 
and denote by 
$\Si_{y;a_i,b_i}=f_i^{-1}(f_i(y))$ and 
$\Si_{y;a_1,b_1,\dots,a_n,b_n}=f^{-1}(f(y))$ 
the cross sections through the point $y$. 

The (Hausdorff) {\em dimension} 
can be characterized in terms of strainers as follows. 
There exists a constant $\bar\theta_{d\half}>0$ such that 
$X$ has dimension $>d$ if and only if some point in $X$ 
admits a $\bar\theta_{d\half}$-straight $d\half$-strainer 
of some positive length. 
(This in turn is implied by the existence of 
a $<\bar\theta_{d\half}$-straight infinitesimal $d\half$-strainer
in some point.) 
These constants need not be extremely small, 
e.g.\ $\bar\theta_{1\half}$ may be chosen arbitrarily 
in $(0,\pihalf)$ 
and $\bar\theta_{2\half}$ in $(0,\pitenth)$. 

Thus, if a sequence of $d$-dimensional pointed Alexandrov balls 
$(X_i,x_i)$ 
with curvature $\geq-1$ 
{\em collapses} to a pointed Alexandrov ball $(X_{\infty},x_{\infty})$ 
of dimension $\leq k<d$, 
$(X_i,x_i)\to(X_{\infty},x_{\infty})$, 
then for any radius $r>0$ the supremum of the lengths of 
$\bar\theta_{k\half}$-straight $k\half$-strainers
in points of $\ol B(x_i,r)$ 
tends to zero as $i\to\infty$. 

\subsubsection{Comparing comparison angles}
\label{sec:strcurvbound}

We will need to compare
the comparison angles of (equilateral) 1-strainers
with respect to 
different model spaces of constant curvature
with curvature values in the interval $[-1,0]$.
We observe that for a triangle with fixed side lengths
comparison angles increase monotonically
with the comparison curvature value.
We are interested
in 1-strainers of bounded length.

Consider a triangle $\Delta$ in euclidean space
with two sides of length 2
and angle $\pi - \theta$ between them.
We define the angle $\alpha(\theta)$
by letting $\pi - \alpha(\theta) < \pi - \theta$
be the comparison angle of $\Delta$
with respect to hyperbolic space of curvature $-1$,
i.e.\ the corresponding angle
of a triangle in hyperbolic space
with the same side lengths as $\Delta$.

\begin{lem}
\label{lem:alpha}
The function $\alpha(\theta)$
is differentiable in $\theta = 0$
with $\alpha'(0^+) = (2\coth 2)^\half < \frac{3}{2}$.
\end{lem}

\proof
Let $l = 4-h$ denote the length
of the third side of $\Delta$,
i.e.\ $l = 4 \cos \frac{\theta}{2}$ and $h = \half \theta^2 + O(\theta^3)$.

By the hyperbolic law of cosines, we have
$\cosh l = (\cosh 2)^2 + (\sinh 2)^2 \cos \alpha = \cosh 4 - (\sinh 2)^2 (1-\cos \alpha) = \cosh 4 + \half (\sinh 2)^2 \alpha^2 + O(\alpha^3)$.
Moreover, $\cosh l = \cosh(4-h) = \cosh 4 - \sinh{4}\ h + O(h^2) = \cosh 4 - \half \sinh{4}\ \theta^2 + O(\theta^3)$.
Combining these equations we obtain that $\lim_{\theta \to 0} \alpha^2/\theta^2 = 2 \coth 2$
and the lemma follows.
\qed

We will frequently use the following
application of the lemma.
Let $(X,x)$ be an Alexandrov space of curvature $\ge -1$.
Suppose that $(a,b)$ is a 1-strainer at $x$
of length $\le 2$
with comparison angle $\cangle_x(a,b)\ge \pi - \theta$
with respect to some comparison curvature value $k \in [-1,0]$;
and hence in particular with {\em euclidean} comparison angle
(with respect to $k = 0$) $\ge \pi - \theta$.
Then for sufficiently small $\theta$,
i.e.\ $\theta \in (0,\theta_0)$ for some universal $\theta_0 > 0$,
Lemma \ref{lem:alpha} implies
that the strainer $(a,b)$
has comparison angle $\ge \pi - \frac{3}{2}\theta$
with respect to the comparison curvature value $-1$.

For future reference, we also compute
that $\alpha(\frac{\pi}{2}) < \frac{3\pi}{4}$
by solving the above equations
with $l = 8^\half$
for $\alpha$.  
In other words, if we have an equilateral 1-strainer
of length $\le 2$ which is $\pihalf$-straight
with respect to some comparison curvature value in $[-1,0]$,
it is still $\frac{3\pi}{4}$-straight
with respect to all other comparison curvature values in $[-1,0]$.

\subsection{Uniform local approximation by cones}
\label{sec:approxcone}

Alexandrov spaces can in every point be arbitrarily well locally approximated 
by their tangent cone if one zooms in sufficiently far. 
We need a quantitative version of this infinitesmial conelikeness, 
that is, we need uniform scales on which one can well approximate by cones. 
(Compare the scaling argument in \cite[3.5]{MT2}.)

\begin{defn}[Local approximation by cones]
\label{def:locapp}
We say that the Alexandrov ball $B(y,1)$ is in the point $z\in B(y,\half)$ 
on the scale $s\leq\half$ 
$\mu$-well {\em approximated by a cone}
if the rescaled pointed ball $s^{-1}\cdot(B(z,s),z)$
has Gromov-Hausdorff distance $<\mu$ 
from the euclidean cone of radius 1 
over some Alexandrov space with curvature $\geq1$
and with base point in the tip of the cone. 
\end{defn}
The base of the approximating cone may be empty, 
in which case the cone is just a point. 


The following result says that Alexandrov balls 
of dimension $\leq d$ or, more generally, 
which are {\em roughly $\leq d$-dimensional} in the sense that 
$\bar\theta_{d\half}$-straight $d\half$-strainers 
must be very short, 
can be arbitrarily well locally approximated by cones of dimension $\leq d$. 
\begin{prop}[Local approximation by cones on uniform scales]
\label{prop:locapproxconealex}
For $d\in\N$ and $\si,\mu>0$ exist scales 
$0<s_{d\half}=s_{d\half}(\si,\mu)<<s_1=s_1(d,\si,\mu)<<\si$ 
such that: 

An Alexandrov ball $B(y,1)$ 
with curvature $\geq-1$ 
and without $\bar\theta_{d\half}$-straight $d\half$-strainers
of length $\geq s_{d\half}$
can in every point $z\in B(y,\half)$
on some scale $s(z)\in[s_1,\si]$ 
be $\mu$-well approximated by a cone of dimension $\leq d$. 
\end{prop}
\proof
Consider a sequence of Alexandrov balls $B(y_k,1)$ 
with curvature $\geq-1$ 
and without $\bar\theta_{d\half}$-straight $d\half$-strainers
of length $\geq\frac{\si}{k}$, 
and suppose that there exist points $z_k\in B(y_k,\half)$ 
such that $B(y_k,1)$ can in $z_k$ not be $\mu$-well approximated by a cone 
of dimension $\leq d$ 
on any scale $s\in[\frac{\si}{k},\si]$. 
The $(B(y_k,1),z_k)$ Gromov-Hausdorff converge 
to a pointed Alexandrov ball 
$(B(y_{\infty},1),z_{\infty})$
with curvature $\geq-1$ and dimension $\leq d$, 
and $z_{\infty}\in B(y_{\infty},\half)$. 
Now $B(y_{\infty},1)$ can 
on a sufficiently small scale $s'<\si$ 
be $\mu$-well approximated in $z_{\infty}$ by 
(the truncation at radius 1 of) its tangent cone.  
Since $\frac{\si}{k}<s'$ for large $k$,
we obtain a contradiction. 
\qed


Throughout the paper we will fix some small value for $\si$,
say $\si=\frac{1}{2010}$. 

\subsection{Islands without strainers}
\label{sec:alexisl}

Building on \ref{prop:locapproxconealex} 
we will now divide 
our roughly $\leq d$-dimensional Alexandrov balls $B(y,1)$ 
with curvature $\geq-1$ 
into two regions 
according to the (non)existence of good 1-strainers on a uniform scale. 
We will show that the points without such strainers 
accumulate in ``islands'' which are uniformly separated from each other.  
\begin{defn}[Hump]
\label{def:humpalex}
For small $\theta>0$ 
we call a point $z\in B(y,\half)$
a $(\theta,\mu)$-{\em hump},
if the base of the approximating cone 
provided by \ref{prop:locapproxconealex} 
has diameter $<\pi-\frac{\theta}{2}$. 
\end{defn}
Let $H=H_{\theta,\mu}\subset B(y,\half)$ 
denote the subset of $(\theta,\mu)$-humps,
and let $S=S_{\theta,\mu}\subset B(y,1)$ denote the subset of points 
admitting
$<\theta$-straight 1-strainers of length 
$>\frac{1}{11}s_1(d,\si,\mu)$,
cf.\ \ref{prop:locapproxconealex}. 
We are interested in the distribution of the set $H-S$ 
for small $\theta$ and $\mu$. 
(Our notation suppresses the dependence on $d$.
Later we will only need the case $d=2$.)

If the approximation accuracy $\mu$ is sufficiently small, 
then humps and non-humps have the following properties. 
\begin{lem}
\label{lem:nhump}
For sufficiently small $\theta>0$
(i.e.\ $0 < \theta \le \theta_0$ for some positive $\theta_0$)
there exists $\mu_0(\theta)>0$ such that 
for $\mu\in(0,\mu_0(\theta)]$ holds: 

(i) A $(\theta,\mu)$-hump $z$ admits no  
$\frac{\theta}{4}$-straight 1-strainers of length $\frac{1}{111}s(z)$, 
but 

(ii) all points in the closed annulus 
$\ol A(z;\frac{1}{10}s(z),\frac{9}{10}s(z))$
do admit $\frac{\theta}{11}$-straight 1-strainers 
of length $>\frac{1}{11}s(z)$,
i.e.\ $\ol A(z;\frac{1}{10}s(z),\frac{9}{10}s(z))\subset S_{\theta,\mu}$. 

(iii) Moreover,
if $y'\in B(y,\half)$ is no $(\theta,\mu)$-hump, 
then it admits $<\theta$-straight 1-strainers of length 
$>\frac{99}{100}s(z)>\frac{1}{11}s_1(d,\si,\mu)$,
i.e.\ $B(y,\half)\subset H_{\theta,\mu}\cup S_{\theta,\mu}$. 
\end{lem}
\proof

Property (i) follows from the fact
that euclidean comparison angles
(with respect to comparison curvature value 0)
are larger than their hyperbolic equivalents
(with respect to comparison curvature value -1).

For the proof of (ii), 
every point in $\ol A(z;\frac{1}{10}s(z),\frac{9}{10}s(z))$
can be made arbitrarily close
to the midpoint of a segment 
of length $> \frac{2}{11}s(z)$
if $\mu$ is chosen sufficiently small.
This implies the existence of a 
$<\theta$-straight 1-strainer of length $> \frac{1}{11}s(z)$
at every point in $\ol A(z;\frac{1}{10}s(z),\frac{9}{10}s(z))$
for sufficiently small $\mu$.


To prove (iii), we observe
that for every $\theta >0$ and sufficiently small $\mu >0$
every point $y'$ which is not a $(\theta,\mu)$-hump $z$
admits a 1-strainer $(a,b)$ of length $> \frac{1}{11}s(y')$
which is $\frac{2}{3}\theta$-straight
as a {\em euclidean} 1-strainer,
i.e.\ with respect to
the comparison curvature value 0.
By Lemma \ref{lem:alpha}
and the discussion afterwards
this strainer then is also
$< \theta$-straight
with respect to comparison curvature value $-1$
for sufficiently small $\theta >0$.
\qed

Now we show that humps can be grouped into finitely many islands.
\begin{prop}
\label{prop:finhumpalex}
Let $d\in\N$ and $\si,\mu,\theta>0$ such that $\mu\leq\mu_0(\theta)$. 
Suppose that $B(y,1)$ is an Alexandrov ball 
with curvature $\geq-1$ 
and without $\bar\theta_{d\half}$-straight $d\half$-strainers
of length $\geq s_{d\half}(\si,\mu)$. 
Then there exist finitely many $(\theta,\mu)$-humps $z_j\in H_{\theta,\mu}$
such that 
\begin{equation*}
B(y,\half)\subset\left(\bigcup_jB(z_j,\frac{1}{10}s(z_j))\right)
\cup S_{\theta,\mu} .
\end{equation*}
Moreover,
$d(z_j,z_k)>\frac{9}{10}s(z_j)$ for $j\neq k$. 
\end{prop}
\proof 
Let $z,z'\in H-S$. 
Then 
$z'\not\in\ol A(z;\frac{1}{10}s(z),\frac{9}{10}s(z))$ 
and 
$z\not\in\ol A(z';\frac{1}{10}s(z'),\frac{9}{10}s(z'))$, 
i.e.\ 
$\frac{1}{s(z)}d(z,z'),$ $\frac{1}{s(z')}d(z,z')
\not\in[\frac{1}{10},\frac{9}{10}]$. 
If $z\in B(z',\frac{1}{10}s(z'))$
but $z'\not\in B(z,\frac{1}{10}s(z))$, 
then 
$\frac{9}{10}s(z)<d(z,z')<\frac{1}{10}s(z')$,  
and so
$B(z,\frac{1}{10}s(z))\subset B(z',\frac{1}{9}s(z'))
\subset B(z',\frac{1}{10}s(z'))\cup S$. 
There can be no infinite sequence of points $z,z',z'',\dots\in H-S$
such that 
$s(z)<\frac{1}{9}s(z')<\frac{1}{9^2}s(z'')<\dots$
because the scales take values in the bounded interval $[s_1,\si]$. 
Let $H_1\subset H-S$ denote the subset of all points $z\in H-S$ 
for which no point $z'\in H-S$ exists with 
$\frac{9}{10}s(z)<d(z,z')<\frac{1}{10}s(z')$. 
We note that 
$\cup_{z\in H-S}B(z,\frac{1}{10}s(z))\subset
\cup_{z\in H_1}B(z,\frac{1}{10}s(z))\cup S$. 

By construction, 
the relation on $H_1$ defined by 
$z\sim z':\Leftrightarrow z'\in B(z,\frac{1}{10}s(z))$ 
is reflexive and symmetric,
and we have that 
$z\not\sim z'\Rightarrow d(z,z')>\frac{9}{10}s(z)$ 
and $>\frac{9}{10}s(z')$. 
To verify that the relation is also transitive, 
suppose that $z'\sim z\sim z''$ and $z'\not\sim z''$. 
Then 
$\frac{9}{20}(s(z')+s(z''))
<d(z',z'')\leq d(z,z')+d(z,z'')<\frac{1}{10}(s(z')+s(z''))$,
a contradiction. 
Thus, ``$\sim$'' is an equivalence relation on $H_1$. 
We call an equivalence class a $(\theta,\mu)$-{\em island}. 
Note that the island inhabited by $z$ is contained in the intersection 
$\cap_{z'\sim z}B(z',\frac{1}{10}s(z'))$. 
Since inequivalent points $z',z''\in H_1$ satisfy
$d(z',z'')>\frac{9}{20}(s(z')+s(z''))$, 
islands are separated, 
$B(z',\frac{9}{20}s(z'))\cap B(z'',\frac{9}{20}s(z''))=\emptyset$. 
In particular, 
there are only finitely many islands. 

If $z\sim z'$ and $s(z')<4s(z)$, 
then 
$B(z',\frac{1}{10}s(z'))\subset B(z,\half s(z))$
$\subset B(z,\frac{1}{10}s(z))\cup S$. 
Let $R\subset H_1$ be a subset 
which contains exactly one representative $z$ from each island 
with almost maximal scale value $s(z)$ (among its fellow islanders). 
Then 
$\cup_{z\in H_1}B(z,\frac{1}{10}s(z))\subset
\cup_{z\in R}B(z,\frac{1}{10}s(z))\cup S$. 
$R$ is finite
and for any two distinct points $z',z''\in R$ holds 
$d(z',z'')>\frac{9}{10}s(z')$. 
Altogether we obtain 
\begin{equation*}
H\subset\cup_{z\in H-S}B(z,\frac{1}{10}s(z))\cup S
\subset\cup_{z\in H_1}B(z,\frac{1}{10}s(z))\cup S
\subset
\cup_{z\in R}B(z,\frac{1}{10}s(z))\cup S
\end{equation*}
\qed

\subsection{The 1-strained region}
\label{sec:1stralex}
We will now study the geometry of the region $S=S_{\theta,\mu}$ 
of points admitting good ($<\theta$-straight) 1-strainers
which was introduced in section~\ref{sec:alexisl}.

%

In the estimates provided in this section 
we will abstain from giving explicit constants 
although this could be done in each case.
Instead we use the symbols $c,c',\dots$ to denote generic positive constants, 
i.e.\ constants which constantly change from estimate to estimate 
and which in each estimate (or assertion) 
take some fixed value 
independent of the other parameters $\theta,l,l',r,\la,\dots>0$ involved. 
The estimates hold
for sufficiently small values of the parameter $\theta$,
i.e.\ there exists some $\theta_0>0$ such that they hold for all 
$\theta\in(0,\theta_0]$. 
By decreasing the upper bound $\theta_0$ for $\theta$ as we go along, 
we can also guarantee that the frequently occuring terms of the form $c\theta$ 
are as small as we wish. 
We will always assume that the upper bound $\theta_0$
is sufficiently small such that the conclusions from
Lemma \ref{lem:nhump} hold.

Throughout this section, 
let $X = B(x,10)$ be an Alexandrov ball with curvature $\geq-1$. 

\subsubsection{Local almost product structure}
\label{sec:strest}

Let $(a,b)$ be a $<2\theta$-straight 1-strainer of length $> (1 - \theta)$ at $x$.
We want to apply the following considerations
to $<\theta$-straight 1-strainers on scale $\frac{1}{11} s_1(d,\si,\mu)$
which we rescale to length $\approx 1$.
They are then not necessarily anymore $<\theta$-straight
with respect to comparison curvature value $-1$,
but by Lemma \ref{lem:alpha}
they still are $<2\theta$-straight.
In this way, we avoid
that our constants depend on
the scale $s_1(d,\si,\mu)$.

The estimates given below express 
that near $x$ there is on a certain small scale 
an almost product structure 
with a one-dimensional factor in the direction of the strainer. 


To begin with, 
the points near $x$ admit 1-strainers close to $(a,b)$ 
of comparable quality and almost the same length. 
More precisely, we have 
\begin{equation}
\label{ineq:strqual}
\cangle_{\cdot}(a,b)>\pi-c\theta \qquad\hbox{ on $B(x,\theta)$}
\end{equation}
with a certain constant $c>1$. 
For future reference,
let us denote by $C_0>1$ a constant 
such that (\ref{ineq:strqual}) holds with $c=C_0$. 


To verify (\ref{ineq:strqual}), 
note that the function $d(a,\cdot)+d(b,\cdot)$ along $xx'$ 
has first derivative $<\theta$ in $x$
and second derivative $<c'$ (in a barrier sense) 
along the whole segment. 
Thus 
$d(a,x')+d(x',b)<d(a,x)+d(x,b)+c''\theta^2<d(a,b)+c'''\theta^2$ 
which implies $\cangle_{x'}(a,b)>\pi-c\theta$.

It follows that 
$f_{a,b}+d(b,\cdot)+const
=-(f_{a,b}-d(a,\cdot))+const
=\half(d(a,\cdot)+d(b,\cdot))+const$
is $c\theta$-Lipschitz on $B(x,\theta)$. 

One can define near $x$ a coarse {\em flow} in the strainer direction 
which replaces the gradient flows of $d(a,\cdot)$ and $d(b,\cdot)$. 
For $t\in(-\theta,\theta)$ and $x'\in B(x,3\theta)$ 
let $\Phi^{a,b}_t(x')$ be the intersection point 
$ax'b\cap f_{a,b}^{-1}(f_{a,b}(x')+t)$. 
The resulting maps $\Phi^{a,b}_t:B(x,3\theta)\to X$ 
are well-defined only up to small ambiguity 
because the broken segments $ax'b$ need not be unique. 
Correspondingly, these maps are in general not continuous. 

However,
they are almost distance non-decreasing, 
$d(\Phi^{a,b}_tx_1,\Phi^{a,b}_tx_2)>
(1-c\theta)d(x_1,x_2)-c'\theta|t|$. 
This follows from triangle comparison applied to 
$\De(x_1,x_2,a)$ or $\De(x_1,x_2,b)$ 
and the fact that 
$|t|\leq 
d(x_i,\Phi^{a,b}_tx_i)<(1+c'\theta)|t|$
because $f_{a,b}$ has slope $\approx1$ along $ax_ib$. 
They are also almost inverse to each other, i.e.\  
$\Phi^{a,b}_{-t}\Phi^{a,b}_t$ 
is $c\theta|t|$-close to the identity
(where it is defined). 
To see this, 
consider the triangle 
$\De(x',\Phi^{a,b}_tx',\Phi^{a,b}_{-t}\Phi^{a,b}_tx')$
and note that 
$\cangle_{\Phi^{a,b}_tx'}(x',\Phi^{a,b}_{-t}\Phi^{a,b}_tx')<c''\theta$ 
by (\ref{ineq:strqual})
and again 
$|t|\leq 
d(\Phi^{a,b}_tx',x'),d(\Phi^{a,b}_tx',\Phi^{a,b}_{-t}\Phi^{a,b}_tx')
<(1+c'\theta)|t|$. 
It follows that we have also an upper bound 
$d(\Phi^{a,b}_tx_1,\Phi^{a,b}_tx_2)
\leq(1+c'''\theta)d(x_1,x_2)+c''''\theta|t|$, 
i.e.\ the $\Phi^{a,b}_t|_{B(x,\theta)}$ are 
$(1+c'''\theta t,c''''\theta t)$-quasi-isometric. 

Let $x_1,x_2\in B(x,\theta)$. 
Since $\angle\geq\cangle$ 
and geodesic triangles in $\Si_{x_i}X$ have circumference $\leq2\pi$, 
(\ref{ineq:strqual}) yields 
$\angle_{x_1}(a,x_2)+\angle_{x_1}(b,x_2), 
\angle_{x_2}(a,x_1)+\angle_{x_2}(b,x_1)<\pi+c\theta$.
Since also 
$\cangle_{x_1}(a,x_2)+\cangle_{x_2}(a,x_1),
\cangle_{x_1}(b,x_2)+\cangle_{x_2}(b,x_1)>\pi-c'\theta$,
we obtain that angles and comparison angles with a strainer direction 
almost coincide, 
\begin{equation}
\label{ineq:ang1}
\angle_{x_1}(a,x_2)-\cangle_{x_1}(a,x_2)<c\theta ,
\end{equation}
and 
\begin{equation}
\label{ineq:gegenw}
|\cangle_{x_1}(a,x_2)+\cangle_{x_1}(b,x_2)-\pi|<c\theta ,\qquad
|\cangle_{x_1}(a,x_2)-\cangle_{x_2}(b,x_1)|<c\theta .
\end{equation}
As a consequence, 
$d(a,\cdot)$ is near $x$ {\em almost affine} along segments 
in the sense that its slope is almost constant. 
More precisely, 
\begin{equation}
\label{ineq:distincr}
d(a,x_1)-d(a,x_2)=d(x_1,x_2)\cos\al
\end{equation}
with some angle $\al$ satisfying 
$|\al-\cangle_{x_1}(a,x_2)|<c\theta$, 
as follows from 
(\ref{ineq:ang1},\ref{ineq:gegenw}) 
and the monotonicity of the cosine 
by integrating the derivative of $d(a,\cdot)$ along $x_1x_2$. 
Estimates of the same form hold for $d(b,\cdot)$ 
and $f_{a,b}$.
(Note that (\ref{ineq:distincr}) and the corresponding estimate 
$d(b,x_1)-d(b,x_2)=d(x_1,x_2)\cos\beta$ with $\beta=\pi-\al'$ 
satisfying $|\al'-\cangle_{x_1}(a,x_2)|<c'\theta$
yield that 
$f_{a,b}(x_1)-f_{a,b}(x_2)=d(x_1,x_2)\cos\al''$
with $\cos\al''=\half(\cos\al+\cos\al')$,
i.e.\ $|\al''-\cangle_{x_1}(a,x_2)|<c''\theta$.)
It follows that there exists a constant $L>0$ such that 
the function 
\begin{equation}
\label{eq:slope}
f_{a,b}-\frac{f_{a,b}(x_2)-f_{a,b}(x_1)}{d(x_2,x_1)}d(x_1,\cdot)
\end{equation}
and the analogous functions derived from $d(a,\cdot)$ and $d(b,\cdot)$
are $L\theta$-Lipschitz continuous along the segment $x_1x_2$.  

In particular, 
the {\em cross sections} of a strainer are 
topological hypersurfaces {\em almost perpendicular} to it: 
If $x_1,x_2\in f_{a,b}^{-1}(t)\cap B(x,\theta)$, 
then (\ref{ineq:ang1},\ref{ineq:gegenw}) imply
\begin{equation}
\label{eq:crosssecperp}
\pihalf-c\theta<\cangle_{x_1}(a,x_2)\leq\angle_{x_1}(a,x_2)
<\pihalf+c\theta 
\end{equation}
because the comparison angles of the triangles 
$\De(x_1,x_2,a)$ amd $\De(x_1,x_2,b)$ almost agree. 
Moreover, 
the functions $d(a,\cdot)$, $d(b,\cdot)$ and $f_{a,b}$ 
are $L\theta$-Lipschitz continuous 
on any segment with endpoints 
in the same (piece of) cross section $f_{a,b}^{-1}(t)\cap B(x,\theta)$. 

The maps $\Phi^{a,b}_t$ can be used to compare cross sections, 
since by their definition they satisfy
$\Phi^{a,b}_t(f_{a,b}^{-1}(t')\cap B(x,2\theta))
\subset f_{a,b}^{-1}(t'+t)$. 

Regarding the {\em size} of cross sections 
of $n$-strainers 
$(a_1,b_1,\dots,a_n,b_n)$ 
we obtain from (\ref{ineq:strqual},\ref{eq:crosssecperp}): 
If points 
$x_1,x_2\in f_{a_1,b_1,\dots,a_n,b_n}^{-1}(t)\cap B(x,\theta)$ 
have distance $l$, 
then they admit $c\theta$-straight $n\half$-strainers of length $l$. 

{\em Connectivity} of cross sections. 
Let $x_1,x_2\in f_{a,b}^{-1}(t)\cap B(x,\theta)$. 
For the midpoint $m$ of (a segment) $x_1x_2$ holds 
$|f_{a,b}(m)-t|<c\theta d(x_1,x_2)$. 
As above, by $amb$ with $f_{a,b}^{-1}(t)$
we find an almost midpoint $y\in f_{a,b}^{-1}(t)$ 
at distance $<(\half+c\theta) d(x_1,x_2)$ 
from $x_1$ and $x_2$. 
Iterating this procedure yields a continuous curve 
in $f_{a,b}^{-1}(t)\cap B(x,3\theta)$ connecting $x_1$ and $x_2$.
(Here one uses the earlier estimates with the parameter $3\theta$ 
instead of $\theta$.) 

To simplify notation, 
let us put $\Si^o_{y;a,b}:=\Si_{y;a,b}\cap B(x,\theta)$. 

\begin{lem}[Projecting to cross sections]
\label{lem:projcr}
Let $y,y_1\in B(x,\theta)$ 
and let $z_1$ be the intersection point 
$ay_1b\cap\Si_{y;a,b}$. 
Then 
\begin{equation}
\label{ineq:projstrabcos}
\left|d(y_1,z_1)-d(y,y_1)\cdot|\cos\cangle_y(a,y_1)|\right|<c\theta d(y,y_1) 
\end{equation}
and
\begin{equation}
\label{ineq:projstrabsin}
|d(y,z_1)-d(y,y_1)\sin\cangle_y(a,y_1)|<c'\theta d(y,y_1) .
\end{equation}
In particular,
\begin{equation}
\label{ineq:projstrabpyth}
|d(y,z_1)^2+d(z_1,y_1)^2-d(y,y_1)^2|<c''\theta d(y,y_1)^2 .
\end{equation}
\end{lem}
\proof
%
We put $l=d(y,y_1)$ and $\al_1=\cangle_y(a,y_1)$. 

We have 
$|f_{a,b}(y_1)-f_{a,b}(y)|
\leq d(y_1,z_1)\leq(1+c\theta)|f_{a,b}(y_1)-f_{a,b}(y))|$
and, due to (\ref{ineq:distincr}) and the remark thereafter, 
$f_{a,b}(y_1)-f_{a,b}(y)=-l\cos\al_1'$ with 
$|\al_1'-\al_1|<c\theta$. 
This yields (\ref{ineq:projstrabcos}). 

To estimate $d(y,z_1)$, 
we consider a comparison triangle for 
$\De(y,y_1,z_1)$. 
In view of (\ref{ineq:gegenw}), 
we may exchange $a$ and $b$,
and therefore assume without loss of generality that $z_1\in y_1a$. 
Then $\al_1>\pihalf-c\theta$,
cf.\ (\ref{eq:crosssecperp}). 
Regarding $\cangle_{y_1}(z_1,y)$, 
we have 
$\cangle_{y_1}(a,y)\leq\cangle_{y_1}(z_1,y)
\leq\angle_{y_1}(z_1,y)=\angle_{y_1}(a,y)$
and hence 
\begin{equation*}
|\cangle_{y_1}(z_1,y)-(\pi-\al_1)|<c\theta
\end{equation*}
because of (\ref{ineq:ang1}) 
and $\pi-c\theta<\cangle_{y_1}(a,y)+\al_1\leq\pi$. 
This information implies (\ref{ineq:projstrabsin}). 
(Whether we use a hyperbolic or euclidean comparison triangle 
to compute the length of the side $y_1'z'$ corresponding to $y_1z$, 
causes only a difference by a factor 
$<\frac{\sinh l}{l}<1+\theta^2<1+c\theta$
(due to the distortion of the exponential map 
for hyperbolic plane up to radius $l$)
and we may therefore work with a euclidean one. 
Then 
$|\frac{1}{l}d(y',z_1')-\sin\cangle_{y_1}(z_1,y)|<
|\cos(\pi-\al_1)-\cos\cangle_{y_1}(z_1,y)|+c\theta
<c'\theta$.)

Finally, (\ref{ineq:projstrabpyth}) is a direct consequence. 
\qed 

\subsection{The roughly $\leq2$-dimensional case}
\label{sec:rough2}

We assume now in addition that 
$B(x,10)$ is {\em roughly $\leq2$-dimensional}
in the sense that 
there are no $\bar\theta_{2\half}$-straight 
$2\half$-strainers of length $\la$ for some (very small) $\la>0$. 

\subsubsection{Cross sections of 1-strainers}
We continue our discussion in section~\ref{sec:1stralex}. 
The next two results express that
the cross sections of good 1-strainers
are now roughly $\leq1$-dimensional. 
\begin{lem}
\label{lem:cl21}
Let $(a,b)$ be a 
$C_0\theta$-straight $>(1-\theta^{\half})$-long 1-strainer 
at a point $y\in B(x,5)$. 
Let $y_1y_2$ be a segment of length $l$ 
with endpoints in $B(y,\frac{\theta}{3})$ and with midpoint $m$. 
Let $y'$ be a point with $d(m,y')>\theta l$. 
Suppose that $f_{a,b}$ is $3L\theta$-Lipschitz 
on the segments $y_1y_2$ and $my'\cap B(m,\theta l)$. 
Then $\angle_m(y_i,y')<c\theta$ 
for $i=1$ or 2, 
if $\la=\la(\theta,l)$ is sufficiently small. 

If also $d(y_1,y')\leq\frac{l}{3}$, 
then $\angle_m(y_1,y')<c\theta$.  
\end{lem}
Proof: 
%
%
Let $z'\in my'$ be the point at distance $\theta l$ from $m$. 
Since $f_{y_1,y_2}$ has slope $\equiv1$ on $y_1y_2$, 
there exists a point $z\in y_1y_2\cap\ol B(m,\theta l)$ 
with $f_{y_1,y_2}(z)=f_{y_1,y_2}(z')$. 
It satisfies 
\begin{equation}
\label{formel3}
|f_{a,b}(z)-f_{a,b}(z')|\leq 6L\theta^2l
\end{equation}
because $f_{a,b}$ is $3L\theta$-Lipschitz on $zm$ and $mz'$.

The Lipschitz assumption also implies that 
$|\angle_{\cdot}(y_i,a)-\pihalf|,|\angle_{\cdot}(y_i,b)-\pihalf|
<c\theta$ 
on (the interior of) $y_1y_2$
and hence 
$|\cangle_{\cdot}(y_i,a)-\pihalf|,|\cangle_{\cdot}(y_i,b)-\pihalf|
<c'\theta$
by (\ref{ineq:ang1}). 
Thus the quadrupel $(a,b,y_1,y_2)$ is a $c'\theta$-straight 
2-strainer at $z$. 

The $1\half$-strainer $(y_1,y_2,z')$ at $z$ is $c\theta$-straight 
by (\ref{eq:crosssecperp}).  

We consider now the $1\half$-strainer $(a,b,z')$ at $z$
and estimate $\cangle_z(a,z')$. 
Since $d(z,z')\leq 2\theta l$ 
and $f_{a,b}-d(a,\cdot)$ is $c\theta$-Lipschitz on $B(y,\theta)$, 
(\ref{formel3}) translates to 
$|d(a,z)-d(a,z')|<c'\theta^2l$. 
With (\ref{ineq:distincr}) follows 
$d(z,z')|\cos\cangle_z(a,z')|<c'\theta^2l+c''\theta d(z,z')
<c'''\theta^2l$.

If $d(z,z')>2c'''\bar\theta_{2\half}^{-1}\theta^2l$
(with the constant $c'''$ from the last estimate), 
then $|\cos\cangle_z(a,z')|<\half\bar\theta_{2\half}$
and $|\cangle_z(a,z')-\pihalf|<\bar\theta_{2\half}$. 
Similarly, $|\cangle_z(b,z')-\pihalf|<\bar\theta_{2\half}$,
and it follows that $(a,b,y_1,y_2,z')$ 
is a $<\bar\theta_{2\half}$-straight $2\half$-strainer at $z$
with length $>2c'''\bar\theta_{2\half}^{-1}\theta^2l$.
This is a contradiction 
if $\la=\la(\theta,l)$ is sufficiently small, 
and we conclude that 
$d(z,z')<c''''\theta^2l$ 
and consequently $\cangle_m(z,z')<c\theta$. 

Suppose without loss of generality that $z\in my_1$. 
Then $\angle_m(z,z')-\cangle_m(z,z')\leq 
\angle_m(y_1,z')-\cangle_m(y_1,z')<c'\theta$, 
where the last inequality follows from (\ref{ineq:ang1})
after rescaling by the factor $l^{-1}$. 
Thus $\angle_m(y_1,y')=\angle_m(z,z')<c''\theta$,
which shows the first assertion. 

For the second assertion,
suppose that $d(y_1,y')\leq\frac{l}{3}$ 
but $\angle_m(y_2,y')<c\theta$.
Then $y'$ is close to $my_2$ or 
$y_2$ is close to $my'$. 
Since $d(y_2,y')\geq\frac{2l}{3}$,
only the second alternative can occur and 
$d(m,y')\gtrapprox\frac{l}{2}+\frac{2l}{3}>l$. 
On the other hand 
$d(m,y')\leq d(m,y_1)+d(y_1,y')<l$, 
a contradiction. 
Thus $\angle_m(y_1,y')<c\theta$.
\qed

\begin{lem}[Roughly one-dimensional cross section]
\label{lem:crosssecclsegnew}
Let $l'\leq l\leq\frac{\theta}{50}$, 
$y \in B(x,5)$ 
and let $(a,b)$ be a 
$C_0\theta$-straight $>(1-\theta^{\half})$-long 1-strainer at $y$. 
Suppose that 
$\diam(\Si^o_{y;a,b})>40l$. 

(i) Then there exists a point $y'\in\Si_{y;a,b}\cap B(y,l)$ 
and a segment $y'm$ of length $\geq 10l$, 
such that $f_{a,b}$ is $3L\theta$-Lipschitz on $y'm$ 
and $\Si_{y;a,b}\cap B(y,l)\subset N_{c\theta l}(y'm)$, 
if $\la=\la(\theta,l)$ is sufficiently small. 

(ii) Moreover, 
if $\diam(\Si_{y;a,b}\cap B(y,l))<\frac{199}{100}l$, 
then $y'$ can be chosen in $B(y,\frac{199}{200}l)$ so that 
$\Si_{y;a,b}\cap B(y',r)\subset N_{c\theta r}(y'm)$ 
for all $r\in[l',l]$, 
if $\la=\la(\theta,l')$ is sufficiently small. 
\end{lem}
\proof
(i)
By assumption,
there exists $q\in\Si_{y;a,b}$ with $d(y,q)=20l$.
(Recall from section~\ref{sec:strest} that $\Si_{y;a,b}$ 
is near $y$ path connected.) 
Let $m$ be the midpoint of $yq$. 
Since $f_{a,b}$ is $L\theta$-Lipschitz on $yq$,
we have that $|f_{a,b}(m)-f_{a,b}(y)|<10L\theta l$.
Thus for every point $z\in\Si_{y;a,b}\cap B(y,l)$
the function $f_{a,b}$ has along the segment $mz$ 
Lipschitz constant 
$<\frac{10}{9}L\theta+L\theta=\frac{19}{9}L\theta<3L\theta$, 
cf.\ (\ref{eq:slope}),  
and \ref{lem:cl21} yields that 
$\angle_m(y,z)<c\theta$. 
It follows that the segments $mz$ have pairwise angles $<2c\theta$ 
and are all contained in the $c'\theta l$-neighborhood 
of an almost longest one among them.
We choose $y'$ as its endpoint. 

(ii) By part (i), 
the piece of cross section $\Si_{y;a,b}\cap B(y,l)$ 
is $c\theta l$-close to a subsegment $y'z\subset y'm$ 
where $z\in y'm$ is a point with $d(y,z)=l$ and 
$|d(z,m)+l-d(y,m)|<c'\theta l$. 
Hence $d(y,y')<(\frac{99}{100}+c''\theta)l$. 
The points in $\Si_{y;a,b}\cap B(y,l)$, 
which are further away from $m$ than $y'$, 
must be $2c\theta l$-close to $y'$, 
i.e.\ $d(m,\cdot)|_{\Si_{y;a,b}\cap B(y,l)}$
assumes a maximum in $\Si_{y;a,b}\cap\ol B(y',2c\theta l)$. 
We replace $y'$ by this maximum and then have 
$d(y,y')<\frac{199}{200}l$. 

Also by (i), $\Si_{y;a,b}$ contains
no $\frac{l}{1000}$-long $\pihalf$-straight 1-strainer at $y'$.

Suppose that $\Si_{y;a,b}$ contains 
a $\theta l'$-long $\pihalf$-straight 1-strainer at $y'$. 
Using (i) at the point $y'$ and on the scale $\theta l'$, 
it follows that 
for sufficiently small $\la=\la(\theta,l')$ 
there exists a segment $\tau$ of length almost $2\theta l'$, say, 
of length $\frac{199}{100}\theta l'$ 
along which $f_{a,b}$ is $c\theta$-Lipschitz 
and such that 
$y'$ lies at distance $<c'\theta^2 l'$ from the midpoint of $\tau$. 
(When applying part (i) to points nearby $y$, 
the 1-strainer $(a,b)$ may only be 
$c\theta$-straight for some constant $c>C_0$ 
and $(1-2\theta^{\half})$-long at these points, 
and we use a version of part (i) with appropriate different constants.) 

We consider $d(m,\cdot)$ along the middle third 
$\tau'$ of $\tau$. 
The function $f_{a,b}$ is $c\theta$-Lipschitz along $\tau$ 
and along all segments $zm$ initiating in interior points $z$ of $\tau'$. 
Lemma \ref{lem:cl21} 
implies for sufficiently small $\la=\la(\theta,l')$ 
that these segments $zm$ 
have angles $<c'\theta$ with $\tau'$. 
This means that $d(m,\cdot)$ has slope $\approx\pm1$ along $\tau'$, 
i.e.\ the directional derivatives of $d(m,\cdot)$ 
in directions tangent to $\tau'$ 
take values in 
$[-1,-1+c''\theta^2)\cup(1-c''\theta^2,1]$. 
(See e.g.\ \cite[\S 11]{BGP} 
for a discussion of directional derivatives of distance functions.)
If at some interior point $z_0$ of $\tau'$ 
the directional derivatives 
in the two antipodal directions tangent to $\tau'$ 
are both negative, 
then $d(m,\cdot)$ decays with slope $\approx-1$ 
along both subsegments of $\tau'$ with initial point $z_0$. 
In particular,
$z_0$ is a maximum of $d(m,\cdot)|_{\tau'}$. 
If such a point $z_0$ does not exist,
then $d(m,\cdot)|_{\tau'}$ is almost affine, 
i.e.\ with respect to an appropriate orientation of $\tau'$ 
it increases with slope $\approx1$ along the whole segment. 

Since $y'$ is a maximum of 
$d(m,\cdot)|_{\Si_{y;a,b}\cap B(y',\frac{l}{200})}$,
it follows that $d(m,\cdot)|_{\tau'}$
attains a maximum at a point $y''\in\tau'$ close to the midpoint of $\tau'$, 
more precisely,
at distance $<c''\theta^2 l'$ from $y'$. 
Furthermore, 
there exist at least two segments $\si_1$ and $\si_2$ connecting $y''$ to $m$ 
whose initial directions $\dot\si_i(0)$ 
(with respect to unit speed parametrizations starting at $y''$) 
are close to the two antipodal directions of $\tau'$ at $y''$. 
Each endpoint of $\tau$ lies at distance $<c'''\theta^2 l'$ 
from one of the two segments $\si_1$ and $\si_2$, 
and hence 
$d(\si_1(\theta l'),\si_2(\theta l'))>\frac{19}{10}\theta l'$. 

Applying part (i) at $y''$ on the scales between $\theta l'$ and $l$ 
yields that for $\la(\theta,l')$ sufficiently small 
the continuous function 
$t\mapsto\frac{1}{t}d(\si_1(t),\si_2(t))$
on $[\theta l',l]$
takes values close to 0 and 2, 
i.e.\ in $[0,c\theta)\cup(2-c\theta,2]$. 
However, by the above, 
it has value $\approx0$ for $t=\frac{l}{1000}$
and value $\approx2$ for $t=\theta l'$, 
a contradiction. 
(Note that the smaller the scale, the smaller $\la$ has to be,
and there exists a $\la$ which serves simultaneously 
for all scales $s\in[\theta l',l]$.) 

Thus $\Si_{y;a,b}$ contains 
no $\theta l'$-long $\pihalf$-straight 1-strainer at $y'$. 
It follows, again by part (i) on the scale $\theta l'$, 
that $\diam(\Si_{y;a,b}\cap\D B(y',\theta l'))<c\theta^2 l'$,
and hence every point $z\in\Si_{y;a,b}$
at distance $d\in[\theta l',l]$ from $y'$ 
has distance $<c\theta d$ from $y'm$. 
This proves the assertion. 
\qed

\begin{lem}[Spreading 1-strainers]
\label{lem:exstrainnew}
Let $l,l',y$ and $(a,b)$ 
be as in \ref{lem:crosssecclsegnew}. 
Let $z,z_1,z_2\in\Si_{y;a,b}\cap B(y,l)$ such that 
$(z_1,z_2)$ is a 
$\pihalf$-straight 1-strainer of length $r\in[l',l]$ at $z$. 
Then there exists a constant $C_1\geq C_0$ 
such that for sufficiently small $\la=\la(\theta,l')$ holds:

(i) The 1-strainer $(z_1,z_2)$ at $z$ is $<C_1\theta$-straight. 

(ii) All points 
$u\in azb\cap B(z,\min(100\theta^{-\3quart}r,\theta))$ 
admit $<C_1\theta$-straight $\frac{r}{2}$-long 1-strainers 
contained in $\Si_{u;a,b}$,
and hence such a 2-strainer contained in $\Si_{u;a,b}\cup aub$).
\end{lem}
\proof
(i) This is a consequence of \ref{lem:crosssecclsegnew}(i) 
applied at the point $z$ on the scale $r$. 

(ii) Now we use the $(1+c\theta,c'\theta t)$-quasi-isometry property 
of the maps $\Phi^{a,b}_t$
up to distance $\approx\theta$ from $y$.
We may assume that $u=\Phi^{a,b}_tz$ 
with $|t|<\min(100\theta^{-\3quart}r,\theta)$ 
and put $u_i=\Phi^{a,b}_tz_i$. 
Then $d(u,u_i)<(1+c\theta)r+c'\theta t$
and $d(u_1,u_2)>(1-c\theta)\frac{99}{50}r-c'\theta t$, 
using that $d(z_1,z_2)\approx2r$ 
according to part (i). 
We obtain a $\pithird$-straight 1-strainer $(u_1,u_2)$ at $u$ 
contained in $\Si_{u;a,b}$ and with length $\approx r$. 
Close to the midpoints of the segments $uu_i$ 
we find a 
$\pihalf$-straight 1-strainer $(u'_1,u'_2)$ at $u$ 
contained in $\Si_{u;a,b}$ and with length $\frac{r}{2}$. 
(Compare the argument for the connectivity of cross sections.) 
Applying part (i) again on the scale $\frac{r}{2}$
yields the $c\theta$-straight $\frac{r}{2}$-long 1-strainer, 
once we can make sure that $\diam(\Si^o_{u;a,b})>20r$. 
But this follows from our assumption that 
$\diam(\Si^o_{y;a,b})>40l\geq40r$ 
by using the maps $\Phi^{a,b}_t$ as before. 
\qed

\begin{rem}
\label{rem:edgyonscales}
As in part (i) of the lemma,
we obtain that
if the 1-strainer $(z_1,z_2)$ is $\frac{3\pi}{4}$-straight,
it is also $<C_1 \theta$-straight
and hence in particular $\pihalf$-straight.
\end{rem}

\subsubsection{Edges}
\label{sec:edgalex}

In view of the local product structure, 
the points $y'$ obtained in \ref{lem:crosssecclsegnew}(ii) 
can be considered as points ``near the edge'' of our space. 
They are characterized by the property that 
they admit no long 1-strainers contained in the cross section. 
We will now investigate the geometry near the edge,
keeping the assumption of rough 2-dimensionality 
from section~\ref{sec:rough2}. 

Good 1-strainers at points near the edge 
must be almost perpendicular to the cross section 
if they are not too short:
\begin{lem}[Almost unique 1-strainers near the edge]
\label{lem:1stredgynew}
Let $l,y$ and $(a,b)$ 
be as in \ref{lem:crosssecclsegnew}. 
Suppose that $(a',b')$ is a 
$c\theta^{\half}$-straight
1-strainer of length $\geq l$ at $y$. 

(i) If $\Si_{y;a,b}$ contains no 
$\pihalf$-straight $c'l$-long 1-strainer at $y$, 
then 
$\angle_y(a,a'),\angle_y(a,b'),\angle_y(b,a'),$ $\angle_y(b,b')
\not\in[\pihunth,\frac{99}{100}\pi]$, 
if $\la=\la(\theta,l)$ is sufficiently small. 

(ii) If $\Si_{y;a,b}$ contains no 
$\pihalf$-straight $\theta l$-long 1-strainer at $y$, 
then 
$\angle_y(a,a'),\angle_y(a,b'),\angle_y(b,a'),$ $\angle_y(b,b')
\not\in[c'\theta^{\half},\pi-c'\theta^{\half}]$, 
if $\la=\la(\theta,l)$ is sufficiently small. 

In both cases, $y$ admits no 
$c\theta^{\half}$-straight $l$-long 2-strainer. 
\end{lem}

\proof
There exists a 
$c\theta^{\half}$-straight $l$-long 1-strainer $(y_1,y_2)$ at $y$. 
(Choose $y_1\in ya'$ and $y_2\in yb'$.) 
We put $\al_i=\cangle_y(a,y_i)$. 

By \ref{ineq:distincr} and the remark afterwards, 
$|f_{a,b}(y_i)-f_{a,b}(y)+l\cos\al_i|<c\theta l$.
Thus for $u_1=ay_1b\cap\Si_{y_2;a,b}$ holds 
$(1-c\theta)d(y_1,u_1)<|f_{a,b}(y_1)-f_{a,b}(y_2)|
\leq|f_{a,b}(y_1)-f_{a,b}(y)|+|f_{a,b}(y)-f_{a,b}(y_2)|$
and 
\begin{equation*}
\frac{1}{l}d(y_1,u_1) < |\cos\al_1|+|\cos\al_2|+c\theta ,
\end{equation*}
compare the proof of \ref{lem:projcr}.

Let $z_i=ay_ib\cap\Si_{y;a,b}$. 
According to \ref{lem:projcr}, we have 
$|d(y,z_i)-l\sin\al_i|<c\theta l$. 
Now the rough one-dimensionality of the cross section $\Si_{y;a,b}$, 
i.e.\ \ref{lem:crosssecclsegnew}(i) applied on the scale $l$, 
and our assumption in part (ii) yield that 
$d(y,y')<c'\theta l$ and 
$|d(z_1,z_2)-|d(y,z_1)-d(y,z_2)||<c''\theta l$. 
Hence 
$|\frac{1}{l}d(z_1,z_2)-|\sin\al_1-\sin\al_2||<c'''\theta$. 
With the metric properties of the maps $\Phi^{a,b}_t$ follows 
\begin{equation*}
\left|\frac{1}{l}d(u_1,y_2)-|\sin\al_1-\sin\al_2|\right|<c\theta .
\end{equation*}
Since $(y_1,y_2)$ is $c\theta^{\half}$-straight, 
we have $d(y_1,y_2)>(2-c\theta)l$,  
and (\ref{ineq:projstrabpyth}) implies 
\begin{equation*}
(|\cos\al_1|+|\cos\al_2|)^2+(\sin\al_1-\sin\al_2)^2>4-c'\theta .
\end{equation*}
Writing $\al_i=\pihalf+\beta_i$, 
the last inequality becomes 
$4\sin^2\frac{|\beta_1|+|\beta_2|}{2}=
(\sin|\beta_1|+\sin|\beta_2|)^2+(\cos\beta_1-\cos\beta_2)^2>4-c'\theta$
and 
\begin{equation*}
\sin\frac{|\beta_1|+|\beta_2|}{2}>1-c\theta .
\end{equation*}
Thus $|\beta_i|>\pihalf-c'\theta^{\half}$ and 
$\cangle_y(a,y_1),\cangle_y(a,y_2),\cangle_y(b,y_1),\cangle_y(b,y_2)
\not\in[c'\theta^{\half},\pi-c'\theta^{\half}]$. 

To pass from comparison angles to angles,
we note that if e.g.\ 
$\cangle_y(a,y_1)<c\theta^{\half}$, 
then 
$\angle_y(a,b')=\angle_y(a,y_2)\geq\cangle_y(a,y_2)>\pi-c'\theta^{\half}$. 
This shows (ii).

Under the assumption of part (i) 
we obtain weaker but still useful estimates. 
Now 
$|d(z_1,z_2)-|d(y,z_1)-d(y,z_2)||<(2c'+c''\theta)l<3c'l$ 
and it follows that 
$\left|\frac{1}{l}d(u_1,y_2)-|\sin\al_1-\sin\al_2|\right|<c''$ 
and $|\beta_i|>\pihalf-c'''$. 
The constant $c'''$ depends on the constant $c'$ in the hypothesis of (i) 
and can be made arbitrarily small 
by choosing $c'$ small enough. 
Assertion (i) follows. 
\qed

\begin{rem}
Similarly, one shows e.g.\ that if 
$\Si_{y;a,b}$ contains no 
$\pihalf$-straight $\theta^{\half}l$-long 1-strainer at $y$, 
then 
$\angle_y(a,a'),\angle_y(a,b'),\angle_y(b,a'),$ $\angle_y(b,b')
\not\in[c'\theta^{\quart},\pi-c'\theta^{\quart}]$, 
if $\la=\la(\theta,l)$ is sufficiently small. 
\end{rem}

We make the following choice of scales to quantify rough edges.

\begin{defn}[Edgy points]\label{def:edgy}
A point $y \in B(x,5)$ is called $\theta$-{\em edgy}
relative to a $< 2\theta$-straight
1-strainer $(a,b)$ of length $> (1 - \theta)$ at $y$
if $\diam(\Si^o_{y;a,b}) > \theta^{\frac{5}{2}}$
and $\Si_{y;a,b}$ contains no $\pihalf$-straight
$\theta^4$-long 1-strainer at $y$.

We say that $y \in B(x,5)$ is $\theta$-{\em weakly edgy}
relative to a $< 2\theta$-straight
1-strainer $(a,b)$ of length $> (1 - \theta)$ at $y$
if $\diam(\Si^o_{y;a,b}) > \half \theta^{\frac{5}{2}}$
and $\Si_{y;a,b}$ contains no $\pihalf$-straight
$2\theta^4$-long 1-strainer at $y$.

We call $y \in B(x,5)$ $\theta$-{\em strongly edgy}
relative to a $< 2\theta$-straight
1-strainer $(a,b)$ of length $> (1 - \theta)$ at $y$
if $\diam(\Si^o_{y;a,b}) > \frac{4}{3} \theta^{\frac{5}{2}}$
and $\Si_{y;a,b}$ contains no $\pihalf$-straight
$\frac{3}{4}\theta^4$-long 1-strainer at $y$.
\end{defn}

\begin{lem}[Points without 2-strainers are close to edgy]
\label{lem:closetoedgy}
Suppose that every point $y \in B(x,\theta)$
admits a $<2\theta$-straight
1-strainer $(a_y,b_y)$ of length $> (1 - \theta)$
and that $(a,b)$ is a $< 2\theta$-straight
1-strainer of length $> (1 - \theta)$ at $x$
such that $\Si^o_{x;a,b}$ has diameter $\ge 2\theta^{\frac{5}{2}}$
and contains no $<C_1 \theta$-straight 2-strainer of length $\theta^4$
centered at $x$.

Then there is a point $z \in B(x,2\theta^4)$
which is $\theta$-strongly edgy
relative to $(a_z,b_z)$.
\end{lem}

\proof
Applying \ref{lem:crosssecclsegnew}
on the scale $l = l' = 2\theta^4$,
we obtain that there is a point
$z \in \Si_{x;a,b} \cap B(x,2\theta^4)$
such that $\Si^o_{z;a,b} \subset \Si^o_{x;a,b}$
has diameter $\ge 2\theta^{\frac{5}{2}}$
and admits no $\pihalf$-straight 1-strainer
of length $\half\theta^4$.
The 1-strainer $(a,b)$ is still
$<C_0 \theta$-straight at $z$.

By assumption, the point $z$
admits a $<2\theta$-straight
1-strainer $(a',b')$ of length $> (1 -\theta)$
By Lemma \ref{lem:1stredgynew},
the angles between the two 1-strainers
$(a,b)$ and $(a',b')$ at $z$ are
(after changing the order of the strainer points if necessary)
less than $c'\theta^\half$.

For any point $z' \in \Si_{z;a',b'}$
we have $\vert \angle_z (a',z') - \pihalf \vert \le c\theta$ (by \ref{eq:crosssecperp}),
$\vert \angle_z (a',z') - \angle_z (a,z') \vert \le c'\theta^\half$
and $\vert \angle_z(a,z') - \cangle_z(a,z') \vert \le c''\theta$ (by \ref{ineq:ang1}).
This implies that $\vert \cangle_z(a,z') - \pihalf \vert \le c'''\theta^\half$.
We now apply \ref{lem:projcr}
and project the cross section
$\Si^o_{z;a',b'}$ to $\Si_{z;a,b}$.
This implies that $\Si_{z;a',b'}$ contains no
$\pihalf$-straight 1-strainer of length $\frac{3}{4}\theta^4$
centered at $z$,
since such a strainer would project
to a 1-strainer
of length $\ge \half\theta^4$
in $\Si_{z;a,b}$
which must be $\pihalf$-straight by
Lemma \ref{lem:exstrainnew}.

Similarly, we apply \ref{lem:projcr}
to project
$\Si^o_{z;a,b}$ to $\Si_{z;a',b'}$.
This implies that
$\diam \Si^o_{z;a',b'} \ge \frac{4}{3}\theta^{\frac{5}{2}}$.
Thus, $z$ is indeed $\theta$-strongly edgy
relative to $(a',b')$.
\qed

%
%
By \ref{lem:1stredgynew}(ii), 
at a $\theta$-edgy point $y$, 
there exist no $\theta^3$-long 
$c\theta^{\half}$-straight 2-strainers 
and any two $\theta^3$-long 
$c\theta^{\half}$-straight 1-strainers 
have angle $<c'\theta^{\half}$
(in the sense that their pairs of directions 
are $c'\theta^{\half}$-Hausdorff close subsets in $\Si_yX$).

The almost uniqueness of 1-strainers at edgy points 
extends to uniform neighborhoods:
\begin{lem}[Almost uniqueness of 1-strainers extends]
\label{lem:auniqext}
Let $y$ be $\theta$-edgy relative to $(a,b)$.
Suppose that $(a_i,b_i)$ are $C_0\theta$-straight 1-strainers 
of lengths $\in(1-\theta^{\half},1+\theta^{\half})$
at points $z_i\in B(y,\theta)$ for $i=1,2$. 
Then 
$\angle_{\cdot}(a_1,a_2),\angle_{\cdot}(a_1,b_2),
\angle_{\cdot}(b_1,a_2),$ $\angle_{\cdot}(b_1,b_2)
\not\in[c'\theta^{\half},\pi-c'\theta^{\half}]$
on $B(y,\theta)$, 
if $\la=\la(\theta)$ is sufficiently small. 
\end{lem}
\proof
The strainers $(a_i,b_i)$ are $c'\theta$-straight at $y$, 
cf.\ (\ref{ineq:strqual}),
and hence $c\theta^{\half}$-straight 
with the constant $c$ as in the hypothesis of \ref{lem:1stredgynew}, 
because $c'\theta<c\theta^{\half}$. 
Applying \ref{lem:1stredgynew}(ii) 
with $l=\theta^3$ in the edgy point $y$ 
yields up to switching $a_1$ and $b_1$ that 
$\angle_y(a_1,a_2),\angle_y(b_1,b_2)<c''\theta^{\half}$. 
Due to our condition on the lengths of the 1-strainers $(a_i,b_i)$
it follows that they are $c'''\theta^{\half}$-Hausdorff close 
(as two point subsets), 
which in turn implies that 
$\angle_{\cdot}(a_1,a_2),\angle_{\cdot}(b_1,b_2)<c''''\theta^{\half}$ 
on $B(y,\theta)$. 
\qed

\begin{lem}[Relative position of nearby edgy points]
\label{lem:reledgynew}
Let $x$ be $\theta$-edgy relative to $(a,b)$. 

(i) If $\la=\la(\theta)$ is sufficiently small, 
then all points in $B(x,\theta^3)$
which are $\theta$-weakly edgy
are contained in 
the $\theta^{\frac{15}{4}}$-neighborhood of $axb$.

(ii) Suppose that every point $y \in B(x,\theta)$
admits a $<2\theta$-straight
1-strainer $(a_y,b_y)$ of length $> (1 - \theta)$.
Then for every point $y' \in axb\cap B(x,\theta^3)$
we have $\diam \Si_{y';a_{y'},b_{y'}} < \frac{3}{4} \theta^2$
or $y'$ lies at distance
$<\theta^{\frac{15}{4}}$
from a point $z$
which is $\theta$-strongly edgy
relative to $(a_z,b_z)$.
\end{lem}
\proof
(i) Let $y\in\Si_{x;a,b}\cap A(x,\half\theta^{\frac{15}{4}},2\theta^3)$. 
Then $\Si_{x;a,b}$ contains a 
$c\theta$-straight $\third\theta^{\frac{15}{4}}$-long 1-strainer at $y_1$, 
cf.\ \ref{lem:crosssecclsegnew}(i).
By \ref{lem:exstrainnew}(ii),
every point $u\in ayb\cap B(y,2\theta^3)$
admits a 
$c\theta$-straight $\sixth\theta^{\frac{15}{4}}$-long 1-strainer
contained in $\Si_{u;a,b}$. 
By the metric properties of the maps $\Phi^{a,b}_t$, 
every point $z\in B(y,\theta^3)$ 
outside the $\theta^{\frac{15}{4}}$-neighborhood of $axb$ 
lies at distance $<c'\theta^4$ from such a point $u$ 
(for some such $y$) 
and therefore admits a 
$c''\theta$-straight $>\seventh\theta^{\frac{15}{4}}$-long 1-strainer 
contained in $\Si^o_{z;a,b}$.

Suppose that $z$ is $\theta$-weakly edgy
with respect to some 1-strainer $(a',b')$.
Then by Lemma \ref{lem:auniqext}
the two 1-strainers $(a,b)$ and $(a',b')$
have angles $\le c'\theta^\half$ at $z$,
and we can project $\Si^o_{z;a,b}$ to $\Si_{z;a',b'}$
as in the proof of Lemma \ref{lem:closetoedgy}
to obtain a contradiction.

(ii) Let $y'\in axb\cap B(x,\theta^3)$. 
Suppose that $\Si^o_{y';a,b}$ contains a 
$\pihalf$-straight $\half\theta^{\frac{15}{4}}$-long 1-strainer at $y$. 
Then $y'=ayb\cap\Si_{x;a,b}$ has distance $<c\theta^4$ from $x$. 
By \ref{lem:exstrainnew}, 
$\Si^o_{x;a,b}$ then contains a 
$c'\theta$-straight $\quart\theta^{\frac{15}{4}}$-long 1-strainer at $y'$,
which is also a 
$\pihalf$-straight $>\fifth\theta^{\frac{15}{4}}$-long 1-strainer at $x$. 
This contradicts the $\theta$-edgyness of $x$. 
Thus $\Si^o_{y';a,b}$ contains no 
$\half\theta^{\frac{15}{4}}$-long $\pihalf$-straight 1-strainer at $y'$.

If $\diam \Si^o_{y';a,b} < \theta^{\frac{9}{8}}$,
by assumption there is a  $< 2\theta$-straight
1-strainer $(a',b')$ of length $> (1 - \theta)$ at $y'$.
Projecting $\Si_{y';a,b}$ to $\Si_{y';a',b'}$
yields that $\diam \Si^o_{y';a',b'} < \frac{3}{4}\theta^2$.

Otherwise, \ref{lem:crosssecclsegnew}(ii) 
applied on the scale $l=\theta^{\frac{15}{4}}$ yields that 
$\Si_{y';a',b'}\cap B(y,2\theta^{\frac{15}{4}})$ 
contains a point $z$
such that $\Si^o_{z;a',b'}$ has diameter $\ge \theta^{\frac{9}{4}}$
and admits no $\pihalf$-straight 1-strainer
of length $\half \theta^4$ centered at $z$.
This implies that $z$ is $\theta$-strongly edgy
relative to some 1-strainer $(a_z,b_z)$
as in the proof of Lemma \ref{lem:closetoedgy}.
\qed

\begin{lem}[Almost parallel cross sections of edges]
\label{lem:parcrossedg}
Let $x$ be a $\theta$-edgy point relative to a
$\theta$-straight 1-strainer $(a,b)$ 
with length $> (1 - \theta)$,
and let $y$ be $\theta$-weakly edgy
relative to another 1-strainer $(a',b')$.
We consider the truncated cross sections
$\check\Si_x := \Si_{x;a,b}\cap B(x,\half\theta^3)$
and $\check\Si_y := \Si_{y;a',b'}\cap B(y,\half\theta^3)$.

(i) Suppose that $\check\Si_y$ and $\check\Si_y$ intersect. 
Then $d(x,y)<c\theta^{\frac{7}{2}}$, 
if $\la=\la(\theta)$ is sufficiently small.

(ii) Suppose that $d(x,y) \le c\theta^{\frac{10}{3}}$.
Then the Hausdorff distance $d_H(\check\Si_x,\check\Si_y)$
is less than $\theta^{\frac{99}{30}}$.
\end{lem}
\proof
(i) Let $z$ be one of the intersection points of the cross sections
$\check\Si_x$ and $\check\Si_y$.
By \ref{lem:auniqext}, 
$f_{a,b}-f_{a',b'}$ is $c\theta^{\half}$-Lipschitz 
on $B(z,\half\theta^3)$
(in fact, on $B(z,\theta)$),
and hence $f_{a,b}$ is $c\theta^{\half}$-Lipschitz 
on $\Si_{y;a',b'}\cap B(z,\half\theta^3)$. 
It follows that 
$|f_{a,b}(x)-f_{a,b}(y)|=|f_{a,b}(z)-f_{a,b}(y)|
<c\theta^{\frac{7}{2}}$. 
Since $y$ is contained in the 
$\theta^{\frac{15}{4}}$-neighborhood of $axb$ 
due to \ref{lem:reledgynew}, 
we obtain that 
$d(y_1,y_2)<c'\theta^{\frac{7}{2}}$.

(ii) Let $z = ayb \cap \Si_{x;a,b}$
and consider a point $u \in \check\Si_y$.
By \ref{lem:auniqext}, 
we have $\angle_y(a,a') \le c'\theta^\half$.
Thus, we can apply \ref{lem:projcr} to project
$\check\Si_y$ to $\Si_{y;a,b}$.
In particular, the point $u$
projects to a point $v$
with $d(u,v) \le c''\theta^{\frac{7}{2}}$
and $\vert d(y,v) - d(y,u)\vert \le c'''\theta^3$.

Next, we apply the coarse flow $\Phi^{a,b}_{-f_{a,b}(y)}$
to transport $\Si_{y;a,b}$ to $\Si_{x;a,b}$.
We have $z = \Phi^{a,b}_{-f_{a,b}(y)}$
and set $w:= \Phi^{a,b}_{-f_{a,b}(y)}(v)$.
Our condition that $d(x,y) \le c\theta^{\frac{10}{3}}$
and the metric properties of the flow yield that
$d(v,w) \le c''''\theta^{\frac{10}{3}}$
and $\vert d(w,z) - d(y,u) \vert \le \theta^{\frac{7}{2}}$.
Finally, Lemma \ref{lem:reledgynew}
implies that $d(z,x) \le \theta^{\frac{7}{2}}$.
Thus, the distance between $w$ and $\check\Si_x$
is at most $2\theta^{\frac{7}{2}}$.
(Here, we use again that close to $x$
the cross section $\Si_{x;a,b}$ is almost 1-dimensional,
i.e.\ close to an interval.)

All in all, we conclude that
$d(u,\check\Si_x) \le \theta^{\frac{99}{30}}$.
By switching the roles of $\check\Si_x$ and $\check\Si_y$,
we similarly obtain that every point in $\check\Si_x$
has distance $\le \theta^{\frac{99}{30}}$
to $\check\Si_y$.
This completes the proof.
\qed
For future reference,
we observe that the lemma also holds true
if we replace $\check\Si_y$
by $B(y,\half\tau\theta^3)$
for some $\tau \in (1-\theta, 1+ \theta)$.
The proof for (i) goes through unchanged
and for (ii) it suffices to observe
that the almost 1-dimensionality of $\Si_{y;a',b'}$
near $y$ (Lemma \ref{lem:exstrainnew})
implies $d_H(\check\Si_y, B(y,\half\tau\theta^3)) < \theta^{\frac{99}{30}}$.

\subsection{Necks}
\label{sec:necalex}
A neck occurs where the connected component of a cross section
has small diameter.

\begin{defn}[Necklike points]
\label{def:necky}
A point $y\in B(x,5)$ is called $\theta$-{\em necklike}
relative to a $<2\theta$-straight
1-strainer $(a,b)$ of length $> (1 - \theta)$
(at $y$), 
if $\diam(\Si_{y;a,b})<\theta^2$.

We say that $y \in B(x,5)$ is $\theta$-{\em weakly necklike}
relative to a $< 2\theta$-straight
1-strainer $(a,b)$ of length $> (1 - \theta)$ at $y$
if $\diam(\Si_{y;a,b}) < 2 \theta^2$
and that it is $\theta$-{\em strongly necklike}
relative to such a strainer
if $\diam(\Si_{y;a,b}) < \frac{3}{4} \theta^2$.
\end{defn}

Our new definition allows us to
reformulate part (ii) of Lemma \ref{lem:reledgynew}:
Suppose that for an edgy point $x$,
every point in $B(x,\theta)$
admits a $<2\theta$-straight 1-strainer of length $> (1 - \theta)$
Then every point $y \in axb\cap B(x,\theta^3)$
is $\theta$-strongly necklike
or has distance
$<\theta^{\frac{15}{4}}$
from a point $z$
which is $\theta$-strongly edgy.

Suppose that $x$ 
is $\theta$-necklike relative to $(a,b)$.
Nearby cross sections have comparable diameters:
Using the metric properties of the maps $\Phi^{a,b}_t$ one sees that 
\begin{equation*}
\diam(\Si_{z;a,b})
<(1+c\theta)\theta^2+c\theta|f_{a,b}(z)-f_{a,b}(y)| 
<c'\theta^2
\end{equation*}
for $z\in B(y,\theta)$.


Every segment of length $>\theta$ initiating in 
$B(y,\theta)$ 
must pass through one of the two cross sections 
$f_{a,b}^{-1}(f_{a,b}(y)\pm\frac{9}{10}\theta)$. 
Hence, 
triangle comparison and (\ref{ineq:ang1})
imply for any point $a'$ 
with $d(y,a')>\theta$ that 
\begin{equation}
\label{eq:necalmuni}
\angle_z(a,a'),\angle_z(b,a')
\not\in[c\theta,\pi-c\theta]
\end{equation} 
for $z\in B(y,\frac{\theta}{2})$. 
In particular, 
any $\pihalf$-straight 1-strainer $(a',b')$ of length $>\theta$
at a point in $B(y,\frac{\theta}{2})$
is $c'\theta$-straight,
and $f_{a',b'}-f_{a,b}$ is $c''\theta$-Lipschitz 
on $B(y,\frac{\theta}{2})$. 

If $x$ is $\theta$-necklike
relative to a 1-strainer $(a,b)$
with $\diam\Si_{x;a,b} < \half \theta^2$,
and if all points in $B(x,\theta)$
admit $<2\theta$-straight 1-strainers
of length $> (1-\theta)$
we can conclude from the above observations 
as in the proof of Lemma \ref{lem:closetoedgy}
that all points in $B(x,\theta)$
are also $\theta$-strongly necklike.

We now deduce 
that cross sections of nearby necklike points 
are almost parallel.
\begin{lem}[Almost parallel cross sections of necks]
\label{lem:almpcrsecnec}
Let $x$ be $\theta$-necklike relative to $(a,b)$. 
Furthermore, let $y \in B(x,\theta)$
be $\theta$-weakly necklike
with respect to a 1-strainer $(a',b')$.

(i) If $\Si_{x;a,b}$ and $\Si_{y;a',b'}$
have nonempty intersection, then 
$d(x,y) < c'\theta^3$

(ii) If $d(x,y) < \theta^{\frac{11}{6}}$,
then the Hausdorff distance 
of $\Si_{x;a,b}$ and $\Si_{y;a',b'}$
is less that $\theta^{\frac{5}{3}}$.
\end{lem}
\proof
The proof is closely related
to the one for edges, i.e.\ \ref{lem:parcrossedg}.

(i)
Let $z$ be one of the intersection points of the cross sections
$\Si_{x;a,b}$ and $\Si_{y;a',b'}$.
By our discussion above,
$f_{a,b}-f_{a',b'}$ is $c\theta$-Lipschitz 
on $B(z,\theta^2)$,
and hence $f_{a,b}$ is $c\theta^{\half}$-Lipschitz 
on $\Si_{y;a',b'}$.
It follows that 
$|f_{a,b}(x)-f_{a,b}(y)|=|f_{a,b}(z)-f_{a,b}(y)|
<c'\theta^3$. 

(ii)
Consider a point $u \in \check\Si_y$.
By \ref{eq:necalmuni}, 
we have $\angle_y(a,a') \le c'\theta$.
When projecting $\Si_{y;a',b'}$ to $\Si_{y;a,b}$,
we map $u$ to a point $v$
with $d(u,v) \le c''\theta^3$
by \ref{lem:projcr}.
The coarse flow $\Phi^{a,b}_{-f_{a,b}(y)}$
transports $v \in \Si_{y;a,b}$ to 
some point $w \in \Si_{x;a,b}$
with $d(v,w) \le c''' \theta^{\frac{11}{6}}$.

This shows that $d(u,\Si_{x;a,b}) \le \theta^{\frac{5}{3}}$.
Again, we switch the roles of $\Si_{x;a,b}$ and $\Si_{y;a',b'}$
to complete the proof.
\qed

\section{Locally volume collapsed 3-orbifolds are graph}

\subsection{Setup and formulation of main result}
\label{sec:setup}

Let $(O,g)$ be a closed connected Riemannian 3-orbifold 
which does {\em not} have nonnegative sectional curvature, 
$sec\not\geq0$. 

\begin{defn}[Curvature scale]
For $-b^2\in[-1,0)$ 
we define the {\em $-b^2$-(sectional) curvature scale} in a point $x\in O$ 
as the maximal radius $\rho_{-b^2}(x)\in(0,\infty)$ such that 
the rescaled ball 
$B_{\rho_{-b^2}(x)^{-2}g}(x,1)=\rho_{-b^2}(x)^{-1}\cdot B_g(x,\rho_{-b^2}(x))$
has sectional curvature $sec\geq-b^2$.
\end{defn}
Note that 
\begin{equation}
\label{ineq:csc}
b\rho_{-1}\leq\rho_{-b^2}\leq\rho_{-1} .
\end{equation}
The function $\rho_{-b^2}$ is continuous on $O$. 
More precisely, 
$\rho_{-b^2}$ {\em does not oscillate too fast} in the sense that 
for $0<\la<1$ holds 
\begin{equation}
\label{eq:osccurvsc}
(1-\la)\rho_{-b^2}(x)\leq\rho_{-b^2}\leq(1+\la)\rho_{-b^2}(x)
\end{equation}
on $B(x,\la\rho_{-b^2}(x))$.

The rescaled balls 
$B_{\rho_{-b^2}(x)^{-2}g}(x,1)$ 
are Alexandrov balls with curvature $\geq-b^2$ and radius $\leq1$ 
in the sense of definition \ref{def:alexb}.  

The purpose of this paper 
is to study the geometry and topology of 3-orbifolds 
which are locally collapsed relative to the curvature scale. 
\begin{defn}[Local volume collapse]
\label{def:volcoll}
Let $v>0$ 
and let $\si:O\to(0,\infty)$ be some 
(not necessarily continuous) function. 
We say that $(O,g)$ is $v$-{\em collapsed} 
at the scale $\si$, 
if for all points $x$ holds 
$\vol(B_{\si(x)^{-2}g}(x,1))<v$, 
equivalently,
$\vol(B_g(x,\si(x)))<v\si(x)^3$. 

If $sec\not\geq0$,
we say that $(O,g)$ is $(v,-b^2)$-{\em collapsed}, 
if it is $v$-collapsed at the scale $\rho_{-b^2}$. 
\end{defn}
Note that 
if $(O,g)$ is locally $v$-collapsed at some scale $\si\leq\rho_{-b^2}$, 
then Bishop-Gromov volume comparison yields that 
it is locally $(v',-b^2)$-collapsed 
with 
$v'=\frac{vol(B_{-b^2}(1))}{vol(B_0(1))}v
\leq\frac{vol(B_{-1}(1))}{vol(B_0(1))}v$
(independent of $-b^2$!). 
Here $B_{-b^2}(1)$ denotes the unit 3-ball with $sec\equiv-b^2$. 

Strongly volume collapsed Riemannian 3-orbifolds are, 
on the scale of their collapse, 
close to Alexandrov spaces of dimension $\leq2$
and with curvature $\geq-b^2$,  
and the volume collapse translates 
into the shortness of $2\half$-strainers, 
cf.\ section~\ref{sec:strainalex}). 
\begin{lem}
\label{lem:collstrain}
For $\la>0$ exists $v=v(\la)>0$ such that: 

If $(O,g)$ is $(v,-b^2)$-collapsed
and $x \in O$,
then $\bar\theta_{2\half}$-straight $2\half$-strainers
in the Alexandrov ball $\rho_{-b^2}(x)^{-1} B(x,\rho_{-b^2}(x))$
of curvature $\ge -b^2 \ge -1$
have length $<\la$.
\end{lem}
\proof
Suppose that the Riemannian 3-orbifolds $(O_i,g_i)$ 
are $(\frac{1}{i},-b_i^2)$-collapsed 
but contain points $x_i$ which admit
such that there are 
$\bar\theta_{2\half}$-straight $2\half$-strainers 
of length $\geq\la$
in the balls $\rho_{-b_i^2}(x)^{-1} B(x, \half \rho_{-b_i^2}(x))$. 
Then the rescaled balls 
$\rho_{-b^2}(x_i)^{-1}B_(x_i,1)$ 
Gromov-Hausdorff subconverge 
to an Alexandrov ball with dimension $\leq2$ 
and curvature $\geq-1$ 
which admits 
$\bar\theta_{2\half}$-straight $2\half$-strainers 
of length $\la$,
a contradiction. 
\qed

\medskip

We will require some additional regularity for our Riemannian orbifolds. 
The conclusions on the global topology of collapsed 3-orbifolds
are valid without this regularity condition, 
but it is technically convenient because 
it avoids the use of 
(an orbifold version of) 
Perelman's Stability Theorem for Alexandrov spaces 
and is expected to be satisfied by the output of the Ricci flow on 3-orbifolds. 
In the next definition, 
$\om_3$ denotes the volume of the euclidean unit 3-ball. 
(Compare \cite[7.4]{Perelman} and \cite[Thm.\ 1.3]{KL_coll}.)

\begin{defn}[Local curvature control]
\label{def:loccc}
Fix numbers $s_0\in\N$, $v_0\in(0,\om_3)$, 
a function $K:(0,\om_3)\to(0,\infty)$ 
and a scale function $\si:O\to (0,\infty)$. 
We say that 
$(O,g)$ 
has $(v_0,s_0,K)$-{\em curvature control} below scale $\si$,  
if the following holds:
If $\vol B(x,r)\geq vr^3$ 
for $v\in[v_0,\om_3)$ and $r\in(0,\si(x)]$, then 
$\|\nabla^sR\|\leq K(v)r^{-2-s}$ on $B(x,r)$, 
equivalently, 
$\|\nabla^sR\|\leq K(v)$ on the rescaled ball $r^{-1}\cdot B(x,r)$ 
for $s=0,\dots,s_0$. 
\end{defn}
We will apply this notion in the following situation.
\begin{lem}
\label{lem:smconv}
Let $(O_i,g_i)$ be a sequence of Riemannian 3-orbifolds as above 
with $(v_i,s_0,K)$-curvature control below scale $\rho_{-b^2}$,
where $v_i\to0$. 
Furthermore, let $x_i\in O_i$ be points and $\la_i\to0$ positive numbers. 

Then, if $s_0$ is sufficiently large, 
the sequence of rescaled pointed orbifolds 
$(\la_i\rho_{-b^2}(x_i))^{-1}\cdot(O_i,x_i)$
subconverges 
either in the Gromov-Hausdorff sense to an Alexandrov space 
with curvature $\geq0$ and dimension $\leq2$, 
or in the ${\mathcal C}^5$-topology to a 
${\mathcal C}^{10}$-smooth complete Riemannian 3-orbifold with $\sec\geq0$. 
\end{lem}
\proof
Clearly, we have Gromov-Hausdorff subconverge to 
a pointed Alexandrov space $(X,x_0)$ 
with curvature $\geq0$ and dimension $\leq3$. 
If $\dim(X)=3$, the convergence can be improved 
using our assumption of local curvature control. 
The approximating pointed 3-orbifolds 
$(\la_i\rho_{-b^2}(x_i))^{-1}\cdot(O_i,x_i)$ 
are uniformly noncollapsed. 
Indeed, for any $r>0$, we have 
$\vol B_{\la_i^{-2}\rho_{-b^2}(x_i)^{-2}g_i}(x_i,r)
>\half\vol_3 B(x_0,r)>0$ 
for large $i$, 
where volume in $X$ is measured 
with respect to the 3-dimensional Hausdorff measure. 
Thus the $(v_i,s_0,K)$-curvature control 
on $(\la_i\rho_{-b^2}(x_i))^{-1}\cdot O_i$
below the scales 
$\la_i^{-1}\to\infty$ applies for large $i$ 
and yields on the balls 
$B_{\la_i^{-2}\rho_{-b^2}(x_i)^{-2}g_i}(x_i,r)$ 
uniform bounds on the curvature tensor and its covariant derivatives 
up to order $s_0$. 
The smoothness of the limit and the convergence follow 
from an orbifold version of well-known compactness results 
for pointed Riemannian manifolds with bounds on curvature 
and some of its derivatives, 
e.g.\ from the following special case of a result in \cite{locstr}:
\begin{thm}
Let $r_0,v_0>0$ 
and let $C:(0,\infty)\to(0,\infty)$ be a function.  
Then the space of pointed complete 
${\mathcal C}^{\infty}$-smooth Riemannian 3-orbifolds $(O,p)$ 
such that $\vol B(p,r_0)\geq v_0$ 
and such that on every ball $B(p,r)$ around $p$ 
the curvature tensor and its covariant derivatives 
up to order 20
on $B(p,r)$ are bounded by the constant $C(r)$,
is precompact in the pointed ${\mathcal C}^{15}$-topology. 
Thus, every sequence of such orbifolds $(O_i,p_i)$
subconverges ${\mathcal C}^5$-smoothly 
to a ${\mathcal C}^{10}$-smooth Riemannian 3-orbifold. 
\end{thm}
This completes the proof of \ref{lem:smconv}. 
\qed

\medskip
The main result of this paper is
(compare \cite[Theorem 7.4]{Perelman},
\cite[Theorem 0.2]{MT2}
and \cite[Theorem 1.3]{KL_coll}):
\begin{thm}
\label{thm:main}
Let $s_0\in\N$ 
and let $K:(0,\om_3)\to(0,\infty)$ be a function. 
If $s_0$ is sufficiently large,
then there exists a constant 
$v_0=v_0(s_0,K)\in(0,\om_3)$ such that: 

If $(O,g)$ is $(v_0,-1)$-collapsed,
has $(v_0,s_0,K)$-curvature control 
below the scale $\rho_{-1}$
and contains no bad 2-suborbifolds,
then $O$ admits a metric with $\sec \geq 0$ 
or can be decomposed
by finitely many surgeries
into components which are spherical or graph.
\end{thm}
\begin{rem}
Unlike in the manifold case we cannot conclude that $O$ is always graph, 
because there are non-graph 3-orbifolds 
admitting nonnegatively curved (e.g.\ spherical or euclidean) metrics.
\end{rem}
One can reduce to the case of a lower diameter bound 
relative to a curvature scale
by suitably rescaling. 
\ref{thm:main} follows from: 
\begin{thm}
\label{thm:mainwdb}
Let $s_0\in\N$ 
and let $K:(0,\om_3)\to(0,\infty)$ be a function. 
If $s_0$ is sufficiently large,
then there exists a constant 
$v_0=v_0(s_0,K)\in(0,\om_3)$ such that: 

If for some $-b^2\in[-1,0)$ 
the orbifold $(O,g)$ is $(v_0,-b^2)$-collapsed, 
satisfies $\rad(O,\cdot)\geq\half\rho_{-b^2}$,
has $(v_0,s_0,K)$-curvature control below the scale $\rho_{-b^2}$, 
and contains no bad 2-suborbifolds,
then $O$ can be decomposed
by finitely many surgeries
into components which are spherical or graph.
\end{thm}

\no
{\em Proof that \ref{thm:mainwdb} implies \ref{thm:main}.}
Suppose that $(O_i)$ is a sequence of $(v_i,-1)$-collapsed orbifolds 
with $(v_i,s_0,K)$-curvature control below the scale $\rho_{-1}$
where $v_i\to0$. 
Then we must show that the $O_i$ satisfy the conclusion of \ref{thm:main} 
for infinitely many $i$.

Note that for all $b\in(0,1]$ 
the $O_i$ are $(b^{-3}v_i,-b^2)$-collapsed 
and have $(v_i,s_0,K)$-curvature control below the scale $\rho_{-b^2}$, 
cf.\ (\ref{ineq:csc}). 
Hence we are done, if for some $-b^2\in[-1,0)$ holds 
$\rad(O_i)\geq\half\rho_{-b^2}$ for all sufficiently large $i$. 

Otherwise, 
after passing to a subsequence,
there exist sequences of numbers $-b_i^2\to0$ and points $x_i\in O_i$ 
such that $\rad(O_i,x_i)<\half\rho_{-b_i^2}(x_i)$.
It follows that 
$\rho_{-b_i^2}\equiv const_i\geq\diam(O_i)$
and we have collapse to the point in the sense of 
$\diam(O_i)\cdot(-\min\sec_{O_i})^{\half}\to0$. 
We rescale and increase the $-b_i^2$ so that $\diam(O_i)=1$ 
and $\min\sec_{O_i}=-b_i^2\to0$. 
Then $\rho_{-b_i^2}\equiv1$. 

If $v'_i=\vol(O_i)\to0$, 
then the $O_i$ are $(v'_i,-b_i^2)$-collapsed 
with $(v_i,s_0,K)$-curvature control below the scales $\rho_{-b_i^2}$ 
and with $\rad (O_i) \ge \half\equiv\half \rho_{-b_i^2}$, 
and we are done by \ref{thm:mainwdb}. 

Otherwise, 
after passing to a subsequence,
we have a lower volume bound $\vol(O_i)\geq v'>0$ and, 
due to the curvature control, 
uniform (global) bounds on the curvature tensor and its covariant derivatives 
up to order $s_0$. 
If $s_0$ is large enough, 
it follows that 
the $O_i$ subconverge, say, in the ${\mathcal C}^5$-topology 
to a ${\mathcal C}^5$-Riemannian 3-orbifold $O_{\infty}$ 
with $\sec\geq0$. 
In particular, 
infinitely many $O_i$ are diffeomorphic to $O_{\infty}$ 
and therefore admit metrics with $\sec\geq0$. 
\qed

\medskip
We fix some arbitrary (more than) sufficiently large value for $s_0$,
say $s_0:=2010$. 

For the rest of this section
we make the following assumption on the orbifolds we work with: 
\begin{ass}
$(O,g)$ is a closed connected Riemannian 3-orbifold 
such that $sec\not\geq0$ 
and $\rad(O,\cdot)\geq\half\rho_{-b^2}$ where $-b^2\in[-1,0)$. 
\end{ass}

Together with Corollary \ref{cor:graphgeom},
Theorem \ref{thm:main} implies the following

\begin{thm}
\label{thm:mainclosed}
Let $s_0\in\N$ 
and let $K:(0,\om_3)\to(0,\infty)$ be a function. 
If $s_0$ is sufficiently large,
then there exists a constant 
$v_0=v_0(s_0,K)\in(0,\om_3)$ such that:

If $(O,g)$
is closed and $(v_0,-1)$-collapsed,
has $(v_0,s_0,K)$-curvature control 
below the scale $\rho_{-1}$
and contains no bad 2-suborbifolds,
then $O$
either admits a metric with $\sec \geq 0$,
or satisfies Thurston's Geometrization Conjecture.
\end{thm}

\subsection{Conical approximation and humps}

Sufficient local volume collapse leads to 
good local approximation by cones of dimension $\leq2$
on scales comparable to the curvature scale. 
Proposition~\ref{prop:locapproxconealex} for $d=2$ and 
Lemma~\ref{lem:collstrain} imply: 
\begin{prop}[Uniform conical approximation 
relative to the curvature scale]
\label{prop:locapproxconeorbi}
For $0<\si,\mu<1$ exist 
a scale $0<s_1=s_1(\si,\mu)<<\si$ 
and a rate of collapsedness $v=v(\si,\mu)>0$
such that the following holds:

If $(O,g)$ is $(v,-b^2)$-collapsed, 
then it can in every point $x$ 
be $\mu$-well approximated 
(cf.\ \ref{def:locapp}) 
on some scale $s(x)\in[s_1\rho_{-b^2}(x),\si\rho_{-b^2}(x)]$ 
by a cone of dimension 1 or 2, 
i.e.\ by the open interval $(-1,1)$,
by the half-open interval $[0,1)$, 
or by the cone of radius 1 over a circle or an interval 
of diameter $\leq\pi$. 
\end{prop}
\proof
The assertion follows with 
$s_1=s_1(\si,\mu):=s_1(2,\si,\mu)$ from 
Proposition~\ref{prop:locapproxconealex} and 
$v=v(\si,\mu):=v(s_{2\half}(\si,\mu))$
from Lemma~\ref{lem:collstrain}. 
\qed

\medskip
We fix some small value $\si\in(0,\tenth)$ for the upper bound 
on the local scales of approximation. 

Next we wish to divide our orbifold 
into regions with and without good 1-strainers. 
The region with good 1-strainers has locally almost product geometry. 
It consists of the points where the bases of approximating cones 
provided by \ref{prop:locapproxconeorbi} have diameter $\approx\pi$. 

For small $\theta>0$, 
let $\mu_0(\theta)>0$ be the constant given by 
\ref{lem:nhump}. 
For $\mu\in(0,\mu_0(\theta)]$, 
we call a point $x\in O$ a $(\theta,\mu,-b^2)$-{\em hump}, 
if $O$ can in $x$ be $\mu$-well approximated on the scale $s(x)$ 
by a cone with base of diameter $<\pi-\frac{\theta}{2}$, 
i.e.\ by the half-open interval $[0,1)$ 
or the cone of radius 1 over a circle or an interval 
with diameter $<\pi-\frac{\theta}{2}$. 
In other words, 
$x$ is a $(\theta,\mu)$-hump of the rescaled ball 
$B_{\rho_{-b^2}(x)^{-2}g}(x,1)$ 
in the sense of Definition~\ref{def:humpalex} for $d=2$. 
If $x$ is no $(\theta,\mu,-b^2)$-hump, 
then $O$ can in $x$ 
be $\mu$-well approximated on the scale $s(x)$ 
by a cone with base of diameter $\geq\pi-\frac{\theta}{2}$, 
i.e.\ by the open interval $(-1,1)$ 
or the cone of radius 1 over a circle or an interval 
with diameter $\geq\pi-\frac{\theta}{2}$. 

We denote by $H=H_{\theta,\mu,-b^2}\subset O$
the subset of $(\theta,\mu,\rho_{-b^2}(x))$-humps 
and by $S=S_{\theta,\mu,-b^2}\subset O$ the open subset of points $x$ 
admitting equilateral 1-strainers 
which are $<\theta$-straight with length in
$(\frac{1}{11}s_1(\si,\mu), \frac{3}{22}s_1(\si,\mu))$
in the Alexandrov ball
$(\rho_{-b^2}(x))^{-1}B(x,\rho_{-b^2}(x))$
of curvature $\ge -1$.

Throughout the following chapters, we will
abbreviate this property by saying that the points $x \in S_{\theta,\mu,-b^2}$
admit $<\theta$-straight (equilateral) 1-strainers
with length in the interval $(\frac{1}{11}s_1(\si,\mu)\rho_{-b^2}(x),\frac{3}{22}s_1(\si,\mu)\rho_{-b^2}(x))$.


Due to our bound on the approximation accuracy $\mu$, 
the implications of \ref{lem:nhump} hold: 
A $(\theta,\mu,-b^2)$-hump $x$ admits no 
$\frac{\theta}{4}$-straight 1-strainers of length $\frac{1}{111}s(x)$, 
but all points in the closed annulus 
$\ol A(x;\frac{1}{10}s(x),\frac{9}{10}s(x))$
do admit $\frac{\theta}{11}$-straight 1-strainers 
of length $>\frac{1}{11}s(x)\geq\frac{1}{11}s_1(\si,\mu)\rho_{-b^2}(x)$, 
i.e.\ 
$\ol A(x;\frac{1}{10}s(x),\frac{9}{10}s(x))\subset S_{\theta,\mu,-b^2}$. 
On the other hand, 
if $x$ is no $(\theta,\mu,-b^2)$-hump, 
then it admits $<\theta$-straight 1-strainers of length 
$>\frac{99}{100}s(x)>\frac{1}{11}s_1(\si,\mu)\rho_{-b^2}(x)$. 
Hence $O=H_{\theta,\mu,-b^2}\cup S_{\theta,\mu,-b^2}$. 
\begin{prop}[cf.\ \ref{prop:finhumpalex}]
\label{prop:finhumporbi}
If $\mu\leq\mu_0(\theta)$ 
and $(O,g)$ is $(v(\si,\mu),-b^2)$-collapsed
(cf.\ \ref{prop:locapproxconeorbi}), 
then there exist finitely many 
$(\theta,\mu,-b^2)$-humps $x_j\in H_{\theta,\mu,-b^2}$
such that 
\begin{equation*}
O=\left(\bigcup_jB(x_j,\frac{1}{10}s(x_j))\right)
\cup S_{\theta,\mu,-b^2} .
\end{equation*}
Moreover,
$d(x_j,x_k)>\frac{9}{10}s(x_j)$ for $j\neq k$. 
\end{prop}
\proof 
We proceed as in the proof of \ref{prop:finhumpalex}.

Again, there can be no infinite sequence of points 
$x_1,x_2,\dots\in H_{\theta,\mu,-b^2}-S_{\theta,\mu,-b^2}$
such that for all $k$
holds $x_k\in B(x_{k+1},\frac{1}{10}s(x_{k+1}))$
and $x_{k+1}\not\in B(x_k,\frac{1}{10}s(x_k))$, 
because then 
$s(x_1)<\frac{1}{9^{k-1}}s(x_k)$. 
But due to the continuity of the curvature scale $\rho_{-b^2}$,
cf.\ (\ref{eq:osccurvsc}), 
and the compactness of $O$ 
the scales take also in this situation values in a bounded interval, 
namely in $[s_1\min\rho_{-b^2},\si\min\rho_{-b^2}]$. 

The rest of the proof goes through without change. 
\qed



%
\subsection{The Shioya-Yamaguchi blow-up}
\label{sec:blowup}
We recall the Shioya-Yamaguchi blow-up argument,
see \cite[\S 3, Key Lemma 3.6]{SY}. 
To simplify things, 
we restrict ourselves to certain special situations. 
Some of our arguments are different;
we also treat some additional cases
not mentioned there explicitely.

\subsubsection{General discussion}
\label{sec:gen}

%

Consider the following situation.
Let $B(p_i,1)$ be a sequence of 
$d$-dimensional Riemannian orbifold balls 
(i.e.\ open metric balls of radius 1 
in complete Riemannian $d$-orbifolds without boundary)
with curvature $sec\geq-1$
which collapse to an Alexandrov ball 
of strictly smaller dimension $1\leq k<d$,
\begin{equation}
\label{coll}
B(p_i,1)\lra X=B(x,1) .
\end{equation}
We suppose furthermore that the collapse limit $X$ is 
(0-){\em conelike} 
in the sense that every segment initiating in $x$ extends to length 1, 
and that the closed balls $\ol B(p_i,\half)$
are {\em not discal}.
Note that the conelikeness of $X$ implies 
that for any fixed $\eps>0$ the annulus $A(p_i,\eps,1-\eps)$ 
contains for large $i$ 
no critical points of the distance function $d(p_i,\cdot)$
and in particular $\ol B(p_i,\half)$ 
is a topological suborbifold with boundary. 

Let $\hat p_i\in B(p_i,1)$ be any sequence of points with 
$d(\hat p_i,p_i)\to0$. 
The distance function $d(\hat p_i,\cdot)$ must have critical values
in $(0,\half)$,
because $\ol B(p_i,\half)$ is not discal. 
Let $\de_i$ be the maximal critical value in $(0,\half)$, 
and let $q_i\in B(p_i,\half)$ be a critical point at distance 
$d(\hat p_i,q_i)=\de_i$ from $\hat p_i$. 
Then $\de_i\to 0$ because $X$ is conelike. 
For any constant $c>1$ holds that 
\begin{equation}
\label{eq:nocrit}
\ol B(p_i,\half)\cong\ol B(\hat p_i,c\de_i)
\end{equation}
for sufficiently large $i$,
i.e.\ the topology of the balls $B(p_i,1)$ 
is concentrated near their centers. 
(Note that, also due to the conelikeness of $X$, 
there exists a common gradient like vector field for 
$d(p_i,\cdot)$ and $d(\hat p_i,\cdot)$,
and so $\ol B(p_i,\half)\cong\ol B(\hat p_i,\half)$.)  

To help revealing the local topology at the $p_i$, 
we form (modulo passing to a subsequence) the {\em blow-up limit}
\begin{equation}
\label{partcoll}
(\de_i^{-1}B(p_i,1),\hat p_i)\lra (Y,y_0) .
\end{equation}
The limit space $Y$ is a noncompact Alexandrov space 
with dimension $\geq k$ and curvature $\geq0$. 
In particular, the Soul Theorem applies. 
Moreover $q_i\to z$ 
with $z$ a critical point of $d(y_0,\cdot)$ at distance $d(y_0,z)=1$. 

If no collapse happens any more in the blow-up limit (\ref{partcoll}),
i.e.\ if $\dim Y=d$, 
then we need topological stability in order to relate the topologies 
of the balls $B(p_i,1)$ and $Y$. 
For instance, 
if $Y$ and the convergence in (\ref{partcoll}) are sufficiently smooth
as will be the case in the situations considered later in the paper, 
say $Y$ is ${\mathcal C}^{10}$-smooth 
and the convergence is ${\mathcal C}^5$-smooth, 
then one can argue as follows. 
There exist $r,\eps>0$ 
and a smooth vector field $V$ on $Y-\ol B(y_0,\frac{r}{2})$
such that for all $y\not\in\ol B(y_0,\frac{r}{2})$
the vector $V(y)$ has angles $\geq\pihalf+\eps$ with all segments $yy_0$, 
compare the proof of the Soul Theorem. 
We regard $\de_i^{-1}B(\hat p_i,2r\de_i)$ as embedded in $Y$ for large $i$. 
In particular, $\hat p_i\to y_0$ 
and $V$ is on a neighborhood of $\D B(y_0,r)$
gradient-like not only for $d(y_0,\cdot)$ 
but also for $\de_i^{-1}d_{B(p_i,1)}(\hat p_i,\cdot)$, 
viewed as a function on part of $Y$. 
Hence $\de_i^{-1}\ol B(\hat p_i,r\de_i)$ is isotopic to 
$\ol B(y_0,r)$ in $Y$, 
and with (\ref{eq:nocrit}) 
we see that $\ol B(p_i,\half)$ is for large $i$ homeomorphic 
to the closed disc bundle in the normal bundle of the soul of $Y$,
in other words, 
to a (small) closed tubular neighborhood 
of the soul. 

The blow-up limit (\ref{partcoll}) is in general 
still a collapse to lower dimension. 
The aim of the following discussion is to find situations 
when the drop of dimension is strictly smaller than 
for the original collapse 
(\ref{coll}). 

The following construction is 
as in the proof of Soul Theorem:
Let $\xi\in\Si_xX$ be a direction at $x$. 
Due to conelikeness
it is represented by a (unique) segment $\si_{\xi}$ emanating from $x$. 
Fix some $t_0\in(0,1)$, say $t_0=\tenth$, and 
let $\si_{\xi}^i\in B(p_i,1)$ be a sequence of points 
converging to the point $\si_{\xi}(t_0)$ on $\si_{\xi}$ 
at distance $t_0$ from $x$.
(Our choice of the $\si_{\xi}(t_0)$ is independent 
of the choice of the $\hat p_i$.) 
The segments $\hat p_i\si_{\xi}^i$ 
subconverge to a 
ray $\rho_{\xi}$ in $Y$ 
emanating from $y_0$. 
Moreover,
the normalized distance functions 
$\de_i^{-1}(d(\si_{\xi}^i,\cdot)-d(\si_{\xi}^i,\hat p_i))$
on the rescaled balls 
$\de_i^{-1}B(p_i,1)$
subconverge (due to Arzela-Ascoli) 
to a concave 1-Lipschitz function $\beta_{\xi}$ on $Y$ 
with $\beta_{\xi}(y_0)=0$ 
which decays along $\rho_{\xi}$ with extremal slope $-1$,
$\beta_{\xi}(\rho_{\xi}(t))=-t$,  
where we use a unit speed parametrization $\rho_{\xi}(t)$ 
starting at $\rho_{\xi}(0)=y_0$. 
In fact,
every point $y\in Y$ is the initial point of a ray $\rho_{\xi}^y$
along which $\beta_{\xi}$ decays with slope $-1$. 
In particular,
the level sets of $\beta_{\xi}$ have no interior points. 
The comparison of $\beta_{\xi}$ with the Busemann function 
$b_{\xi}=\lim_{t\to\infty}(d(\rho_{\xi}(t),\cdot)-d(\rho_{\xi}(t),y_0))$ 
associated to the ray $\rho_{\xi}$ 
is given by the inequality 
\begin{equation*}
\beta_{\xi}\leq b_{\xi} .
\end{equation*}
To verify this, 
let $\si_{\xi}^i(a)$ denote the point on $\hat p_i\si_{\xi}^i$ 
at (unrescaled) distance $a\de_i$ from $\hat p_i$. 
Then 
$d(\si_{\xi}^i,\cdot)-d(\si_{\xi}^i,\hat p_i)
\leq d(\si_{\xi}^i(a),\cdot)-d(\si_{\xi}^i(a),\hat p_i)$
for large $i$ 
and hence 
$\beta_{\xi}\leq d(\rho_{\xi}(a),\cdot)-a$
for all $a>0$. 
Letting $a\to\infty$ yields the inequality. 
As a consequence, 
the convex suplevel sets of $\beta_{\xi}$
are smaller (not larger) than the corresponding suplevels of $b_{\xi}$. 

Since the $q_i$ are critical for $d(\hat p_i,\cdot)$,
we have that 
$\tilde\angle_{q_i}(\hat p_i,\si_{\xi}^i)\leq\pihalf$
and
$\liminf_{i\to\infty}
\de_i^{-1}(d(\si_{\xi}^i,q_i)-d(\si_{\xi}^i,\hat p_i))
\geq0$.
So 
\begin{equation*}
\beta_{\xi}(z)\geq0 
\end{equation*}
and $z$ is contained in the totally convex subset 
$\cap_{\xi}\{\beta_{\xi}\geq0\}=\{\min_{\xi}\beta_{\xi}\geq0\}$.

The blow-up {\em expands} $\Si_xX$ in the following sense. 
When doing the above construction for two directions $\xi$ and $\xi'$ 
at the same time, 
one obtains for every point $y\in Y$ 
a pair of rays $\rho_{\xi}^y$ and $\rho_{\xi'}^y$ satisfying
\begin{equation}
\label{ineq:blexpand}
\angle_{Tits}(\rho_{\xi}^y,\rho_{\xi'}^y)
:=\lim_{t\to\infty}\cangle_y(\rho_{\xi}^y(t),\rho_{\xi'}^y(t))
\geq\angle_x(\xi,\xi').
\end{equation}
The construction can be done for any finite subset $A\subset\Si_xX$ 
and hence yields weakly expanding maps $\eps_{y,A}:A\to\tits Y$. 

We will use the following observations: 
\begin{lem}
\label{lem:blowup}
(i) 
The blow-up limit $Y$ is not isometric to a euclidean space. 

(ii) 
If $\Si_xX$ contains an embedded unit $l$-sphere
(i.e.\ with $\sec\equiv1$), 
then so does $\tits Y$ and 
$Y$ splits off an $\R^{l+1}$-factor. 
If $\Si_xX$ contains an embedded unit $l$-hemisphere, 
then so does $\tits Y$
and $Y$ contains an isometrically embedded copy 
of the $(l+1)$-dimensional euclidean halfspace
(and in particular splits off an $\R^l$-factor).
%
\end{lem}
\proof
(i) $d(y_0,\cdot)$ has a critical point (at distance 1). 

(ii) Choose $A$ as the union of $l+1$ pairs of antipodes 
which span the embedded unit sphere 
(corresponding to coordinate axes).
Then the expanding map $A\to\tits Y$ must be an isometric embedding
and the assertion follows from the Splitting Theorem. 
The second assertion follows similarly by applying the first assertion 
to the boundary $(l-1)$-sphere of the embedded $l$-hemisphere. 
\qed

\subsubsection{The case of flat conical limits with dimension $\leq2$}
\label{sec:sy2dlim}
We apply the general discussion above in certain special situations. 
Note that always $\dim Y\geq\dim X$.
We aim now to achieve that $\dim Y>\dim X$  
by making a good choice of the $\hat p_i$. 

\medskip
{\em Collapse to a flat $k$-disc.} 
Suppose that $X$ is isometric to the euclidean unit $k$-disc, $k\geq1$. 
By \ref{lem:blowup}, 
$Y$ splits off an $\R^k$-factor and $Y\not\cong\R^k$.
Hence always $\dim(Y)>k$, 
independently of the choice of the $\hat p_i$.

\medskip
{\em Noses: Collapse to the half-open interval.} 
Suppose that $X=[0,1)$ with $x=0$. 
There is a unique direction $\xi$ at $x=0$.
We choose $\hat p_i$ as a ``tip'' of the nose,
i.e.\ as a maximum of 
$d(\si_{\xi}^i,\cdot)$. 
Then $\hat p_i\to x$ and the choice of the $\hat p_i$ 
is admissible in the sense that $d(p_i,\hat p_i)\to0$. 
The base point $y_0$ is a maximum of $\beta_{\xi}$
and hence $\beta_{\xi}(z)=\beta_{\xi}(y_0)=0$. 
If $\dim(Y)=1$, then $Y$ is a halfline 
since $Y\not\cong\R$ by \ref{lem:blowup}(i), 
and $\beta_{\xi}=b_{\xi}$ has a unique maximum.
This contradicts $z\neq y_0$. 
Thus $\dim(Y)\geq2$.

Note that $Y$ contains a flat half-strip, 
but nevertheless its geometry is in general not rigid. 

\medskip
{\em Collapse to a flat 2-disc with cone point or a sector.} 
Suppose that $X$ is the cone of radius 1 
over a circle or interval with diameter $<\pi$. 
We generalize the cases of noses to humps 
by adapting the argument in \cite[\S 3]{SY} to this case. 

Let $A\subset\Si_xX$ be a finite subset such that 
$\sum_{\xi\in A}d(\si_{\xi}(t_0),\cdot)$ 
has a unique maximum in $x$. 
We choose the $\hat p_i$ as maxima of the corresponding functions 
$\sum_{\xi\in A}d(\si_{\xi}^i,\cdot)$. 
Then $\hat p_i\to x$. 
The point $y_0$ is a maximum of 
$\sum_{\xi\in A}\beta_{\xi}$. 
It follows that $\beta_{\xi}(z)=0$ for all $\xi\in A$, 
and the totally convex subset 
$\cap_{\xi\in A}\{\beta_{\xi}\geq0\}=
\cap_{\xi\in A}\{\beta_{\xi}=0\}$ 
containing $y_0$ and $z$ 
has positive 
dimension. 
In particular, $\dim Y\geq2$. 

Suppose that $\dim Y=2$. 
Then $\cap_{\xi\in A}\{\beta_{\xi}=0\}$ is one-dimensional.
Let $y$ be an interior point of a segment $y_0z$. 
Since the rays $\rho_{\xi}^y$ for $\xi\in A$ are perpendicular to $y_0z$,
there can be at most two of them, 
$|A|\leq2$. 
Since we are free to choose $|A|$ with any cardinality,
we obtain a contradiction. 
Thus $\dim Y\geq3$. 

An argument analogous to the last one shows furthermore 
that the soul of $Y$ must have codimension $\geq2$. 
In particular, if $\dim Y=3$, then $\dim soul(Y)\leq1$. 

\medskip
{\em Collapse to the flat 2-halfdisc.} 
Suppose that $X$ is the flat unit halfdisc in 
$\{u\in\R^2:u_2\leq0\}$ 
centered at $x=0$.
(This case has not been treated explicitely in \cite[\S 3]{SY}. 
There, blow-up limits have been obtained under the assumption that 
$\diam(\Si_xX)<\pi$.) 

\ref{lem:blowup} implies that 
$Y$ contains a flat halfplane, but $Y\not\cong\R^2$. 
In particular, 
$Y$ splits metrically as $Y\cong\R\times W$. 
If $\dim Y=2$, 
then $W$ is a halfline and $Y$ a flat halfplane. 
If $\dim Y=3$, 
then $W$ is a noncompact Alexandrov surface with curvature $\geq0$. 
We may assume that $y_0\in 0\times W$. 
The critical points of $d(y_0,\cdot)$ lie on $0\times W$. 

If we denote by $\eta^{\pm}\in\Si_xX$ 
the directions pointing to $(\pm1,0)$, 
then for any $y\in Y$ the rays $\rho^y_{\eta^{\pm}}$ have angle $\pi$ at $y$, 
cf.\ (\ref{ineq:blexpand}), 
and their union is the line $\R\times w$ through $y$. 
Moreover, $\{\beta_{\eta^+}=\beta_{\eta^-}\}=0\times W$,
and it is the Gromov-Hausdorff limit of the bisectors 
$\{d(\si_{\eta^+}^i,\cdot)=d(\si_{\eta^-}^i,\cdot)\}$. 

For any direction $\xi\in\Si_xX$ holds 
$\angle_y(\rho^y_{\xi},\rho^y_{\eta^{\pm}})=\angle_x(\xi,\eta^{\pm})$ 
because 
$\pi=
\tangle(\rho^y_{\eta^+},\rho^y_{\xi})+\tangle(\rho^y_{\xi},\rho^y_{\eta^-})
\geq
\angle_y(\rho^y_{\eta^+},\rho^y_{\xi})+\angle_y(\rho^y_{\xi},\rho^y_{\eta^-})
\geq\angle_x(\eta^+,\xi)+\angle_x(\xi,\eta^-)=\pi$ 
also by (\ref{ineq:blexpand}). 
Let $\eta\in\Si_xX$ denote the (bisector) direction pointing to $(0,-1)$. 
Then the rays $\rho^y_{\eta}$ are orthogonal to $\rho^y_{\eta^{\pm}}$ 
and contained in layers $t\times W$. 

We choose now $\hat p_i$ as a maximum of $d(\si_{\eta}^i,\cdot)$ 
on the bisector 
$\{d(\si_{\eta^+}^i,\cdot)=d(\si_{\eta^-}^i,\cdot)\}$.
Again $\hat p_i\to x$, 
and $y_0$ is a maximum of $\beta_{\eta}$ on 
$\{\beta_{\eta^+}=\beta_{\eta^-}\}$. 
If $Y$ is a flat halfplane, 
then $\{\beta_{\eta^+}=\beta_{\eta^-}\}$ is a halfline 
and $y_0$ its endpoint. 
This is a contradiction 
because $d(y_0,\cdot)$ has critical points 
and $y_0$ cannot lie on the boundary of the halfplane $Y$. 
Thus $\dim Y\geq3$.

The next observation narrows down the possibilities for a 
3-dimensional blow-up limit $Y$.
\begin{lem}
\label{lem:1dsoul}
If $\dim Y=\dim W +1=3$ and $\dim soul(W)=1$ 
(and hence $W$ is a quotient of the flat cylinder),
then $W$ must be one-ended. 
\end{lem}
\proof
To see this, assume the contrary. 
Then $W$ splits off a line,
i.e.\ $W\cong\R\times F^1$ with a connected closed 1-orbifold $F^1$,
and we may assume that $y_0\in 0\times F^1$. 
Since $y_0$ is a maximum of $\beta_{\eta}$, 
we have $\beta_{\eta}(s,f)=-|s|$. 
There exist two unit speed rays $\rho_i:[0,\infty)\to 0\times W$ 
starting from $y_0$ in antipodal directions,
$\angle_y(\dot\rho_1(0),\dot\rho_2(0))=\pi$, 
such that $\beta_{\eta}(\rho_i(t))=-t$. 
From every point in $Y-0\times F^1$ starts a unique ray 
along which $\beta_{\eta}$ decays with slope 1, 
and thus for $s>0$ holds 
$\rho_{\eta}^{\rho_i(s)}(t)=\rho_i(s+t)$. 
It follows that there exist points 
$x_{ij}\in\{d(\si_{\eta^+}^i,\cdot)=d(\si_{\eta^-}^i,\cdot)\}$ 
such that, with respect to the rescaled metrics, 
the segments $x_{ij}\si_{\eta}^i$ 
converge to the ray $\rho_j$. 
In particular, $\de_i^{-1}d(\hat p_i,x_{ij})\to0$. 
On the other hand, without rescaling,
the two sequences of segments converge to the same segment 
$x\si_{\eta}(t_0)$. 

It follows by continuity 
that there exist points $z_{ij}\in x_{ij}\si_{\eta}^i$
such that 
$d(\si_{\eta}^i,z_{i1})=d(\si_{\eta}^i,z_{i2})$ 
and $\cangle_{\hat p_i}(z_{i1},z_{i2})=\pithird$. 
We put $l_i=d(z_{i1},z_{i2})$. 
Then $d(\hat p_i,z_{ij})\to0$ and $\de_i^{-1}d(\hat p_i,z_{ij})\to\infty$. 
Moreover, 
$\de_i^{-1}|d(\hat p_i,z_{ij})-l_i|\to0$. 
Let $m_i$ be the midpoints of segments $z_{i1}z_{i2}$. 

Triangle comparison applied to the triangles 
$\De(z_{i1},z_{i2},\si_{\eta}^i)$ yields 
$\angle_{z_{ij}}(\si_{\eta}^i,m_i)\gtrapprox\pihalf$, 
and $d(m_i,z_{ij}\si_{\eta}^i)\gtrapprox\frac{l_i}{2}$, 
whereas comparison at 
$\De(z_{i1},z_{i2},\hat p_i)$ yields 
$\angle_{z_{ij}}(\hat p_i,m_i)\gtrapprox\pithird$ 
and $d(m_i,z_{ij}\hat p_i)\gtrapprox\frac{\sqrt{3}}{4}l_i$. 
It follows that 
$l_i^{-1}d(m_i,x_{ij}\si_{\eta}^i)$ is bounded away from 0 
and $\cangle_{\hat p_i}(z_{ij},m_i)\geq\phi_0>0$ for large $i$. 
The segments $\hat p_im_i$ subconverge to a ray $\rho$ in $Y$ 
with initial point $y_0$ and 
$\angle_{y_0}(\rho_j,\rho)\geq\phi_0$. 

Comparison at the triangles 
$\De(x_{ij},\si_{\eta}^i,\si_{\eta^{\pm}}^i)$ yields that 
$\liminf\cangle_{x_{ij}}(z_{ij},\si_{\eta^{\pm}}^i)\geq\pihalf$. 
Rescaling with the factors $l_i^{-1}\to\infty$ 
and taking into account that a Gromov-Hausdorff limit splits off a line 
shows that in fact 
$\cangle_{x_{ij}}(z_{ij},\si_{\eta^{\pm}}^i)\to\pihalf$ 
and furthermore 
$\cangle_{x_{ij}}(m_i,\si_{\eta^{\pm}}^i)\to\pihalf$. 
Thus $\angle_{y_0}(\rho_{\eta^{\pm}},\rho)=\pihalf$. 
In view of $\angle_{y_0}(\rho_j,\rho)\geq\phi_0$,
this is a contradiction. 
\qed

\medskip
Combining the discussion of collapse in the various special cases, we obtain:
\begin{prop}[Blow-up limits of strictly larger dimension]
\label{prop:blowuplargerdim}
If $X=B(x,1)$ is a flat cone of dimension $\leq2$,
then the base points $\hat p_i$ can be chosen so that 
$\dim Y>\dim X$.
\end{prop}
In fact, the arguments in the special cases above only used the 
conelikeness of $X$ and the geometry of $\Si_xX$,
and hence the conclusion of \ref{prop:blowuplargerdim} holds more generally
when $X$ is conelike of dimension $\leq2$. 

%
\subsection{Strainers}
\label{sec:gradl}
\subsubsection{Position relative to the singular locus}

A Riemannian orbifold is a local Alexandrov space 
with very special singularity structure. 
The existence of a strainer in a point implies a certain regularity. 
More precisely,
for sufficiently small $\theta>0$, 
a point admits a $\theta$-straight $m$-strainer 
if and only if its link splits off a join factor 
isometric to the $(m-1)$-dimensional unit sphere, 
i.e.\ 
if and only if the singular stratum containing it has dimension $\geq m$. 
The strainer must be almost tangent to the singular stratum. 

Thus an interior point in a Riemannian 3-orbifold $O$ 
admits 3-strainers if and only if it is regular,
admits 2-strainers if and only if it is regular or a reflector boundary point,
and admits 1-strainers if and only if it is no singular vertex.
For instance, $S_{\theta,\mu,-b^2}\cap O^{(0)}=\emptyset$. 

\subsubsection{Gradient-like vector fields}
\label{sec:gradlike}
We recall that 
for a point $p\in O$ 
the distance function $d(p,\cdot)$ has directional derivatives. 
The derivative $\D_vd(p,\cdot)$ 
in the direction of a unit tangent vector $v\in\Si_xO$ 
equals $-\cos(v,D_{x,p})$ 
where $D_{x,p}\subset\Si_xO$ is the compact subset 
of the directions of all segments $xp$. 
The function $v\mapsto\D_vd(p,\cdot)$ 
on the unit tangent bundle of $O$ 
is lower semicontinuous outside $\Si_pO$, 
and there its suplevel sets are open.
(This argument also works
for the distance function from a compact subset of $O$,
e.g.\ from a component of $\D O$ if $O$ has boundary.)

For $x\neq p$ and $c\in[0,1]$ 
the subset 
$\{v\in\Si_xO:\D_vd(p,\cdot)\geq c\}$
is totally convex and has diameter $\leq2\arccos c$
(e.g.\ because it has angular distance $\geq\pi-\arccos c$ from $D_{x,p}$). 
Such a subset for $c>0$ can be nonempty only if 
$x\not\in O^{(0)}$, 
because this requires that 
$\diam(\Si_xO)\geq\pi-\arccos c>\pihalf$.
If it is nonempty for a singular point $x\in O^{(1)}\cup O^{(2)}$, 
then it contains a singular direction at $x$. 

By a standard construction using a partition of unity 
it follows that for $c\in[0,1)$ 
there exists a smooth vector field $X$ 
on the open subset 
$\{x\neq p:\exists v\in\Si_xO\hbox{ with }\D_vd(p,\cdot)>c\}$
which is tangent to the singular locus 
and satisfies 
$\D_Xd(p,\cdot)>c$. 
One calls such a vector field {\em gradient-like} for $d(p,\cdot)$. 

For distinct points $a$ and $b$ 
there exist gradient-like vector fields 
for $d(a,\cdot)$ and $d(b,\cdot)$ 
on the open set $S_{a,b}=\{\angle_{\cdot}(a,b)>\pihalf\}\subset O-\{a,b\}$. 
(Note that the function $\angle_{\cdot}(a,b)$ is lower semicontinuous
on $O-\{a,b\}$.) 
More precisely, 
for $\phi\in(0,\pihalf]$ 
exists a gradient-like vector field $X$ for $d(a,\cdot)$
on $\{\angle_{\cdot}(a,b)>\pi-\phi\}$ 
with $\angle_{\cdot}(a,X)>\pi-\phi$. 
Such a vector field satisfies 
$\angle_{\cdot}(b,-X)
\geq\angle_{\cdot}(a,b)-\angle_{\cdot}(a,-X)>\pi-2\phi$ 
and $\angle_{\cdot}(b,X)<2\phi$. 
(Note that $-X$ is defined, 
since $X$ is tangent to the singular locus.) 
Thus, if $\phi\leq\piquart$ 
then $\D_Xd(b,\cdot)<0$ and $\D_Xf_{a,b}>0$, 
i.e.\ $X$ is gradient-like for $f_{a,b}$. 
 

If $(a,b)$ is a $\theta$-straight 1-strainer at $p$ 
for sufficiently small $\theta$, 
then its cross section $\Si_{x;a,b}$ is near $p$ 
a topological 2-suborbifold, 
because a gradient-like vector field for $f_{a,b}$ 
has local cross sections  
through $p$ which are smooth 2-suborbifolds, 
and any two local cross sections can be isotoped to each other 
using the flow. 

%
%
%

\subsubsection{Local bilipschitz charts and fibrations by cross sections}
\label{sec:chafib}

Suppose that 
$(a_1,b_1,a_2,b_2,a_3,b_3)$ is a $<\theta$-straight 3-strainer 
at a regular point $x\in O$ 
for some small $\theta>0$. 
Then there exist unit vectors $v_1,v_2,v_3\in\Si_xO$ 
with $\angle_x(a_i,v_i)>\pi-\theta$. 
It follows that 
$\angle_x(b_i,-v_i)>\pi-2\theta$, 
$|\angle_x(a_i,v_j)-\pihalf|,|\angle_x(b_i,v_j)-\pihalf|<2\theta$ 
and 
$|\angle_x(v_i,v_j)-\pihalf|<3\theta$ for $i\neq j$. 
Let $X_i$ be arbitrary commuting smooth vector fields near $x$ 
with $X_i(x)=v_i$. 
By continuity, on a sufficiently small neighborhood of $x$ 
they have length $\approx1$ 
and satisfy the same angle inequalities
$\angle_{\cdot}(a_i,X_i)>\pi-\theta$ 
and their implications. 
Thus, they are almost orthogonal, 
$|\angle(X_i,X_j)-\pihalf|<3\theta$ for $i\neq j$, 
and gradient-like for the functions $f_{a_i,b_i}$ 
associated to the 1-substrainers, 
$\D_{X_i}f_{a_i,b_i}>\half(\cos\phi+\cos2\phi)>1-c\theta^2$, 
and 
$|\D_{X_j}f_{a_i,b_i}|<\sin2\theta<2\theta$ for $i\neq j$. 
The $X_i$ are the coordinate vector fields for some local coordinates, 
and it follows that $(f_{a_1,b_1},f_{a_2,b_2},f_{a_3,b_3})$ 
restricts to a bilipschitz homeomorphism 
from a neighborhood of $x$ onto an open subset of $\R^3$. 
(Compare the discussion in
\cite[\S 11.8]{BGP} und \cite[2.4.1]{MT2}.)

An analogous argument can be carried out for a 
$2\half$-strainer $(a_1,b_1,a_2,b_2,a_3)$ 
at a reflector boundary point $x\in O^{(2)}$. 
Then the 2-substrainer $(a_1,b_1,a_2,b_2)$ 
is almost tangent to the reflector boundary. 
One uses an orbifold chart at $x$, 
lifts the functions $f_{a_1,b_1},f_{a_2,b_2}$ and constructs $f_{a_3,b_3}$ 
by lifting $a_3$ to a 1-strainer in the chart. 
Furthermore, one adapts the smooth vector fields $X_i$ to the 
reflection $\iota$ on the chart, 
i.e.\ constructs them so that 
$\iota^*X_1=X_1$,  $\iota^*X_2=X_2$ 
and $\iota^*X_3=-X_3$. 
One obtains a local bilipschitz homeomorphism 
to the 3-dimensional halfspace with reflector boundary. 

Consider now a $<\theta$-straight 2-strainer $(a_1,b_1,a_2,b_2)$ 
at a point $x$. 
That its cross section $\Si_{x;a_1,b_1,a_2,b_2}$ 
is near $x$ a bilipschitz 1-suborbifold,
can be seen as follows. 
If $x$ is a regular point, 
then one can choose an $\approx\theta$-straight 1-strainer 
$(a',b')$ contained in $\Si_{x;a_1,b_1,a_2,b_2}$ and, 
with respect to the local bilipschitz coordinates near $x$ 
provided by the 3-strainer $(a_1,b_1,a_2,b_2,a',b')$, 
the cross section $\Si_{x;a_1,b_1,a_2,b_2}$ 
is a coordinate line. 
If $x$ is singular, then it must be a reflector boundary point. 
One chooses a point $a'\in\Si_{x;a_1,b_1,a_2,b_2}$ near $x$ 
and works with the bilipschitz coordinates provided by the 
$2\half$-strainer $(a_1,b_1,a_2,b_2,a')$. 

Suppose that $C$ is a compact connected component of $\Si_{x;a_1,b_1,a_2,b_2}$ 
such that $(a_1,b_1,a_2,b_2)$ is a 
$c\theta$-straight 2-strainer at all points of $C$. 
Then $C$ is a closed 1-suborbifold 
and hence homeomorphic to $S^1$ or the mirrored interval $I^1$. 
The map $(f_{a_1,b_1},f_{a_2,b_2})$ 
yields a product fibration of a neighborhood of $C$
by cross sections of the 2-strainer. 
This can be seen as follows 
using the bilipschitz coordinates 
near the points $y\in C$. 
The distance functions $f_{a_1,b_1},f_{a_2,b_2}$ given by the 2-strainer 
and the auxiliary distance function $f_{a_3,b_3}$ near $y$ 
are normalized so that $f_{a_i,b_i}(y)=0$. 
For sufficiently small $\eps_i>0$ 
the map $(f_{a_1,b_1},f_{a_2,b_2})$ 
yields near $y$ a product fibration of the box 
$\{|f_{a_i,b_i}|\leq\eps_i\,\forall i\}$ 
over the rectangle $[-\eps_1,\eps_1]\times[-\eps_2,\eps_2]$. 
By covering $C$ with finitely many such boxes 
one obtains the fibration of a neighborhood. 


In the above discussion, one or both of the functions 
$f_{a_1,b_1}$ and $f_{a_2,b_2}$ 
can be replaced by $d(a_i,\cdot)$, 
because their directional derivatives differ only slightly. 
(Recall that $f_{a_i,b_i}-d(a_i,\cdot)$ is $c'\theta$-Lipschitz 
in the region where $(a_i,b_i)$ is a $c\theta$-straight 1-strainer,
cf.\ section~\ref{sec:strest}). 

\subsection{A decomposition according to the coarse stratification}
\label{sec:decomposrs}

We start by formulating the collapse assumption 
on our orbifolds $(O,g)$ 
needed in this section 
and the quantities involved in it. 

The parameter $\theta$ 
(straightness of 1-strainers) 
is required to be small, 
$\theta\in(0,\theta_0]$,
where $\theta_0$ is sufficiently small
for the arguments
in sections \ref{sec:1stralex} and \ref{sec:rough2} to apply
(cf. the discussion at the beginning of section \ref{sec:1stralex}).
The upper bound for $\theta$ will be decreased several times 
during our later arguments.
The parameter $\mu$ 
(accuracy of conical approximation) 
needs to be sufficiently small so that 
the conclusions of \ref{lem:nhump} 
regarding the existence of $\theta$-straight 1-strainers apply, 
$\mu\leq\mu_0(\theta)$ 
with the constant 
$\mu_0(\theta)$ from there.
The parameter $\mu$ determines 
(together with the fixed parameter $\si$) 
via \ref{prop:locapproxconeorbi} 
the bound $s_1(\si,\mu)$ 
and the scale 
$\hat s_{\mu,-b^2}:=\eleventh s_1(\si,\mu)\rho_{-b^2}$. 
Conical approximation in all points $x\in O$ 
on scales 
$s(x)\in[s_1(\si,\mu)\rho_{-b^2}(x),\si\rho_{-b^2}(x)]$ 
holds if $(O,g)$ is $(v(\si,\mu),-b^2)$-collapsed 
with the constant $v(\si,\mu)>0$ from \ref{prop:locapproxconeorbi}. 

In order to make the results on edgy points in section~\ref{sec:rough2} 
available,
we need to rule out the existence of 
$\bar\theta_{2\half}$-straight $2\half$-strainers 
with length $\la(\theta)\hat s_{\mu,-b^2}(x)$ at all points $x$ 
for a constant $\la(\theta)$ which is sufficiently small 
so that \ref{lem:auniqext}, \ref{lem:reledgynew} 
and \ref{lem:parcrossedg} hold. 
This is achieved by requiring $(O,g)$ to be 
$(v(\theta,\mu),-b^2)$-collapsed 
for a suitable small constant $v(\theta,\mu)>0$, 
cf.\ \ref{lem:collstrain}. 
We make it so small 
that it is smaller than the constant $v(\si,\mu)$ mentioned above.

Furthermore, 
in order to obtain sufficient collapse on the scale 
$\theta^4\hat s_{\mu,-b^2}(x)$ for small $\theta$, 
we ask that 
$v(\theta,\mu)\leq(\theta^5s_1(\si,\mu))^3$. 
We will also assume that $(O,g)$ has 
$(v(\theta,\mu),s_0,K)$-curvature control below scale $\rho_{-b^2}$. 

\subsubsection{The 2-strained region}
\label{sec:2str}

We define the {\em 2-strained region}
$R_{\theta,\mu,-b^2}\subset O$ 
as the open set consisting of all points $x$ which admit 
$<C_1\theta$-straight 2-strainers of length 
$>\theta^4\hat s_{\mu,-b^2}(x)$, 
where $C_1$ is the constant from \ref{lem:exstrainnew}. 
Note that 
$R_{\theta,\mu,-b^2}\cap O^{(sing)}\subset O^{(2)}$. 

We will show that for sufficiently small $\theta$ 
the 2-strained region admits metrically an almost product fibration
with short fibers almost orthogonal to the 2-strainers. 
(The smallness of the parameter $\mu$ is not important at this point, 
because we are not yet using conical approximation 
from \ref{prop:locapproxconeorbi}.) 

We begin with a local approximation result:
\begin{prop}
\label{cor:locmod2str}
For $\eps>0$ exists 
$\theta_1=\theta_1(\eps)>0$ 
such that: 

Let $\theta\leq\theta_1$ and $\mu\leq\mu_0(\theta)$. 
If $(O,g)$ is $(v(\theta,\mu),-b^2)$-collapsed 
with $(v(\theta,\mu),s_0,K)$-curvature control below scale $\rho_{-b^2}$, 
and if $(a_1,b_1,a_2,b_2)$ is a 
$C_1\theta$-straight $\theta^4\hat s_{\mu,-b^2}(x)$-long 
2-strainer at $x$, 
then $\diam(\Si^o_{x;a_1,b_1,a_2,b_2})^{-1}\cdot (O,x)$
is $\eps$-close in the pointed ${\mathcal C}^5$-topology 
to the product $\R^2\times F^1$ of the euclidean plane 
with a connected closed 1-orbifold of diameter 1. 
\end{prop}
\proof
Let $-b_i^2\in[-1,0)$ 
and $\theta_i,\mu_i>0$ such that 
$\theta_i\to0$ and $\mu_i\leq\mu_0(\theta_i)$. 
Suppose that the orbifolds $(O_i,g_i)$ are 
$(v(\theta_i,\mu_i),-b_i^2)$-collapsed 
with $(v(\theta_i,\mu_i),s_0,K)$-curvature control 
below the scales $\rho_{-b_i^2}$,
and that the 
$(a^i_1,b^i_1,a^i_2,b^i_2)$ are $C_1\theta_i$-straight 
$\theta_i^4\hat s_{\mu_i,-b_i^2}(x_i)$-long 2-strainers 
at the points $x_i$. 
We have to show that the conlusion of the proposition holds for large $i$. 

Consider the cross sections 
$\Si^o_i=\Si^o_{(x_i;a^i_1,b^i_1,a^i_2,b^i_2)}$. 
For any point $x_i\neq y_i\in\Si^o_i$, 
the $2\half$-strainer $(a^i_1,b^i_1,a^i_2,b^i_2,y_i)$ at $x_i$ 
is $\bar\theta_{2\half}$-straight for large $i$. 
Since $v(\theta,\mu)\leq(\theta^5s_1(\si,\mu))^3$, 
it follows that 
$\la_i:=(\theta_i^4\hat s_{\mu_i,-b_i^2}(x_i))^{-1}\diam(\Si^o_i)\to0$, 
cf.\ \ref{lem:collstrain}. 
The rescaled pointed orbifolds 
$(\la_i\theta_i^4\hat s_{\mu_i,-b_i^2}(x_i))^{-1}\cdot(O_i,x_i)$
Gromov-Hausdorff subconverge to 
a pointed Alexandrov space $(X,x_0)$ 
with dimension $\leq3$ and curvature $\geq0$. 
Moreover, the broken segments $a^i_1x_ib^i_1$ and $a^i_2x_ib^i_2$
(sub)converge to two perpendicular lines through $x_0$, 
and hence $X$ splits metrically as a product $\R^2\times\Si$. 

To see that the rescaled cross sections 
$(\la_i\theta_i^4\hat s_{\mu_i,-b_i^2}(x_i))^{-1}\cdot\Si^o_i$ 
subconverge to $\Si$,
we note that 
$(\la_i\theta_i^4\hat s_{\mu_i,-b_i^2}(x_i))^{-1}
(d(a^i_j,\cdot)-d(a^i_j,x_i))$ 
and $(\la_i\theta_i^4\hat s_{\mu_i,-b_i^2}(x_i))^{-1}
(d(b^i_j,\cdot)-d(b^i_j,x_i))$ 
subconverge (due to Arzela-Ascoli) 
to concave 1-Lipschitz functions $\al_j$ and $\beta_j$ 
with $\al_j(x_0)=\beta_j(x_0)=0$. 
The concavity of the sums $\al_j+\beta_j$ 
implies together with the triangle inequality 
that the functions $\al_j$ and $\beta_j$ are constant on fibers 
$pt\times\Si$. 
Since for any sequence points $y_i\in\Si^o_i$ holds that 
$(\la_i\theta_i^4\hat s_{\mu_i,-b_i^2}(x_i))^{-1}
(d(a^i_j,y_i)-d(a^i_j,x_i)),
(\la_i\theta_i^4\hat s_{\mu_i,-b_i^2}(x_i))^{-1}
(d(b^i_j,y_i)-d(b^i_j,x_i))\to0$, 
compare (\ref{eq:crosssecperp}), 
it follows that $\Si^o_i\to\Si$. 
Therefore $\Si$ is a compact Alexandrov space 
with diameter 1 and dimension 1. 

Since in the blow-up limit 
$(\la_i\theta_i^4\hat s_{\mu_i,-b_i^2}(x_i))^{-1}\cdot(O_i,x_i)\to(X,x_0)$
we have no dimension drop, 
after passing to another subsequence, 
the convergence can be improved to ${\mathcal C}^5$-smooth convergence
and $F^1=\Si$ is a connected closed 1-orbifold,
cf.\ \ref{lem:smconv}. 
\qed

\medskip
Note that in the situation of the proposition, 
$\diam(\Si_{x;a_1,b_1,a_2,b_2})<<\theta^4\hat s_{\mu,-b^2}(x)$.
Moreover,
for $x\in R_{\theta,\mu,-b^2}$
the cross sections $\Si^o_{x;a_1,b_1,a_2,b_2}$ of different 
$C_1\theta$-straight $\theta^4\hat s_{\mu,-b^2}(x)$-long
2-strainers $(a_1,b_1,a_2,b_2)$ at $x$ 
have Hausdorff distance $<<\diam(\Si^o_{x;a_1,b_1,a_2,b_2})$
and almost equal diameters,
and we define the {\em width} $w(x)$ at $x$ 
as the infimum of these diameters. 

The fiber direction of a local approximation as in \ref{cor:locmod2str}
yields a smooth line field 
which is {\em almost vertical} in the sense that 
it is perpendicular to the stratum $O^{(2)}$ 
and almost perpendicular to sufficiently long segments. 
Any two such local line fields almost agree on the overlaps 
of their domains of definition, 
and using a partition of unity 
we can combine such local line fields 
to a global almost vertical line field 
$L=L_{\theta,\mu,-b^2}$ on $R_{\theta,\mu,-b^2}$. 
More precisely, 
for (small) $\nu>0$ 
exists $\theta_1=\theta_1(\nu)>0$ 
such that the following holds: 
If $\theta\leq\theta_1$, $\mu\leq\mu_0(\theta)$ and 
if $(O,g)$ is $(v(\theta,\mu),-b^2)$-collapsed with 
$(v(\theta,\mu),s_0,K)$-curvature control below scale $\rho_{-b^2}$, 
then for $x\in R_{\theta,\mu,-b^2}$ 
the line $L(x)$ has angle $>\pihalf-\nu$ 
with any segment of length 
$10\nu^{-1}w(x)$ initiating in $x$,  
and in particular with any segment of length 
$\theta^4\hat s_{\mu,-b^2}(x)>>w(x)$. 
In fact,
a line field $L_{\theta,\mu,-b^2}$ with these properties 
can be constructed on the slightly larger open set 
$\hat R_{\theta,\mu,-b^2}:=
\cup_{x\in R_{\theta,\mu,-b^2}}B(x,\frac{1}{\nu}w(x))$.
The reflector boundary $O^{(2)}$ is, where it meets $R_{\theta,\mu,-b^2}$, 
{\em almost horizontal} in the sense that it is almost tangent 
to sufficiently long segments and, in particular, 
to sufficiently long 2-strainers. 

The trajectories of $L$ starting in a point $x\in R_{\theta,\mu,-b^2}$
move almost orthogonally to sufficiently long segments 
and therefore remain close 
to the cross section $\Si^o_{x;a_1,b_1,a_2,b_2}$ 
of a suitable 2-strainer as above 
for length at least $>>w(x)$. 

If $x\in O^{(2)}\cap R_{\theta,\mu,-b^2}$, 
then the trajectory is orthogonal to $O^{(2)}$ in $x$ 
and reaches $O^{(2)}$ again after length $\approx w(x)$. 
More generally, all trajectories intersecting $B(x,10w(x))$ 
have length $\approx w(x)$ and connect reflector boundary points. 
The construction of an almost product fibration 
by mirrored intervals 
close to $O^{(2)}\cap R_{\theta,\mu,-b^2}$ 
is therefore immediate. 

Away from the reflector boundary, 
the $L$-trajectories starting in points $x\in R_{\theta,\mu,-b^2}$ 
almost close up after length $\approx 2w(x)$. 
The question whether $L$ can be globally perturbed to an integrable line field 
with closed trajectories of lengths $\approx 2w(x)$
has been treated in \cite[\S 4.2]{MT2} in a very similar setting. 
The discussion there uses only the control on finitely many derivatives 
of the curvature tensor and goes through without change
in the situation considered here.
One obtains the following result 
which can be considered as a version of a special case of 
Yamaguchi's Fibration Theorem \cite{Yamaguchi} 
``without a priori given base''. 
\begin{prop}[Almost vertical fibration of the 2-strained region,
cf.\ {\cite[Prop.\ 4.4]{MT2}}]
\label{prop:fib}
For $\nu>0$ exists $\theta_1=\theta_1(\nu)>0$ 
such that:

If $\theta\leq\theta_1$, $\mu\leq\mu_0(\theta)$ and 
if $(O,g)$ is $(v(\theta,\mu),-b^2)$-collapsed with 
$(v(\theta,\mu),s_0,K)$-curvature control below scale $\rho_{-b^2}$, 
then there exists an open subset $U$, 
$R_{\theta,\mu,-b^2}\subseteq U\subseteq\hat R_{\theta,\mu,-b^2}$, 
such that every connected component of $U$ is the total space
of a smooth orbifold fibration 
with fiber 
$S^1$ or the mirrored interval $I^1$, 
and all fibers have angle $<\nu$ 
with the almost vertical line field $L_{\theta,\mu,-b^2}$. 
\end{prop}

The local fibrations provided by the product approximations 
in \ref{cor:locmod2str} have only a finite degree of regularity,
but the global fibration of $U$ obtained by interpolating 
these local fibrations can be smoothed. 
The smoothness will however not be important to us. 

From now on,
we fix some small positive value of $\nu$
(e.g.\ $\nu = \frac{1}{2010}$)
and set $\theta_1 = \theta_1(\nu)$.
Moreover, whenever an orbifold is sufficiently collapsed,
we implicitely fix a fibration as in \ref{prop:fib}.

%
\subsubsection{Edges}
\label{sec:edg}

\begin{defn}
\label{def:edgyorbi}
We define a point $x\in S_{\theta,\mu,-b^2}$ to be 
$(\theta,\mu,-b^2)$-{\em edgy} 
relative to an equilateral $<\theta$-straight 
1-strainer $(a,b)$ at $x$
with length in $(\hat{s}_{\mu,-b^2}(x),\frac{3}{2}\hat{s}_{\mu,-b^2}(x))$
if $\diam(\Si^o_{x;a,b}) > \theta^\frac{5}{2}\hat s_{\mu,-b^2}(x)$ and 
if $\Si_{x;a,b}$ contains no $\frac{\pi}{2}$-straight 
$\theta^4\hat s_{\mu,-b^2}(x)$-long 1-strainers at $x$, 
compare \ref{def:edgy}. 
\end{defn}


Let us briefly discuss
the correspondence between the two
definitions \ref{def:edgy} and \ref{def:edgyorbi}.
First, we recall from our discussion
in section \ref{sec:strcurvbound}
that a $(\theta,\mu,-b^2)$-edgy point $x$
relative to a $<\theta$-straight 1-strainer $(a,b)$
with length in $(\hat{s}_{\mu,-b^2}(x),\frac{3}{2}\hat{s}_{\mu,-b^2}(x))$
is also $\theta$-edgy in the space
$(\hat{s}_{\mu,-b^2}(x))^{-1} \cdot B(x,\rho_{-b^2}(x))$
of curvature $\ge -1$
in the sense of Definition \ref{def:edgy}:
In the rescaled space,
the strainer $(a,b)$ is $<2\theta$-straight
of length $> (1-\theta)$
and the cross section
$\Si_{x;a,b}$ cannot contain any
$\pihalf$-straight 1-strainer of length $\theta^4$.

Conversely, suppose that an equilateral 1-strainer
$(a,b)$ is
$<\theta$-straight at $x$ 
with length in $(\hat{s}_{\mu,-b^2}(x),\frac{3}{2}\hat{s}_{\mu,-b^2}(x))$
and that in the rescaled space
$\hat{s}_{\mu,-b^2}(x)^{-1} \cdot B(x,\rho_{-b^2}(x))$
$x$ is $\theta$-edgy relative to $(a,b)$
(again in the sense of Definition \ref{def:edgy}).
Then Remark \ref{rem:edgyonscales}
together with the discussion at the end
of \ref{sec:strcurvbound}
implies that $x$ is also $(\theta,\mu,-b^2)$-edgy.
For suppose that $\Si_{x;a,b}$ contains
a $\frac{\pi}{2}$-straight
$\theta^4\hat s_{\mu,-b^2}(x)$-long
1-strainer at $x$.
Then after rescaling this strainer
is still $<\frac{3\pi}{4}$-straight
and hence in fact $<C_1\theta$-straight
which is a contradiction.

By construction, for any $x \in O$ there are no
$\bar\theta_{2\half}$-straight $2\half$-strainers
of length $\la(\theta)$
in the rescaled space
$\hat{s}_{\mu,-b^2}(x)^{-1} B(x,\half\rho_{-b^2}(x))$
Moreover, estimate \ref{eq:osccurvsc} implies
that on the ball $B(x,\theta\hat{s}_{\mu,-b^2}(x))$
we have the estimate
$\hat{s}_{\mu,-b^2}/\hat{s}_{\mu,-b^2}(x) \in (1-\theta,1+\theta)$.

This allows us to generalize the above arguments
to the following results:
If $y \in B(\theta\hat{s}_{\mu,-b^2}(x))$
is $(\theta,\mu,-b^2)$-edgy relative to a 1-strainer $(a,b)$
then it is $\theta$-weakly in the rescaled space
$\hat{s}_{\mu,-b^2}(x)^{-1} \cdot B(x,\rho_{-b^2}(x))$
of curvature $\ge -1$.
Conversely, suppose that $y \in B(\theta\hat{s}_{\mu,-b^2}(x))$
admits an equilateral $<\theta$-straight
1-strainer $(a,b)$ with length in $(\hat{s}_{\mu,-b^2}(y),\frac{3}{2}\hat{s}_{\mu,-b^2}(y))$
such that $y$ is $\theta$-strongly edgy
relative to $(a,b)$ in $\hat{s}_{\mu,-b^2}(x)^{-1} \cdot B(x,\rho_{-b^2}(x))$.
Then $y$ is also $(\theta,\mu,-b^2)$-edgy.

Hence our results from section \ref{sec:edgalex}
can be applied to our present situation and used to control
e.g.\ the relative position of edgy points.

\medskip

We will globally construct tubes along the ``coarse edges''.
Our previous discussion applies 
provided that $\theta$ and $\mu$ are sufficiently small
(i.e.\ $\theta \le \theta_1$ and $\mu \le \mu_0(\theta)$),
and that $(O,g)$ is $(v(\theta,\mu),-b^2)$-collapsed
with $(v(\theta,\mu),s_0,K)$-curvature control
below scale $\rho_{-b^2}$.

For every $(\theta,\mu,-b^2)$-edgy point $x\in O$, 
let $(a_x,b_x)$ be some equilateral $<\theta$-straight 1-strainer 
with length in $(\hat s_{\mu,-b^2}(x),\frac{3}{2}\hat s_{\mu,-b^2}(x))$
relative to which $x$ is edgy.
(Since 1-strainers at $x$ are almost unique by \ref{lem:1stredgynew},
it does not matter which one we use.) 
Associated to it is the truncated cross section 
$\check\Si_x=\Si_{x;a_x,b_x}\cap 
B(x,\half\theta^3\hat s_{\mu,-b^2}(x))$. 

For $j=1,2,3$ we consider the fibers 
$\ga_x^j=
\Si_{x;a_x,b_x}\cap B(x,\frac{j}{8}\theta^3\hat s_{\mu,-b^2}(x))$ 
of the partial (topological) product fibration of $\Si_{x;a_x,b_x}$ 
induced by $d(x,\cdot)$.
(Compare the discussion in section \ref{sec:chafib}.) 
The $\ga_x^j$ and also their 
$\frac{1}{16}\theta^3\hat s_{\mu,-b^2}(x)$-neighborhoods
are contained in the 2-strained region $R_{\theta,\mu,-b^2}$, 
and they are almost vertical in the sense that they are isotopic 
to a fiber of the fibration given by \ref{prop:fib} 
by a small isotopy, say, supported on the 
$\theta^4\hat s_{\mu,-b^2}(x)$-neighborhood of $\ga_x^j$.

An almost unique 1-strainer at a $(\theta,\mu,-b^2)$-edgy point $x$
locally defines a vector field on the ball 
$B(x,\theta \hat{s}_{\mu,-b^2})$
as in section \ref{sec:gradlike}.
By interpolating these fields by a partition of unity,
we obtain a smooth vector field $L = L_{edge}$ tangent to the singular locus 
whose open domain of definition contains the 
balls $B(x,\theta \hat s_{\mu,-b^2}(x))$ 
around the $(\theta,\mu,-b^2)$-edgy points $x$, 
and which is on these balls almost parallel to the 1-strainers $(a_x,b_x)$,
meaning that 
$\angle_{\cdot}(L_{edge},a_x),\angle_{\cdot}(L_{edge},b_x)
\not\in[c\theta^{\half},\pi-c\theta^{\half}]$
on $B(x,\theta \hat s_{\mu,-b^2}(x))$, 
see \ref{lem:auniqext}. 
In particular, we have 
$|\angle_{\cdot}(L_{edge},x)-\pihalf|<c'\theta^{\half}$ 
on the 
$\seventh\theta^3 \hat s_{\mu,-b^2}(x)$-neighborhood of $\ga_x^2$. 

We choose a maximal subfamily of pairs $(x,\check\Si_x)$ 
such that the corresponding subset $\eps$ of 
$(\theta,\mu,-b^2)$-edgy points $x$ 
is separated in the sense that for any two distinct points 
$x_1,x_2\in\eps$ holds 
$d(x_1,x_2)>\theta^{\frac{10}{3}}\hat s_{\mu,-b^2}(x_1)$. 
By \ref{lem:parcrossedg}, 
the $\check\Si_x$ for $x\in\eps$ are pairwise disjoint. 

We call $x_1,x_2\in\eps$ {\em adjacent}, if 
$d(x_1,x_2)<4\theta^{\frac{10}{3}}\hat s_{\mu,-b^2}(x_1)$
and if the arc in the $L_{edge}$-trajectory (leaf) 
connecting $x_1$ to its intersection point with $\check\Si_{x_2}$
does not meet the other cross sections $(x,\check\Si_x)$ 
for $x\in\eps-\{x_1,x_2\}$. 
By \ref{lem:reledgynew},
we have $d(x_1,x_2)\approx\theta^{\frac{10}{3}}\hat s_{\mu,-b^2}(x_1)$
unless points on the segment $x_1x_2$
have no $\theta$-straight $\hat{s}_{\mu,-b^2}$-long 1-strainers
or have such strainers with cross sections
of diameter $\approx \theta^{\frac{5}{2}}$.
The relation of adjacency generates an equivalence relation on $\eps$, 
and we call an equivalence class a 
{\em chain} of $(\theta,\mu,-b^2)$-edgy points. 
Each chain can be given a linear or cyclic order. 

Consider two adjacent edgy points $x_1,x_2\in\eps$. 
The cross sections $\check\Si_{x_i}$ have Hausdorff distance 
$<\theta^{\frac{99}{30}}\hat s_{\mu,-b^2}(x_1)$
by \ref{lem:parcrossedg}. 
Using the integral curves of the line field $L_{edge}$, 
we flow the 1-orbifolds $\ga_{x_1}^j$ 
from $\check\Si_{x_1}$ into $\check\Si_{x_2}$. 
Their images ${\ga'}_{x_1,x_2}^j$ in $\check\Si_{x_2}$ 
are $<\theta^{\frac{9}{20}+\frac{99}{30}}\hat s_{\mu,-b^2}(x_1)
<\theta^{\frac{7}{2}}\hat s_{\mu,-b^2}(x_1)$-close 
to the $\ga_{x_2}^j$, 
compare the discussion of the maps $\Phi^{a,b}_t$ 
in section~\ref{sec:strest}. 
Moreover, 
since the isotopy of $\ga_{x_1}^j$ to ${\ga'}_{x_1,x_2}^j$ 
takes place inside a 
$\theta^{\frac{98}{30}}\hat s_{\mu,-b^2}(x_1)$-ball 
contained in $R_{\theta,\mu,-b^2}$, 
${\ga'}_{x_1,x_2}^j$ is isotopic by a small isotopy 
to the almost vertical fibers of the fibration of $R_{\theta,\mu,-b^2}$, 
and therefore it must inside $\check\Si_{x_2}$ be homotopic to $\ga_{x_2}^j$, 
and hence isotopic by a small isotopy. 
Thus the trace of the isotopy of $\ga_{x_1}^j$ to ${\ga'}_{x_1,x_2}^j$ 
can be adjusted 
(by a small isotopy supported near $\ga_{x_2}^j$) 
to a 2-suborbifold $S_{x_1,x_2}^j$
homeomorphic to $\cong\ga_{x_1}^j\times[0,1]$, 
contained in 
$A(x;(\frac{j}{8}-\hunth)\theta^3\hat s_{\mu,-b^2}(x),
(\frac{j}{8}+\hunth)\theta^3\hat s_{\mu,-b^2}(x))$,
lying between $\check\Si_{x_1}$ and $\check\Si_{x_2}$, 
and with boundary $\ga_{x_1}^j\cup\ga_{x_2}^j$. 
By concatenating the $S_{x_1,x_2}^j$, 
we obtain three disjoint embedded 2-suborbifolds $S^j$ 
following along the chains of edgy points. 
We call the region contained between $S^1$ 
and the two final cross section of the chain
a {\em tube} along the coarse edge.

Note that again by \ref{lem:reledgynew},
a chain can only end inside a $(\theta,\mu,-b^2)$-hump
or if the cross sections to $C_0 \theta$-straight $\hat{s}_{\mu,-b^2}$-long 1-strainers
have diameter $\approx \theta^{\frac{5}{2}}$.
It is also possible that an edgy point has no adjacent edgy points.
We simply discard such isolated cross sections.

\medskip

The simplest {\em interface} between chains of edgy points
(the ``coarse edge'') 
and the rest of $O$ arises for a 
{\em cyclic chain} $\kappa\subseteq\eps$ of $(\theta,\mu,-b^2)$-edgy points. 
The parts $S^j_{\kappa}\subset S^j$ corresponding to $\kappa$ 
are then closed 2-suborbifolds. 
Let $T^j_{\kappa}$ denote the tube 
containing $\kappa$ and bounded by $S^j_{\kappa}$,
and let $A^{i,j}_{\kappa}:=T^j_{\kappa}-int(T^i_{\kappa})$ 
for $1\leq i<j\leq3$. 
The $T^j_{\kappa}$ and $A^{i,j}_{\kappa}$ 
are compact 3-suborbifolds  
compatibly fibering over the circle, 
and $A^{1,3}_{\kappa}\subset R_{\theta,\mu,-b^2}$. 
According to \ref{prop:topcredg}, 
the fiber of $T^j_{\kappa}$ is a compact 2-orbifold 
with Euler characteristic $\chi\geq0$ and one boundary component. 
Note that $S^2_{\kappa}$ separates $S^1_{\kappa}$ and $S^3_{\kappa}$. 

We wish to replace $S^2_{\kappa}$ by a 2-suborbifold 
which is {\em vertically saturated}, 
i.e.\ saturated with respect to the fibration of $R_{\theta,\mu,-b^2}$, 
cf.\ \ref{prop:fib}. 
To do so, we take a vertically saturated compact connected 3-suborbifold $W$ 
such that $S^2_{\kappa}\subset int(W)$ and $W\subset int(A^{1,3}_{\kappa})$. 
Then $W$ separates $S^1_{\kappa}$ and $S^3_{\kappa}$, 
and according to Lemma \ref{lem:separates} below, 
one of the boundary components of $W$ 
separates $S^1_{\kappa}$ and $S^3_{\kappa}$. 
We denote this boundary component by $S^{2,v}_{\kappa}$.

Then as a consequence of Lemma \ref{lem:wald}, 
$S^{2,v}_{\kappa}$
is isotopic to the $S^j_{\kappa}$. 
In other words,
we can isotope $S^2_{\kappa}$ by an isotopy supported in 
$int(A^{1,3}_{\kappa})$ so that it becomes vertically saturated 
and the fibrations on it 
induced by $R_{\theta,\mu,-b^2}$ and $T^2_{\kappa}$ match. 

\medskip

\begin{lem}
\label{lem:separates}
Let $\Si$ be a connected closed 2-orbifold without singular points
and let $W\subset\Si\times[0,1]$ be a compact connected 3-suborbifold
disjoint from $\Si\times0$ and $\Si\times1$ and separating them.
Then some component of $\D W$ separates them also.
\end{lem}
\proof

For the purpose of this lemma,
we consider $\Si$ and $\Si'$ as compact manifolds, possibly with boundary.

Based on the existence and uniqueness of smooth structures on 3-manifolds
and the uniqueness up to isotopy of smooth structures on 2-manifolds,
we know that there exists a smooth structure on $\Si\times[0,1]$
with respect to which the embedded topological 2-submanifold $\D W$ is a smooth submanifold.
(Cut along $\D W$, put a smooth structure and glue again after adjusting the induced smooth structures on the boundaries by an isotopy in a collar.
See e.g.\ \cite{Munkres}, \cite{Whitehead}, \cite{Epstein}.)

Moreover, given a smooth structure on $\Si$ (and hence on $\Si\times[0,1]$),
there exists a homeomorphism of $\Si\times[0,1]$
carrying  the embedded topological 2-submanifolds $\D W$ to a smooth submanifold.
Hence we may work without loss of generality in the smooth category.

Let $V$ denote the component of $\Si\times(0,1)-int(W)$ containing $\Si\times0$.
Then $\Si'=V\cap W$ separates $\Si\times0$ and $\Si\times1$,
and we have to show that it is connected.
Suppose the contrary,  i.e.\ that it decomposes as the disjoint union
$\Si'=\Si'_1\cup\Si'_2$  of closed 2-submanifolds.
Then there exists an embedded circle $\ga$ in $\Si\times(0,1)$
which intersects $\Si'_1$ once transversally.
This is absurd because $\ga$ can be homotoped into $\Si\times0$.
\qed 

\medskip

The following fact is for tori a simple special case of a result of Waldhausen \cite[2.8]{Waldhausen}.
The arguments in the non-orientable and orbifold cases are similar.

\begin{lem} \label{lem:wald}
(i) Let $\Si$ and $\Si'$ be closed surfaces,
each of which is homeomorphic to the 2-torus $T^2$  or to the Klein bottle $K^2$.
Suppose that $\Si'\subset\Si\times[0,1]$ is embedded
so that it is disjoint from $\Si\times0$ and $\Si\times1$  and separates them.
Then $\Si'$ is isotopic to $\Si\times0$ and $\Si\times1$.

(ii) The same conclusion holds  if $\Si$ and $\Si'$ are closed 2-orbifolds,
each of which is homeomorphic to the annulus $Ann^2$  or the M\"obius strip $M\ddot ob^2$
with reflector boundary.
\end{lem}

\proof
Again we can assume without loss of generality
that $\Si'$ is a smooth surface, respectively, suborbifold.

(i) Cut open $\Si\times[0,1]$ along an annulus $A$ so as to obtain a solid torus $[0,1]\times[0,1]\times S^1$.
After adjusting $\Si'$ we can assume that it intersects $A$ transversally in circles.
Because $\Si'$ separates $\Si\times[0,1]$, every circle which is null-homotopic in $A$
bounds a 2-ball in $\Si'$.
(Otherwise, by the orientability and irreducibility of the solid torus,
such a circle would decompose $\Si'$ into an annulus,
and we could compress $\Si'$ to a circle or a point.)
Using irreducibility again, we can isotope $\Si'$
such that it intersects $A$ only in circles
which are not null-homotopic and decompose $\Si'$ into annuli.
Moreover, we can assume that these circles are vertical
in the fibration of the full torus by circles.
Hence every annulus component of $\Si' - A$
can be isotoped to be vertical, too. (Cf.\ \cite[2.4]{Waldhausen}.)
Thus, $\Si'$ can be isotoped to be vertical
in a fibration of $\Si \times [0,1]$ by circles.
This clearly implies (i).

(ii) Without loss of generality, we can assume that every boundary component of $\D \Si'$
is horizontal in the product $\Si \times [0,1]$.
In the case where $\Si$ is homeomorphic to $Ann^2$,
so is $\Si'$ and the two boundary components of $\Si'$
lie above different boundary components of $\Si$.

As above, we cut open $\Si \times [0,1]$ along a 2-ball $B$ to obtain a 3-ball $[0,1]\times[0,1]\times[0,1]$.
Again, we arrange that $\Si'$ is transversal to $B$ and hence intersects it
in null-homotopic circles
or in intervals connecting opposite sides of $B \cong [0,1]\times[0,1]$
or one side to itself.

Circles in $B\cap \Si'$ must again bound 2-balls in $\Si'$
and hence can be removed by suitable isotopies.
(If $\Si'$ is a Moebius band, this follows from the orientability of the 3-ball.
If $\Si'$ is an annulus which is decomposed into two annuli
by a component of $\Si' \cap B$,
we would obtain a compression disc for $\Si$
which is equally impossible.)

Similarly, it is impossible that an interval component of $\Si' \cap B$
connects one side of $B$ to itself:
This can only occur if $\Si$ was a Moebius band and $\Si'$ an annulus.
However, it follows in this case that $\Si'$
can be compressed to the $\D \Si \times [0,1]$.

This implies that without loss of generality $\Si' \cap B$ consists
of one or two intervals connecting opposite sides of the square $B$.
Since they decompose $\Si'$ into 2-balls,
claim (ii) now follows.
\qed

\subsubsection{Necks}
\label{sec:necorbi}

Throughout this section, we assume that $\theta < \theta_1$
and $\mu < \mu_0(\theta)$ are chosen sufficiently small,
and that $(O,g)$ is $(v(\theta,\mu),-b^2)$-collapsed
with $(v(\theta,\mu),s_0,K)$-curvature control
below scale $\rho_{-b^2}$.

We define a point $x\in O$ to be 
$(\theta,\mu,-b^2)$-{\em necklike} 
relative to an equilateral $<\theta$-straight 
1-strainer $(a,b)$ at $x$
with length in $(\hat s_{\mu,-b^2}(x), \frac{3}{2}\hat s_{\mu,-b^2}(x))$ 
if $\diam(\Si^o_{x;a,b})<\theta^2\hat s_{\mu,-b^2}(x)$, 
compare \ref{def:necky}. 
We call the open subset $N_{\theta,\mu,-b^2}\subseteq S_{\theta,\mu,-b^2}$ 
of $(\theta,\mu,-b^2)$-necklike points 
the {\em necklike region} of $O$. 
It does not contain any singular vertices, 
$O^{(0)}\cap N_{\theta,\mu,-b^2}=\emptyset$.

As in the beginning of the previous section,
we verify that a point $x \in O$ is $(\theta,\mu,-b^2)$-necklike
relative to a 1-strainer $(a,b)$
if and only if it is $\theta$-necklike
in the rescaled space
$(\hat{s}_{\mu,-b^2}(x))^{-1} \cdot B(x,\rho_{-b^2}(x))$
in the sense of Definition \ref{def:necky}.

Similarly, if $y \in B(x,\rho_{-b^2}(x))$
is $(\theta,\mu,-b^2)$-necklike
it is $\theta$-weakly necklike relative to $(a,b)$ in 
$(\hat{s}_{\mu,-b^2}(x))^{-1} \cdot B(x,\rho_{-b^2}(x))$.
If $y$ admits an equilateral $<\theta$-straight
1-strainer $(a,b)$ with length in $(\hat s_{\mu,-b^2}(y), \frac{3}{2}\hat s_{\mu,-b^2}(y))$
and is $\theta$-strongly necklike relative to $(a,b)$
in $(\hat{s}_{\mu,-b^2}(x))^{-1} \cdot B(x,\rho_{-b^2}(x))$,
it is also $(\theta,\mu,-b^2)$-necklike.

Thus we can again use the results from \ref{sec:edgalex}
to control the existence and relative position of necklike points.
\medskip

As in \ref{sec:edg},
we construct a smooth line field $L_{neck}$ 
tangent to the singular locus 
whose open domain of definition contains the balls 
$B(x,\theta^{\frac{3}{2}}\hat s_{\mu,-b^2}(x))$ 
around all $(\theta,\mu,-b^2)$-necklike points $x$, 
and which is on every such ball 
almost parallel to the equilateral 1-strainers $(a,b)$ at $x$ 
with length in $(\hat s_{\mu,-b^2}(x),\frac{3}{2}\hat s_{\mu,-b^2}(x))$, 
i.e.\ $\angle_{\cdot}(L_{neck},a),\angle_{\cdot}(L_{neck},b)<\theta^\half$ 
on that ball. 
The line fields $L_{neck}$ and $L_{edge}$ can be matched 
in the overlap of their domains of definition. 

%

For a $<\theta$-straight equilateral 1-strainer $(a,b)$
with length in $(\hat s_{\mu,-b^2}(x),\frac{3}{2}\hat s_{\mu,-b^2}(x)$, 
$\Si_{x;a,b}$ is a closed topological 2-suborbifold 
almost perpendicular to $L_{neck}$,
i.e.\ it has angle $>\pihalf-c''\theta$ with it. 
We call $\Si_{x;a,b}$ 
a {\em neck cross section} through the point $x$. 

If two neck cross sections 
$\Si_{x_1;a_1,b_1}$ and $\Si_{x_2;a_2,b_2}$ 
intersect, then they have Hausdorff distance 
$<c\theta^3\hat s_{\mu,-b^2}(x)$ 
by \ref{lem:almpcrsecnec}. 
If they are disjoint but also not too far apart from each other,
say if they have Hausdorff distance 
$<\theta^{\frac{5}{3}}\hat s_{\mu,-b^2}(x)$, 
then one can move one of the cross sections to the other 
along the trajectories of $L_{neck}$, 
and therefore the $\Si^o_{x_i;a_i,b_i}$
are in this case topologically {\em parallel},
i.e.\ they bound a product suborbifold $\cong\Si^o_{x_i;a_i,b_i}\times[0,1]$. 

Among all neck cross sections, 
we choose a maximal subfamily $\nu$ 
such that any two distinct cross sections 
$\Si_{x_1;a_1,b_1}$ and $\Si_{x_2;a_2,b_2}$ in $\nu$ 
have Hausdorff distance 
$>\theta^{\frac{5}{2}}\hat s_{\mu,-b^2}(x_1)$. 
In particular, they are disjoint. 
Due to the compactness of $O$, $\nu$ is finite. 
Inside $\nu$, 
we form equivalence classes of topologically parallel cross sections. 
Each equivalence class
has a linear or cyclic order.  
Any two successive cross sections in it 
are topologically parallel 
and bound a cylinder (``segment'') 
homeomorphic to the product of one of them with the compact interval. 
By concatenating these pieces,
the equivalence class yields an embedded {\em neck} in $O$ 
which fibers over the interval or over the circle, 
unless the equivalence class consists only of a single neck cross section. 
Such an isolated neck cross section 
has diameter $\approx\theta^2\hat s_{\mu,-b^2}(x)$ 
(with respect to some point $x$ in it); 
it is contained in the union of the 2-strained and edgy regions, 
and we simply disregard it. 
Any two necks are disjoint.
The union $N$ of all necks is a compact 3-suborbifold. 

A {\em cyclic} neck is a closed 3-orbifold and hence fills out $O$ entirely. 
The topology of cyclic necks will be determined later.

A {\em linear} neck has two boundary components 
which are neck cross sections. 
We call an end of the neck {\em thick} 
if its boundary $\Si_{x;a,b}$ has diameter 
$>\theta^{\frac{49}{20}}\hat s_{\mu,-b^2}(x)$, 
and {\em thin} otherwise. 
If the end is thin,
then nearby $\Si_{x;a,b}$, 
according to our construction e.g.\ at distance 
$<\theta^{\frac{49}{20}}\hat s_{\mu,-b^2}(x)$,
must exist points outside the 1-strained region $S_{\theta,\mu,-b^2}$, 
i.e.\ a $(\theta,\mu,-b^2)$-hump,
cf.\ \ref{prop:finhumporbi}. 
The interface between a thin end of a neck and a hump 
will be discussed later. 

\medskip

Let $\Si_{x;a,b}$ be a neck cross section with 
$\diam(\Si_{x;a,b})
>\theta^{\frac{49}{20}}\hat s_{\mu,-b^2}(x)$. 
Then 
every point in $\Si_{x;a,b}$ is $(\theta,\mu,-b^2)$-edgy
or belongs to the 2-strained region $R_{\theta,\mu,-b^2}$. 

Let $y\in\Si_{x;a,b}\cap R_{\theta,\mu,-b^2}$. 
Then $\Si_{x;a,b}$ contains a 
$C_1\theta$-straight 1-strainer $(z_1,z_2)$ at $y$ 
with length $\theta^4\hat s_{\mu,-b^2}(y)$, 
cf.\ \ref{lem:crosssecclsegnew} and \ref{lem:exstrainnew}(i). 
As discussed in section~\ref{sec:chafib}, 
the portion 
$A_y=\Si_{x;a,b}\cap
\{|f_{z_1,z_2}|\leq\tenth\theta^4\hat s_{\mu,-b^2}(y)\}
\cap B(y,\theta^4\hat s_{\mu,-b^2}(y))$
of the cross section 
fibers over a compact interval 
with fibers the $f_{z_1,z_2}$-level sets. 
These are embedded 1-suborbifolds $\cong S^1$ or $I^1$. 
Let $\ga_y=\Si_{x;a,b}\cap f_{z_1,z_2}^{-1}(0)=\Si^o_{y;a,b,z_1,z_2}$
denote the central fiber. 

To combine these local fibrations to a global one,
we choose a maximal family $F$ of $\ga_y$'s 
so that any two distinct $\ga_{y_1},\ga_{y_2}\in F$ have Hausdorff distance 
$>\frac{1}{100}\theta^4\hat s_{\mu,-b^2}(y_1)$. 
The family $F$ is finite, since $\Si_{x;a,b}$ is compact. 
If $\ga_{y_1}$ and $\ga_{y_2}$ have Hausdorff distance 
$<\frac{9}{100}\theta^4\hat s_{\mu,-b^2}(y_1)$, 
then $\ga_{y_2}$ separates $A_{y_1}$ 
and is isotopic inside $A_{y_1}$ (by a small isotopy) 
to a fiber of the above fibration of $A_{y_1}$.
We call $\ga_{y_1}$ and $\ga_{y_2}$ {\em adjacent},
if they are not separated inside  $A_{y_1}$ by another $\ga_y\in F$. 
In this case, they have Hausdorff distance 
$\approx\frac{1}{100}\theta^4\hat s_{\mu,-b^2}(y_1)$.
It follows that the $A_y$ for all $\ga_y\in F$ 
can be simultaneously isotoped (by small isotopies) 
so that their fibrations match afterwards.
This yields a fibration of part of $\Si^o_{x;a,b}$ 
and, if $\Si_{x;a,b}\subset R_{\theta,\mu,-b^2}$, a global fibration. 

If $e\in\Si_{x;a,b}-R_{\theta,\mu,-b^2}$ 
is a $(\theta,\mu,-b^2)$-edgy point, 
then the discussion in section~\ref{sec:edg} implies that 
$\D B(e,\half\theta^3\hat s_{\mu,-b^2}(e))$ lies in $R_{\theta,\mu,-b^2}$
and can be slightly isotoped inside $\Si_{x;a,b}$ 
to match the fibration obtained so far
or, vice versa, the fibration can be adapted 
so that $\D B(e,\half\theta^3\hat s_{\mu,-b^2}(e))$ 
becomes a fiber. 
Let us call $\ol B(e,\half\theta^3\hat s_{\mu,-b^2}(e))$
a {\em cap} of $\Si_{x;a,b}$. 
Since $\Si_{x;a,b}$ is compact and connected, 
it must have two disjoint caps which contain all 
$(\theta,\mu,-b^2)$-edgy points. 

In the case when the neck cross section lies entirely 
in the 2-strained region, 
$\Si_{x;a,b}\subset R_{\theta,\mu,-b^2}$, 
we can sandwich it between two nearby neck cross sections 
and proceed as in section~\ref{sec:edg} (for $S^2_{\kappa}$)
to isotope it 
by an isotopy supported nearby (i.e.\ in the sandwich) 
so that it becomes vertically saturated. 

If the neck cross section $\Si_{x;a,b}$ contains edgy points 
and hence two caps,
then we can coordinate the fibration of $\Si_{x;a,b}$ 
with the fibration of the tubes $T^j$ along the coarse edges, 
cf.\ section~\ref{sec:edg}, 
so that the intersection $T^j\cap\Si_{x;a,b}$
consists of two fibers of $T^j$, 
namely one for each cap of $\Si_{x;a,b}$. 
(Here, we refer to the fibration of $T^j$ 
by compact 2-orbifolds with one boundary component.) 
This is achieved by 
perturbing $\Si_{x;a,b}$
by a small isotopy
(using the flow of the vector field $L_{neck}$)
based near $T^j\cap\Si_{x;a,b}$
until it coincides with
the closest tube cross sections
of the $T^j$.
Moreover,
by a small isotopy of the fibration on $\Si_{x;a,b}$ minus the two caps, 
we can arrange that the intersections 
$S^j\cap\Si_{x;a,b}$ are fibers of the fibration of $\Si_{x;a,b}$. 
(In the latter step, we just use that any two noncontractible 
simple closed curves in an annular 2-orbifold are isotopic.) 

\subsubsection{Humps}
\label{sec:humpint}

We keep our assumption that 
$(O,g)$ is $(v(\theta,\mu),-b^2)$-collapsed
with $(v(\theta,\mu),s_0,K)$-curvature control
below scale $\rho_{-b^2}$
for sufficiently small $\theta,\mu>0$.
Then the discussion of 
sections~\ref{sec:2str}, \ref{sec:edg} and \ref{sec:necorbi} applies. 

Let $x\in O$ be a $(\theta,\mu,-b^2)$-hump
as defined after \ref{prop:locapproxconeorbi}.
This means that on the scale
$s(x)\in[s_1(\si,\mu)\rho_{-b^2}(x),\si\rho_{-b^2}(x)]$
of uniform conical approximation
provided by \ref{prop:locapproxconeorbi},
$O$ is $\mu$-well approximated in $x$
by a flat disc of radius 1 with cone point of angle $\leq2\pi-\theta$
or by a flat sector of radius 1 with angle $\leq\pi-\frac{\theta}{2}$.
(This includes the half-open interval $[0,1)$
as the degenerate case of the disc with cone angle $0$.)

The closed ball $\ol B(x_i,\half s(x_i))$ is a compact 3-suborbifold, 
since its boundary is almost orthogonal to 
radial (with respect to $x$) 
$\frac{\theta}{11}$-straight 1-strainers 
of length $\frac{1}{11}s(x)$, 
cf.\ \ref{lem:nhump}, 
and hence a closed 2-suborbifold. 

If $\diam(\D B(x,\half s(x)))$ is not too small, e.g.\ if 
$\diam(\D B(x,\half s(x)))
>\theta^{\frac{49}{20}} \hat{s}_{\mu,-b^2}(x)$, 
then all points in $\D B(x,\half s(x))$ are edgy or 2-strained and, 
as above in section~\ref{sec:necorbi} for neck cross sections,
we can construct a 1-dimensional fibration 
on most or all of $\D B(x,\half s(x))$. 
If the conical approximation in $x$ is by a disc with a cone point, 
then $\D B(x,\half s(x))\subset R^{O}_{\theta,\mu,-b^2}$ 
and the fibration is global. 
If the approximation is by a sector, 
then $\D B(x,\half s(x))$ contains edgy points 
close to the edges of the sector; 
these and the complement of the fibered region are covered by two caps 
whose boundaries are fibers. 
Since the edge cross sections 
($\Si^o_{x;a,b}$, cf.\ section~\ref{sec:edg}) 
associated to edgy points in $\D B(x,\half s(x))$ 
are also almost orthogonal to the radial direction (with respect to $x$), 
they can be embedded into $\D B(x,\half s(x))$ 
using an almost radial gradient like flow. 
As in section~\ref{sec:necorbi}, 
we can coordinate the fibration of $\D B(x,\half s(x))$ 
with the fibration of the tubes $T^j$ along the coarse edge 
so that the intersection $T^j\cap\D B(x,\half s(x))$ 
consists of two fibers of $T^j$, 
one for each cap of $\D B(x,\half s(x))$, 
and the intersections 
$\D B(x,\half s(x))$ are fibers of the fibration of 
$\D B(x,\half s(x))$.
We say that the hump $x$ has a {\em thick end}.

On the other hand, 
if $\diam(\D B(x,\half s(x)))$ is not too large, e.g.\ if 
$\diam(\D B(x,\half s(x)))
<\theta_i^{\frac{401}{200}} \hat{s}_{\mu,-b^2}(x)$,
then $\D B(x_i,\half s(x))$ 
is via an almost radial flow isotopic to a neck cross section 
associated to a $\frac{\theta}{11}$-straight 1-strainer 
of length $>\eleventh s(x)\geq\eleventh \hat{s}_{\mu,-b^2}(x)$
at a point in $\D B(x,\half s(x))$, 
and we have a {\em neck-hump} interface.
Such an interface corresponds to a 
{\em thin} end of a neck (as defined in section \ref{sec:necorbi}).

\medskip

We now have constructed a covering of the orbifold $O$
by finitely many humps $B(x_i,\half s(x_i)$,
necks, tubes and the fibration of $R_{\theta,\mu,-b^2}$.
This follows from Lemma \ref{lem:closetoedgy}:
Consider a point $x$ which is not contained in the balls
$B(x_i,\frac{1}{4}s(x_i)$ for the humps $x_i$;
it admits an equilateral $<\theta$-straight 1-strainer
with length $(\hat{s}_{\mu,-b^2}(x),\frac{3}{2}\hat s_{\mu,-b^2}(x))$.
If the diameter of its cross section is $\ge \theta^{\frac{9}{4}}$
and $x \not\in R_{\theta,\mu,-b^2}$,
the lemma implies that $x$
is contained in a tube with no end near $x$.

We now adjust the boundaries of humps, necks and tubes
to our fibration of the 2-strained part $R_{\theta,\mu,-b^2}$.
We have already done this in sections \ref{sec:edg} and \ref{sec:necorbi}
for cyclic chain of edgy points and for thick ends of necks
which contain no edgy points.
We now proceed analogously for $\D B(x_i,\half s(x_i))$
for all $(\theta,\mu,-b^2)$-humps with
$\diam(\D B(x,\half s(x_i)))
>\theta^{\frac{49}{20}} \hat{s}_{\mu,-b^2}(x_i)$
and conical approximation by a disc with a cone point.

At this point, we discard all necks which are entirely contained
in the union of humps, tubes and $R_{\theta,\mu,-b^2}$.
By our previous discussion,
every thick end of a hump or a neck meets $R_{\theta,\mu,-b^2}$.
If they contain no edgy points, we have already isotoped them
(by an isotopy supported nearby)
to a vertically saturated 2-suborbifold.

Every thick end (of a hump or a neck) containing edgy points
meets precisely two linear tubes along the coarse edge.
By our discussion in section \ref{sec:edg},
these tubes can only end deep inside a neck or a hump,
and hence a finite time after leaving the end
intersect another hump or neck in a thick end.
By the compactness of $(O,g)$,
such a sequence of tubes and thick ends
must eventually close up.
The union $C$ of all humps, necks and tubes
contained in such a closed chain
is a topological 3-suborbifold
with one boundary component $\D_0 C$
a topological 2-suborbifold homeomorphic to
$T^2$, $K^2$, $Ann^2$ or $M\ddot ob^2$,
as follows from our discussion
in sections \ref{sec:edg} and \ref{sec:necorbi}.
Of course, $C$ may have other boundary components
if it contains at least one neck.

For every chain $C$, we now extend the suborbifolds $S^i_j$
(for the different tubes $T_j \subset C$)
to closed 2-suborbifolds $S^i_C$ isotopic (in $C$) to $\D_0 C$,
e.g.\ by forming suitable unions
with the first three neck cross sections
of a neck in $C$,
and similarly for humps.
The 2-suborbifold $S^2_C$ is then
contained in $R_{\theta,\mu,-b^2}$, and
we can apply our Waldhausen-like arguments
from section \ref{sec:edg}
to isotope it (in the region between $S^1_C$ and $S^3_C$)
so that it becomes vertically saturated
with respect to the fibration of $R_{\theta,\mu,-b^2}$.

After performing this isotopy,
we cut off the tubes $T_j$ by a suitable cross section
such that the fibrations on $S^2_j$
induced by $T_j$ and $R_{\theta,\mu,-b^2}$ match.
The complement of all humps, necks and tubes
is now a saturated subset
of the fibration of $R_{\theta,\mu,-b^2}$
(see remark after Definition \ref{def:edgy}).

We can cut off tubes by smooth cross sections
transversal to the singular locus
by using the uniqueness of differentiable structures
(compare the discussion in section \ref{sec:edg}.)
Then the components of our decomposition
of a 3-orbifold $(O,g)$
have piecewise smooth boundary
and their interiors are
disjoint open smooth 3-suborbifold.

\medskip

Let us sum up our progress so far:
For every $0 < \theta < \theta_1$
and $0 < \mu, < \mu_0(\theta)$
the following holds:
If a 3-orbifold $(O,g)$
is $(v(\theta,\mu),-b^2)$-collapsed
with $(v(\theta,\mu),s_0,K)$-curvature control,
it admits a decomposition
(according to its coarse stratification)
into topological 3-suborbifolds
with disjoint interiors
and piecewise smooth boundary,
namely into $(\theta,\mu,-b^2)$-humps,
necks, tubes and total spaces of orbifold fibrations
with 1-dimensional fibers.

%
\subsection{Local topology}
\label{sec:loctop}

In this section, we determine the topological structures
of the components of the decomposition
we constructed in the previous section.
More precisely, we will determine
the topology of cross sections
to tubes and necks
and the topological type of humps.

\subsubsection{Tube and neck cross sections}

In this section we prove
that after decreasing $\theta$
further if necessary,
we can control the topological type
of the cross sections to tubes and necks
in sufficiently collapsed 3-orbifolds.

The following proposition is related to
an argument in the appendix of \cite{FukayaYamaguchi};
see  also \cite[4.24]{MT2} for a simplification
 of the special case needed here.

\begin{prop}[Topology of edge cross sections]
\label{prop:topcredg}
There exists $\theta_2>0$ such that 
for $\theta\in(0,\theta_2]$ and $\mu\in(0,\mu_0(\theta)]$ holds: 

If $(O,g)$ is $(v(\theta,\mu),-b^2)$-collapsed 
with $(v(\theta,\mu),s_0,K)$-curvature control below scale $\rho_{-b^2}$, 
and if $x\in O$ is $(\theta,\mu,-b^2)$-edgy relative to a 
$<\theta$-straight 1-strainer $(a,b)$ with length in
$(\hat s_{\mu,-b^2}(x),\frac{3}{2}\hat s_{\mu,-b^2}(x))$, 
then the truncated cross section 
$\Si_{x;a,b}\cap\ol B(x,\half\theta^3\hat s_{\mu,-b^2}(x))$ 
is a connected compact 2-suborbifold with one boundary component 
and Euler characteristic $\chi\geq0$. 
\end{prop}
\proof
Let $-b_i^2\in[-1,0)$
and let $\theta_i,\mu_i$ be sequences of small positive numbers
$\theta_i\to0$ and $\mu_i \le \mu_0(\theta)$.
Suppose that the orbifolds $(O_i,g_i)$ are
$(v(\theta_i,\mu_i),-b_i^2)$-collapsed
with $(v(\theta_i,\mu_i),s_0,K)$-curvature control below scale $\rho_{-b^2_i}$,
and that the points $x_i\in O_i$ 
are $(\theta_i,\mu_i,-b^2_i)$-edgy
relative to 
$<\theta_i$-straight 1-strainers $(a_i,b_i)$ with lengths in
$(\hat s_{\mu_i,-b_i^2}(x_i),\frac{3}{2}\hat s_{\mu_i,-b_i^2}(x_i))$. 

We consider the neighborhoods of the points $x_i$ 
on the scales $\theta_i^3\hat s_{\mu_i,-b_i^2}(x_i)$. 
The rescaled pointed orbifolds 
$(\theta_i^3\hat s_{\mu_i,-b_i^2}(x_i))^{-1}\cdot(O_i,x_i)$ 
Gromov-Hausdorff subconverge (collapse) 
to a 2-dimensional Alexandrov space with curvature $\geq0$ 
which splits off a line. 
In view of \ref{lem:crosssecclsegnew}, 
this limit is the pointed flat halfplane with base point on the boundary. 
The cross sections $\Si_{x_i;a_i,b_i}$ 
converge to the cross sectional ray of the halfplane 
through the base point. 


Let $z_i\in\Si_{x_i;a_i,b_i}$ 
with $d(x_i,z_i)=\theta_i^3 \hat{s}_{\mu_i,-b^2_i}(x_i)$. 
Then $(a_i,b_i,x_i,z_i)$ is a $c\theta_i$-straight 2-strainer 
near the intersection 
$\ga_i=
\Si_{x_i;a_i,b_i}\cap\D B(x_i,\half\theta_i^3 \hat{s}_{\mu_i,-b^2_i}(x_i))$,
for instance 
in the $\theta_i^4 \hat{s}_{\mu_i,-b^2_i}(x_i)$-neighborhood of it. 
Note that 
$(\theta_i^4 \hat{s}_{\mu_i,-b^2_i}(x_i))^{-1}\diam(\ga_i)\to0$,
because 
$v(\theta_i,\mu_i)\leq(\theta_i^5s_1(\si,\mu_i))^3$. 

The following consideration applies for sufficiently large $i$. 
From the discussion in \ref{sec:chafib} 
we know that 
a neighborhood (of at least comparable size) of $\ga_i$ 
is fibered by the level sets of the $\R^2$-valued map 
$(f_{a_i,b_i},d(x_i,\cdot))$. 
In particular,
$\ga_i\subset\Si_{x_i;a_i,b_i}$ 
is a connected closed 1-suborbifold.
This fibration exists in fact on a larger region, 
for instance on a neighborhood of 
$\Si_{x_i;a_i,b_i}\cap 
\ol A(x_i,\hunth\theta_i^3 \hat{s}_{\mu_i,-b^2_i}(x_i),
\frac{99}{100}\theta_i^3 \hat{s}_{\mu_i,-b^2_i}(x_i))$. 
In particular,
$d(x_i,\cdot)$ yields a product fibration (topologically) of 
$\Si_{x_i;a_i,b_i}\cap 
\ol A(x_i,\hunth\theta_i^3 \hat{s}_{\mu_i,-b^2_i}(x_i),
\frac{99}{100}\theta_i^3 \hat{s}_{\mu_i,-b^2_i}(x_i))$ 
over a compact interval. 
We note that as a consequence e.g.\ 
$\Si_{x_i;a_i,b_i}\cap B(x_i,\3quart\theta_i^3 \hat{s}_{\mu_i,-b^2_i}(x_i))$ 
(deformation) retracts onto 
$\Si_{x_i;a_i,b_i}\cap
\ol B(x_i,\half\theta_i^3 \hat{s}_{\mu_i,-b^2_i}(x_i))$. 
Using the flow of a gradient like vector field $X_i$ 
for the 1-strainer $(a_i,b_i)$ as constructed in section~\ref{sec:gradlike}, 
we obtain a homotopy of 
$B_i=\ol B(x_i,\half\theta_i^3 \hat{s}_{\mu_i,-b^2_i}(x_i))$
into 
$\Si_{x_i;a_i,b_i}\cap B(x_i,\3quart\theta_i^3 \hat{s}_{\mu_i,-b^2_i}(x_i))$
relative $\Si_i=\Si_{x_i;a_i,b_i}\cap B_i$, 
and together with the retraction 
$\Si_{x_i;a_i,b_i}\cap B(x_i,\3quart\theta_i^3 \hat{s}_{\mu_i,-b^2_i}(x_i))
\to\Si_i$ 
a retraction $r_i:B_i\to\Si_i$.
It is a retraction in the orbifold sense 
because $X_i$ is tangential to the singular locus. 
We have $\D\Si_i=\ga_i$. 

If $B_i$ is homeomorphic to the closed 3-ball, 
then $\pi_1(\Si_i)\cong1$ due to the retraction $r_i$, 
and hence $\Si_i$ is a closed 2-disc. 
More generally, 
if $B_i$ is discal 
then $r_i$ lifts to an equivariant retraction 
$\tilde r_i:\tilde B_i\to\tilde\Si_i$ 
of manifold covers.
As before, it follows that $\tilde\Si_i$ is a 2-disc 
and thus $\Si_i$ is discal. 

Suppose now that (after passing to a subsequence) 
none of the $B_i$ is discal. 
We then determine the possible topological types of the $B_i$ 
using the Shioya-Yamaguchi blow-up argument. 
Since the $O_i$ also have 
$(v(\theta_i,\mu_i),s_0,K)$-curvature control below scale $\rho_{-b_i^2}$, 
according to our discussion in section~\ref{sec:sy2dlim}, 
$B_i$ is for large $i$ homeomorphic to the product $[0,1]\times\Si'_i$ 
of the compact interval with 
a connected compact 2-orbifold with Euler characteristic $\chi\geq0$ 
and one boundary component. 
Namely, 
after a suitable choice of base points, 
all blow-up limits are 3-dimensional of the form 
$(Y,y_0)=(\R\times W,(0,w_0))$ 
with a noncompact ${\mathcal C}^{10}$-smooth 2-orbifold $W$ with $\sec\geq0$, 
cf.\ \ref{prop:blowuplargerdim} and \ref{lem:smconv}. 
The topology of $W$ is restricted by the 
(orbifold version of the) Soul Theorem. 
If $soul(W)$ is a point, then $W$ is discal. 
If $\dim soul(W)=1$, 
then $W$ is a one-ended quotient of the flat cylinder $S^1\times\R$, 
cf.\ \ref{lem:1dsoul}. 
The relation between the topologies of $\Si'_i$ 
(for a subsequence yielding the blow-up limit)
and $W$ is that $W$ is homeomorphic to the interior of $\Si'_i$. 

Knowing that 
$B_i\cong [0,1]\times\Si'_i$ for large $i$, 
we derive the topology of the truncated cross sections $\Si_i$ 
using the embeddings 
$\Si_i\subset B_i\cong[0,1]\times\Si'_i$ and the retractions $r_i$ 
as before. 
If $\Si'_i$ is discal, then we saw above that also $\Si_i$ is discal 
(and $\Si_i\cong\Si'_i$). 
Otherwise, 
$\Si'_i$ is finitely covered by an annulus $\tilde\Si'_i$ 
and $r_i$ lifts to an equivariant retraction 
$\tilde r_i:[0,1]\times\tilde\Si'_i\to\tilde\Si_i$
of smooth finite covers. 
Since the composition 
$\pi_1(\tilde\Si_i)\to\pi_1(\tilde\Si'_i)\cong\Z
\buildrel(\tilde r_i)_*\over\to\pi_1(\tilde\Si_i)$ 
of induced maps of fundamental groups 
is the identity, 
it follows that $\pi_1(\tilde\Si_i)\cong\Z$ or 0
and $\Si_i$ is finitely covered by a 2-disc or an annulus. 
(We do not worry about excluding the case of the disc here.) 
\qed

\medskip

It follows from the proposition
and our discussion in section \ref{sec:necorbi}
that for sufficiently small $\theta, \mu >0$,
{\em thick} ends of necks
in $(v(\theta,\mu),-b^2)$-collapsed orbifolds
with $(v(\theta,\mu),s_0,K)$-curvature control
have Euler characteristic $\chi \ge 0$.
Thus, we alredy control the topological structure
of necks with at least one thick end.
The following result generalizes
this to arbitrary necks,
i.e. cyclic necks or linear necks with two thin ends.

\begin{prop}[Topology of neck cross sections]
\label{prop:topcrneck}
There exists $\theta_3>0$ such that 
for $\theta\in(0,\theta_3]$ and $\mu\in(0,\mu_0(\theta)]$ holds: 

If $(O,g)$ is $(v(\theta,\mu),-b^2)$-collapsed 
with $(v(\theta,\mu),s_0,K)$-curvature control below scale $\rho_{-b^2}$, 
and if $x\in O$ is $(\theta,\mu,-b^2)$-necklike relative to a 
$<\theta$-straight 1-strainer $(a,b)$ with length in
$(\hat s_{\mu,-b^2}(x),\frac{3}{2}\hat s_{\mu,-b^2}(x))$, 
then the corresponding cross section 
$\Si_{x;a,b}$ 
is a closed 2-suborbifold with one boundary component 
and Euler characteristic $\chi\geq0$. 
\end{prop}

\proof
Let $-b_i^2\in[-1,0)$
and let $\theta_i,\mu_i$ be sequences of small positive numbers
$\theta_i\to0$ and $\mu_i \le \mu_0(\theta)$.
Suppose that the orbifolds $(O_i,g_i)$ are
$(v(\theta_i,\mu_i),-b_i^2)$-collapsed
with $(v(\theta_i,\mu_i),s_0,K)$-curvature control below scale $\rho_{-b^2_i}$,
and that the points $x_i\in O_i$ 
are $(\theta_i,\mu_i,-b^2_i)$-necklike
relative to 
$<\theta_i$-straight 1-strainers $(a_i,b_i)$ with lengths in
$(\hat s_{\mu_i,-b_i^2}(x_i),\frac{3}{2}\hat s_{\mu_i,-b_i^2}(x_i))$.

We let $d_i < \theta_i^2$ denote 
the diameter of the cross sections
$(\hat{s}_{\mu_i,-b^2_i}(x_i))^{-1}\cdot \Si_{x_i;a_i,b_i}$
and rescale by $2d_i^{-1}$.
Then after passing to a subsequence,
the orbifolds $(\half d_i \hat{s}_{\mu_i,-b^2_i})^{-1}\cdot (O_i,x_i)$
converge to an Alexandrov space $(Y,y)$
of curvature $\ge 0$ which splits off a line.
The factor of $Y$ orthogonal to the line
is the limit of the rescaled cross sections
to the 1-strainers $(a_i,b_i)$,
and hence has diameter 1 and is a compact
Alexandrov space of curvature $\ge 0$
and dimension 1 or 2.

If $\dim(Y) = 3$, by Lemma \ref{lem:smconv}
$Y$ is a $\mathcal{C}^{10}$-smooth orbifold
and the convergence can be improved to $\mathcal{C}^5$-smooth.
It follows that $Y$ is isometric to $\Si \times \R$
for some closed 2-orbifold $\Si$
of Euler characteristic $\chi(\Si) \ge 0$
and diameter 2.
For sufficiently large $i$,
we have an embedding of $\Si$
into $O_i$ which is transversal to the
gradient-like vector field $X_i$
for the strainer $(a_i,b_i)$.
This implies that $\Si$ is isotopic
to $\Si_{x_i;a_i,b_i}$.

If $\dim(Y) = 2$, $Y$ must be isotopic
to $[-1,1] \times \R$ or $S^1 \times \R$.
Hence there are constants $\phi_i \to 0$
such that every point in $(\half d_i \theta_i \hat{s}_{\mu_i,-b^2_i})^{-1}\cdot (O_i,x_i) \cap B(x_i,100)$
either admits a $C_1\phi_i$-straight 2-strainer of length $\phi_i^4$
or lies within $\phi_i^3$ of a point $z$
admitting a $C_0 \phi_i$-straight 1-strainer $(a',b')$ of length 1
such that $\Si_{z;a',b'}$ has diamter $\ge \phi_i^\frac{5}{2}$
and admits no $\frac{\pi}{2}$-straight
1-strainer of length $\phi_i^4$.

If $Y$ is isometric to $[-1,1] \times \R$,
we construct two tubes of diameter $\approx \phi_i^3$
along $(\half d_i \theta_i \hat{s}_{\mu_i,-b^2_i})^{-1}\cdot (O_i,x_i) \cap B(x_i,100)$
as in section \ref{sec:edg}
corresponding to the two edges of $Y$.
As described in section \ref{sec:necorbi},
we can now decompose $\Si_{x_i;a_i,b_i}$
into an annular part admitting a fibration
by 1-dimensional orbifolds
and two caps isotopic
to cross sections of the two tubes.
Note that for any $v > 0$ 
and sufficiently large $i$,
the balls $(\half d_i \theta_i \hat{s}_{\mu_i,-b^2_i})^{-1}\cdot (O_i,x_i) \cap B(z,1)$
centered at points $z \in (\half d_i \theta_i \hat{s}_{\mu_i,-b^2_i})^{-1}\cdot (O_i,x_i) \cap B(x_i,10)$
uniformly have volume $< v$
with $(v,s_0,K)$-curvature control on scale 1
because of $1 << \rho_{-b^2_i}$.
Thus we can proceed as in the proof of Proposition \ref{prop:topcredg}
to deduce that for sufficiently large $i$
the cross sections to both tubes
are compact 2-orbifolds
with Euler characteristic $\chi \ge 0$
and one boundary component.
This implies $\chi(\Si_{x_i;a_i,b_i}) \ge 0$.

If on the other hand 
$Y$ is isometric to $S^1 \times \R$,
it follows as in section \ref{sec:necorbi}
that for sufficiently large $i$
the cross sections $\Si_{x_i;a_i,b_i}$
admit a global fibration
by embedded 1-dimensional orbifolds
and hence are toric (have Euler characteristic $\chi = 0$).
\qed

\subsubsection{Humps}
\label{sec:topstrhumps}

For the remaining part of the proof, we fix $\bar{\theta} > 0$ 
such that the results of section \ref{sec:decomposrs}
and Propositions \ref{prop:topcredg} and \ref{prop:topcrneck} apply.
Thus whenever $\mu < \mu_0(\bar{\theta})$ and
a 3-orbifold $(O,g)$ is $(v(\bar{\theta},\mu),-b^2)$-collapsed
with $(v(\bar{\theta},\mu),s_0,K)$-curvature control,
it admits a decomposition
according to its coarse stratification
and we have control over the cross sections
of all tubes and necks.

In order to determine the local topology of humps
we will improve the quality
of our conical approximations,
i.e.\ make $\mu$ sufficiently small.
Again, we will adjust the upper bound for $\mu$
in several steps.

We say that a $(\bar{\theta},\mu,-b^2)$-hump $x \in (O,g)$
is a {\em thick hump}
if $O$ can in $x$ be $\mu$-well approximated on scale $s(x)$
by a flat cone with a base
of diameter $\in (\frac{\pi}{4},\pi - \frac{\bar{\theta}}{2})$.
In particular, this excludes the case 
of conical approximation by 
the 1-dimensional cones $(-1,1)$ or $[0,-1)$.
A thick hump must have a thick end
in the sense of section \ref{sec:humpint}.

\begin{prop}[Topological type of thick humps]
\label{prop:thickhumps}
There exists $0 <\mu_1 < \mu_0(\bar{\theta})$ such that:

Let $0 < \mu < \mu_1$. 
If $(O,g)$ is $(v(\bar{\theta},\mu),-b^2)$-collapsed
with $(v(\bar{\theta},\mu),s_0,K)$-curvature control
and if $x \in O$ is a thick $(\bar{\theta},\mu,-b^2)$-hump,
then $B(x,\half s(x))$ is discal or solid toric.
\end{prop}

\proof
This is an application 
of the Shioya-Yamaguchi blow-up
discussed in section \ref{sec:sy2dlim}.
Let $b^2_i \in [-1,0)$ and $\mu_i \to 0$.
Suppose that the orbifolds $(O_i,g_i)$
are $(v(\bar{\theta},\mu_i),-b_i^2)$-collapsed
with $(v(\bar{\theta},\mu_i),s_0,K)$-curvature control
and that the points $x_i \in O_i$
are thick $(\bar{\theta},\mu_i,-b_i^2)$-humps,
i.e.\ $\mu_i$-well approximated on scale $s(x_i)$
by flat cones $C_i$ with bases
of diameter $\in (\frac{\pi}{4},\pi - \frac{\bar{\theta}}{2})$.

A subsequence of the cones $C_i$ converges
to a flat cone $C_\infty$ with a base
of diameter $\in [\frac{\pi}{4},\pi - \frac{\bar{\theta}}{2}]$.
It follows that a subsequence
of the rescaled balls
$s(x_i)^{-1}\cdot B(x_i,s(x_i))$
also converges to $C_\infty$
in the Gromov-Hausdorff sense.

Unless infinitely many of the balls $B(x_i,s(x_i))$ are discal,
by our results in section \ref{sec:sy2dlim}
there is a sequence of rescaling factors $\delta_i \to 0$
such that the sequence $(\delta_i s(x_i))^{-1}\cdot B(x_i,s(x_i)$
Gromov-Hausdorff subconverges to a 3-dimensional limit space
$(Y,y)$ of curvature $\ge 0$.
As discussed in \ref{sec:gen},
$(Y,y)$ is actually a $\mathcal{C}^{10}$-smooth 3-orbifold
and the convergence can be improved to $\mathcal{C}^5$-smooth.

The soul of the blow-up limit $Y$
must be a point or 1-dimensional
since it cannot be 2-dimensional by \ref{sec:sy2dlim}.
Hence $Y$ must be either discal or solid toric.
Again by \ref{sec:gen},
this implies that for sufficiently large $i$,
the balls $B(x_i,\half s(x_i))$ are also
either discal or solid toric.
\qed

\medskip

We can ``read off''
the topological type of a thick hump
from the components of the decomposition it intersects.
Let $\mu > 0$ be sufficiently small
and suppose that the orbifold $(O,g)$
is sufficiently volume collapsed
with curvature control
such that it admits a decomposition 
according to its coarse stratification
and that Proposition \ref{prop:thickhumps} holds.
Let $x \in O$ be a thick hump
of this decomposition.

If $O$ is in $x$ $\mu$-well approximated on scale $s(x)$
by a flat cone over a circle,
equivalently if $\D B(x,\half s(x))$
is contained in $R_{\bar{\theta},\mu,-b^2}$,
it follows that $\D B(x,\half s(x))$
admits a fibration by 1-dimensional fibers
and hence cannot be spherical.
This implies that the hump $x$
is a solid toric 3-suborbifold
bounded by a vertically saturated
component of a 1-fibered
component of the decoposition of $O$.

If the conical approximation of $O$ in $x$
is by a flat sector,
$\D B(x,\half s(x))$ intersects
precisely two tubes with cross sections
$\Si_1, \Si_2$.
Both cross sections have Euler characteristic
$\chi \ge 0$.
Since $\D B(x,\half s(x))$ can be decomposed
into a union of $\Si_1, \Si_2$ 
and an annular (1-fibered) component, we have
$\chi(\D B(x,\half s(x)) = \chi(\Si_1) + \chi(\Si_2)$. 
Thus if $\chi(\Si_1) = \chi(\Si_2) = 0$,
it follows that $\chi(\D B(x,\half s(x)) = 0$
and hence that the hump $x$
is again solid toric.
Conversely, if at least one
of the $\Si_i$ is discal,
the hump $x$ must be discal as well.

\medskip

We cannot expect that the arguments
from the proof of the last Proposition
also work for
not necessarily thick humps,
i.e.\ humps with conical approximation by
cones with a base of arbitrarily small diameter.

In this case, it is possible that
the $s(x_i)^{-1}\cdot B(x_i,s(x_i))$
collapse to a 1-dimensional cone
and that the rescaled blow-ups
$(\delta_i s(x_i))^{-1}\cdot B(x_i,s(x_i))$
only Gromov-Hausdorff converge to a
2-dimensional Alexandrov space $(Y,y)$
of non-negative curvature
(see section \ref{sec:sy2dlim}).
We will however see that in this case
we can again apply our arguments
from section \ref{sec:decomposrs}
to obtain a decomposition
of the humps $B(x_i,s(x_i))$
with respect to the scale $\delta_i s(x_i)$
such that no thin humps or necks
occur in the decomposition.
When investigating collapse to the 2-dimensional space $Y$,
we operate on the scale $\delta_i s(x_i)$
rather than on the natural curvature scale $\rho_{-b_i^2}$.
Equivalently, the rescaled orbifolds
$(\delta_i s(x_i))^{-1} \cdot (O_i,x_i)$
collapse to $Y$ on scale $1$.
Note that we have already encountered
a similar situation
(in a very restricted setting)
in the proof of Proposition \ref{prop:topcrneck}.

Throughout the following considerations,
we will always assume that $\rho_{-b^2} >> 1$.
This implies $\sec \ge -b^2 \ge -1$ on balls of radius 1.
Moreover, it means that
$(v,s_0,K)$-curvature control on scale $\rho_{-b^2}$
implies $(v,s_0,K)$-curvature control on scale $1$.

We ommit $-b^2$ in our notation
to indicate that we work on scale 1 rather than $\rho_{-b^2}$.
For instance, we say that a 3-orbifold is
$v$-collapsed at a point $p$
if $\vol B(p,1) < v$.
Similarly, for $\theta, \mu > 0$ we define the set
$R_{\theta,\mu}$ as the set of all points
admitting $<C_1 \theta$-straight 2-strainers
of length $> \theta^4 s_1(\si,\mu)$.
Similarly, we define $(\theta,\mu)$-edgy points,
$(\theta,\mu)$-necklike points
and $(\theta,\mu)$-humps.

We can adapt our previous results
to our new setting of collapse at scale 1.
More precisely, we have

\begin{lem}
\label{lem:scale1}
There are $\hat{\theta}>0$
and $0 < \hat{\mu} < \mu_0(\hat{\theta})$ such that:

Let $(O,g)$ be a 3-orbifold with $\sec \not\geq 0$ and $x \in O$.
Suppose that for some $R >0$
the orbifold $O$ is
on the ball $B(x,6R)$
$(v(\hat{\theta},\hat{\mu}))$-collapsed
with $(v(\hat{\theta},\hat{\mu}),s_0,K)$-curvature control
on scale $1 << \rho_{-b^2}$,
and that $\diam O \ge 6R$.
Then the following hold:

\begin{itemize}
 \item
The orbifold $O$ can in every $y \in B(x,5R)$ 
be $\hat{\mu}$-well approximated 
on some scale $s(y)\in[s_1(\si,\hat{\mu}),\si]$ 
by a cone of dimension 1 or 2.
 \item
There exists an open subset $U$,
$R_{\hat\theta,\hat\mu} \cap B(x,3R) \subseteq U \subset B(x,4R)$
such that every connected component of $U$ is the total space
of a smooth orbifold fibration 
with fiber $S^1$ or the mirrored interval $\bar{I}^1$, 
and all fibers have angle $<\nu$ 
with the almost vertical line field $L_{\hat\theta,\hat\mu}$. 
 \item
If $y \in B(x,4R)$ is $(\hat{\theta},\hat{\mu})$-edgy
relative to an equilateral
$\hat{\theta}$-straight 1-strainer $(a,b)$
with length in $(s_1(\si,\hat{\mu}),\frac{3}{2}s_1(\si,\hat{\mu}))$, 
then the truncated cross section 
$\Si_{y;a,b}\cap\ol B(y,\half\hat\theta^3 s_1(\si,\hat{\mu}))$
is a connected compact 2-suborbifold with one boundary component 
and Euler characteristic $\chi\geq0$. 
 \item
If $y \in B(x,4R)$ is a thick $(\hat{\theta},\hat{\mu})$-hump,
then $B(y,\half s(y))$ is discal or solid toric.
\end{itemize}
\end{lem}

\proof
The proof works exactly as for Propositions
\ref{prop:locapproxconeorbi}, \ref{prop:fib}, \ref{prop:topcredg} and \ref{prop:thickhumps}
since in all of these proofs
we rescale by the collapse scale anyway.
\qed

\medskip

We now return to our original discussion
of the topological structure of general humps
in a $(v(\bar\theta,\mu),-b^2)$-collapsed 3-orbifold
with $(v(\bar\theta,\mu),s_0,K)$-curvature control.

\begin{prop}
\label{prop:allhumps}
There exists $0 <\mu_2 < \mu_0(\bar{\theta})$ such that:

Let $0 < \mu < \mu_2$. 
If $(O,g)$ is $(v(\bar{\theta},\mu),-b^2)$-collapsed
with $(v(\bar{\theta},\mu),s_0,K)$-curvature control
and if $x \in O$ is any $(\bar{\theta},\mu,-b^2)$-hump,
then one of the following holds:
\begin{enumerate}
 \item
$B(x,\half s(x))$ is discal or solid toric.
 \item
$B(x,\half s(x))$ is has the topological type of $(\Si \times [-1,1])/\Z_2$
with $\Si$ a closed 2-orbifold with $\chi(\Si) \ge 0$
and $\Z_2$ operating as a reflection on $[-1,1]$.
 \item
$B(x,\half s(x))$ admits
a decomposition as in section \ref{sec:decomposrs}
into a 1-fibered part, tubes, humps
and precisely one neck
containing $A(x,\frac{1}{4} s(x), \half s(x))$.
The cross section of this neck
is a closed 2-orbifold with $\chi \ge 0$.
Finally, all the humps 
occuring in this decomposition are discal or solid toric.
\end{enumerate}
\end{prop}

\proof
Again, let $b^2_i \in [-1,0)$ and $\mu_i \to 0$.
Suppose that the orbifolds $(O_i,g_i)$
are $(v(\bar{\theta},\mu_i),-b_i^2)$-collapsed
with $(v(\bar{\theta},\mu_i),s_0,K)$-curvature control
and that the points $x_i \in O_i$
are any $(\bar{\theta},\mu_i,-b_i^2)$-humps,
i.e.\ $\mu_i$-well approximated on scale $s(x_i)$
by flat cones $C_i$ with bases
of diameter $< \pi - \frac{\bar{\theta}}{2}$.

A subsequence of the cones $C_i$ converges
to some flat cone $C_\infty$ with a base
of diameter $< \pi - \frac{\bar{\theta}}{2}$,
and thus a subsequence
of the rescaled balls
$s(x_i)^{-1}\cdot B(x_i,s(x_i))$
also converges to $C_\infty$
in the Gromov-Hausdorff sense.

If the cone $C_\infty$ is 2-dimensional,
the proof proceeds as for Proposition \ref{prop:thickhumps}
to show that for sufficiently large $i$
we are in the first case of the proposition.

We now suppose that $C_\infty$ is 1-dimensional.
It is then isometric to
the half-open interval $[0,1)$ with cone point $\{0\}$.
Moreover we suppose that for infinitely many $i$
the ball $B(x_i,s(x_i))$ is not discal.

By our discussion in \ref{sec:sy2dlim}
we therefore have
rescaling factors $\delta_i \to 0$
such that the sequence $(\delta_i s(x_i))^{-1}\cdot B(x_i,s(x_i))$
Gromov-Hausdorff subconverges 
to a non-compact limit Alexandrov space $(Y,y)$
of curvature $\ge 0$
and dimension 2 or 3.

In case $\dim(Y) = 3$
$Y$ is again a $\mathcal{C}^{10}$-smooth orbifold
and we can improve the convergence
to $\mathcal{C}^5$-smooth.
Depending on the dimension of its soul,
$Y$ is discal, solid toric
or diffeomorphic to $(\Si \times [-1,1])/\Z_2$
with $\chi(\Si) \ge 0$
and $\Z_2$ operating on $[-1,1]$ by a reflection.
(We can exclude a product structure
since the $B(x_i,\half s(x_i))$ and hence $Y$ are one-ended.)
To finish this case,
we note again that
$B(x_i,\half s(x_i))$ is homeomorphic to $Y$
for sufficiently large $i$
by our discussion in \ref{sec:gen}.

We are now left with the
case where $C_\infty$ is isometric to $[0,1)$
and the pointed blow-up limit $(Y,y)$
is 2-dimensional.
Remember from \ref{sec:gen}
that we have a concave 1-Lipschitz function
$\beta_\xi$ on $Y$ coming from 
the unique direction at $\{0\} \in [0,1)$.
By construction, $y$ is a maximum of $\beta$
with $\beta(y) = 0$.

We observe two important properties
of the space $Y$
with respect to its curvature bound $\ge 0$:

\medskip

(i) For every point $z \in Y$,
there is a $\frac{\pi}{2}$-straight
1-strainer of length $\half$
centered at $z$.

More precisely,
let $\rho_\xi^z$ be a ray
of maximal $\beta_\xi$-decay
emanating from $z$
and let $y' \in Y$ be a maximum
of $\beta_\xi$,
i.e.\ $\beta_\xi(y') = 0$.
Recall that by construction,
there is a critical point $x$
at distance 1 from $y$
with $\beta_\xi(x) = 0$.
Concavity of $\beta_\xi$ implies
$\beta_\xi = 0$
on the whole segment $yx$ of length 1.
Hence we can choose $y'$
such that $d(z,y') \ge \half$.
Since $\beta_\xi$ is 1-Lipschitz,
we have $d(\rho_\xi^z(t),y') \ge \vert \beta_\xi(\rho_\xi^z(t)) \vert \ge t$
which implies that
for arbitrarily small $\eps > 0$
and sufficiently large $t > t_0(\eps)$
we have $\cangle_z (y',\rho_\xi^z(t)) \ge \frac{\pi}{2} - \eps$.
This implies property (i).

\medskip

(ii) There is a radius $R$
(depending on $\bar{\theta}$ and $Y$)
such that for every point $x \in Y$
with $d(z,y) = r \ge R$
there is a $\bar{\theta}$-straight
1-strainer $zxz'$ of length $r$.

Otherwise, we could find a sequence of points
$z_i \to \infty$
with $d(z_{i+1},y) \ge 2d(z_i,y)$
such that $\cangle_{z_i}(y,z_j) \le \pi - \frac{\bar\theta}{4}$
for $i < j$.
After passing to a subsequence,
we can assume that moreover
$\cangle_{z_j}(y,z_i) \le \frac{\bar\theta}{8}$:
We only have to make sure
that $d(y,z_i)$ is growing sufficiently fast,
i.e.\ $d(y,z_{i+1}) \ge \la d(y,z_i)$
for some $\la(\bar\theta)$.
But this implies
that $\angle_y (x_i,x_j) \ge \cangle_y (x_i,x_j) \ge \frac{\bar\theta}{4}$
for all $i \ne j$
which is absurd.

\medskip

Since we have by construction
$(\delta_i s(x_i))^{-1} \cdot B(x_i,s(x_i)) \to (Y,y)$,
for sufficiently large $i$
the $6R$-balls in $(\delta_i s(x_i))^{-1}\cdot B(x_i,s(x_i))$
are $v(\hat{\theta},\hat{\mu})$-collapsed
with $(v(\hat{\theta},\hat{\mu}),s_0,K)$-curvature control
on scale $1 << \rho_{-b^2}$.
Thus Lemma \ref{lem:scale1} applies to these balls.

In particular, for these $i$ every point
$z \in B_{(\delta_i s(x_i))^{-1}g_i}(x_i,6R)$
can be $\hat\mu$-well be approximated
on scale $s(z) \in [s_1(\si,\hat\mu),\si]$
by a flat cone.

If we make $i$ sufficiently large
(so that $d_{GH} (B_{(\delta_i s(x_i))^{-1}g_i}(x_i,6R), B(y,6R))$
becomes sufficiently small),
we also have conical approximation
on the ball $B(y,6R) \in Y$
of slightly lower quality,
say $2\hat\mu$-good approximation.
On the other hand,
it follows from properties (i) and (ii)
and the fact that $Y$ contains a flat strip of width 1
(see \ref{sec:sy2dlim})
that the diameters of approximating cones
must be $> \frac{\pi}{4}$ on $B(y,6R)$
and $> \pi - \frac{\hat\theta}{2}$
on $A(y,R,6R)$.
Hence for $i$ sufficiently large
we can deduce that the same bounds
on the diameters of approximating cones hold on the balls
$B_{(\delta_i s(x_i))^{-1}g_i}(x_i,6R)$.

We can now apply Proposition \ref{prop:finhumpalex}
to obtain a decomposition of $B_{(\delta_i s(x_i))^{-1}g_i}(x_i,6R)$
into finitely many $(\hat\theta,\hat\mu)$-humps
and the set $S^i_{\hat\theta,\hat\mu}$
of points admitting 
$\hat\theta$-straight 1-strainers with length in
$(\frac{1}{11} s_1(\si,\hat\mu), \frac{3}{22}\frac{1}{11} s_1(\si,\hat\mu))$.
Note that all humps must lie in $B_{(\delta_i s(x_i))^{-1}g_i}(x_i,R)$
and are {\em thick}.

We now proceed as in Lemma \ref{lem:scale1}
and construct a covering of 
$B_{(\delta_i s(x_i))^{-1}g_i}(x_i,3R)$
by the total spaces of fibrations with 1-dimensional fibers,
tubes and thick humps as in \ref{sec:decomposrs}.
Note that all occuring humps are
discal or solid toric.

By construction, the region $N = A_{(\delta_i s(x_i))^{-1}g_i}(x_i,2R,6R)$
is homeomorphic to a product of the interval and
$\D B(x_i,s(x_i))$.
We add it as a ``neck'' to our covering
of $B_{(\delta_i s(x_i))^{-1}g_i}(x_i,3R)$
and note that all points on its inner boundary
$\D B_{(\delta_i s(x_i))^{-1}g_i}(x_i,2R)$
are edgy or 2-strained.
Moreover, property (ii) of $(Y,y)$
implies that we can find
$<\theta$-straight 1-strainers
of length in $(\frac{1}{11}s_1(\si,\hat{\mu}),\frac{3}{22}s_1(\si,\hat{\mu}))$ 
which are almost orthogonal to $\D B_{(\delta_i s(x_i))^{-1}g_i}(x_i,2R)$
at every point
$z \in \D B_{(\delta_i s(x_i))^{-1}g_i}(x_i,2R)$
As in \ref{sec:humpint},
this allows us to construct a fibration on most or all of
$\D B_{(\delta_i s(x_i))^{-1}g_i}(x_i,2R)$,
which is either global
or capped off by the cross sections
of two tubes.
In other words, we treat
$\D B_{(\delta_i s(x_i))^{-1}g_i}(x_i,2R)$
as a thick end of the ``neck'' $N$.

We now adjust the interfaces
in $B_{(\delta_i s(x_i))^{-1}g_i}(x_i,3R)$
as in \ref{sec:humpint}
and extend $N$ radially away from $x_i$
to $\D B(x_i,\half s(x_i))$
using an gradient-like vector field
for $d(x_i,.)$.
(Recall that there are no critical points
for $x_i$ at distance $> \delta_i s(x_i)$.)
Thus for sufficiently large $i$
we have constructed a decomposition
as in the third case of the proposition.
\qed

%
\subsubsection{Proof of the main result}

In this section we complete the proof of
Theorem \ref{thm:mainwdb}.
As explained in section \ref{sec:setup},
Theorem \ref{thm:mainwdb} together with
Corollary \ref{cor:graphgeom}
implies our main result Theorem \ref{thm:mainclosed}.

In addition to $\bar\theta > 0$
from section \ref{sec:topstrhumps},
we fix some $0 < \bar\mu \le \mu_2$
with $\mu_2$ as in Proposition \ref{prop:allhumps}.
We set $v = v(\bar\theta,\bar\mu)$.
Then our discussion so far implies the following:

Let $(O,g)$ be a closed connected 3-orbifold
with $\sec \not\geq 0$
and $\rad(O) \ge \half \rho_{-b^2}$
for some $-b^2 \in [-1,0)$
and which contains no bad 2-suborbifold.
Suppose that $(O,g)$ is
$(v,-b^2)$-collapsed
with $(v,s_0,K)$-curvature control
on scale $\rho_{-b^2}$.
Then $O$ admits a decomposition
according to its coarse stratification
into finitely many components
of the following kind:
total spaces of orbifold fibrations
with 1-dimensional fibers,
tubes and necks with cross sections of 
Euler characteristic $\chi \ge 0$,
and humps which are
solid toric, 3-discal
or homeomorphic to
$(\Si \times [-1,1])/\Z_2$
as in Proposition \ref{prop:allhumps}.

There are three possibilities for
every end of a neck or hump:
If the end is thin, a hump ends in a neck and vice versa.
If the end is thick and meets no tubes,
it intersects one of the components with 1-dimensional fibration.
In this case the end is toric and vertically saturated
with respect to this fibration.
Finally, there is the possibility that
a thick end meets precisely two tubes
and one of the components with 1-dimensional fibration.
The boundary component of such an end
can then be further decomposed
into cross sections of the two tubes
and an annular part between them
which again is vertically saturated
with respect to the 1-dimensional fibration.
The end can be spherical or toric
depending on the Euler characteristic
of the tube cross sections.

In order to complete the proof of Theorem \ref{thm:mainwdb}
it therefore suffices to show
that after performing a finite number of surgeries
such a decomposition
can be simplified to a graph decomposition.

Because $O$ admits no bad 2-suborbifolds,
the cross sections of all necks must be
spherical or toric.
We perform surgery across all necks
with spherical cross section.
If such a neck bounds a discal component
or one of type $(\Si \times [-1,1])/\Z_2$
for some closed spherical 2-orbifold $\Si$
(corresponding to a hump at a thin end of the neck)
the resulting summmand is a finite quotient of $S^3$.
Similarly, if $O$ is a cyclic neck
with spherical cross section,
it is decomposed by one surgery
into a finite quotient of the 3-sphere.
We also perform surgery along the boundaries
of all components of type $(\Si \times [-1,1])/\Z_2$
with spherical $\Si$
which are not adjacent to a neck,
i.e.\ coming from humps with thick ends.

The orbifold $O$ may now be disconnected.
We discard all spherical summands.
Every remaining summand admits a decomposition as above
without necks with spherical cross sections
or humps of type $(\Si \times [-1,1])/\Z_2$
for spherical $\Si$.

Let $V$ be a 3-discal component of this decomposition.
It meets two coarse edges of $O$;
let $T$ and $T'$ be the corresponding tubes.
(We do not exclude the case $T = T'$.)
Due to Euler characteristic reasons,
at least one of the tube cross sections $\Si_T$, $\Si_{T'}$
must be discal.
If both are discal, then the two cross sections
must be homeomorphic 
because otherwise $\D V$ would be a bad 2-suborbifold of $O$.
In this case, $V$ is homeomorphic to $\Si_T \times [0,1]$
and we can replace the union
$T \cup V \cup T'$
by a single tube with cross section $\Si_{T}$,
thereby simplifying the decomposition. 

Because of the finiteness of the decomposition of $O$,
after repeating this step 
a finite number of times we can assume that
no 3-discal component of the decomposition of $O$
meets two tubes with discal cross section.

Consider now a tube $T$
with discal cross section.
If $T$ is cyclic, it is homeomorphic
to a fibration over $S^1$
with discal fiber
and hence to a solid toric 3-orbifold with boundary.
If $T$ is linear, it ends in
two 3-discal components $V_1$ and $V_2$
such that the other tubes ending in the $V_i$
have {\em annular} cross section.
In this case, the union $V_1 \cup T \cup V_2$
is again solid toric,
cf.\ the discussion after Proposition \ref{prop:nonnegsec3orbi}.
By considering these solid toric suborbifolds
as components of our decomposition,
we therefore can assume
that all tubes occuring in the decomposition of $O$
have annular cross section.

We recall from sections \ref{sec:edg} and \ref{sec:humpint}
that for each remaining tube $T$ (which now must have annular cross section)
intersecting the total space $U$
of a fibration with 1-dimensional fiber,
the two fibrations of $U$ and $T$
match on the 2-suborbifold $T \cap U$.
Since every annular 2-orbifold
inherits an orbifold fibration with 1-dimensional fiber
from the fibration of the annulus by circles,
we can extend the fibration
of $U$ to a Seifert fibration of $T \cup U$.

We have now obtained
a decomposition along disjoint embedded toric 2-suborbifolds
into components which
are total spaces of orbifold Seifert fibrations,
solid toric suborbifolds,
necks with toric cross section
and components of type $(\Si \times [-1,1])/\Z_2$
with toric $\Si$.
This is a graph decomposition by definition
(cf. section \ref{sec:decos}).
The proof of Theorem \ref{thm:mainwdb} is now complete.
\qed

\section{An extension to the case with boundary}
In this section we extend the results of the previous one
to a somewhat larger class of
volume collapsed 3-orbifolds.

We define a {\em hyperbolic orbifold cusp}
to be a complete 3-orbifold with boundary
which is isometric to the quotient
of a horoball in hyperbolic 3-space
by a cocompact isometric group action.
Thus, a hyperbolic orbifold cusp
is diffeomorphic to $\Si^2 \times [0,\infty)$
for some toric orbifold $\Si^2$
(by Bieberbach's theroem).
With the construction of the Ricci flow with surgery in mind
(cf. \cite{Perelman} and \cite{KL_coll} for orientable manifolds),
we will consider hyperbolic orbifold cusps
with sectional curvature equal to $-\frac{1}{4}$.

\begin{defn}[Almost cuspidal ends]
A Riemannian 3-orbifold $(O,g)$ with boundary has
$(v,s_0)$-almost cuspidal ends
if for every component $C \subset \partial O$
there is a hyperbolic orbifold cusp $X_C$
such that the pairs $(N_{100}(C),C)$ and $(N_{100}(\partial X_C),\D X_C)$
habe distance $\le v$ in the $\mathcal{C}^{s_0}$-topology.
\end{defn}

The following theorem generalizes 
Theorem \ref{thm:mainclosed}
to locally volume collapsed 3-orbifolds
with almost cuspidal ends
(compare again \cite[Theorem 7.4]{Perelman},
\cite[Theorem 0.2]{MT2}
and \cite[Theorem 1.3]{KL_coll}).

\begin{thm}
\label{thm:mainbdry}
Let $s_0\in\N$ 
and let $K:(0,\om_3)\to(0,\infty)$ be a function. 
If $s_0$ is sufficiently large,
then there exists a constant 
$v_0=v_0(s_0,K)\in(0,\om_3)$ such that:

If $(O,g)$
is closed or compact with $(v_0,s_0)$-almost cuspidal ends,
is $(v_0,-1)$-collapsed,
has $(v_0,s_0,K)$-curvature control 
below the scale $\rho_{-1}$
and contains no bad 2-suborbifolds,
then $O$
is either closed and admits a metric with $\sec \geq 0$,
or satisfies Thurston's Geometrization Conjecture.
\end{thm}

%

\proof
Throughout the following proof, we choose
$s_0$, $\bar{\theta}$, $\bar{\mu}$ and $v = v(\bar{\theta},\bar{\mu})$
as in the proof of Theorem \ref{thm:mainwdb}.

If $(O,g)$ is closed, $(v,-1)$-collapsed,
has $(v,s_0,K)$-curvature control
below scale $\rho_{-1}$
and contains no bad 2-suborbifolds,
we have already shown that the theorem holds.

We therefore now suppose that
$(O,g)$ has at least one ($(v,s_0)$-cuspidal) end.
In this case, we first observe that
$\rho_{-1}(x) \approx 4$ near a boundary component $C$,
say on $A(C,10,90)$.
This means that there are points $x \in O$ with
$\diam O >> 2 \rho_{-1}(x)$.
Hence collapse to a point cannot occur
and we can work with $-b^2 = -1$.
(In other words, we do not need
to make use of the more general setting
of Theorem \ref{thm:mainwdb}.)

After decreasing $v$ if necessary,
we obtain that every point $x$ 
close to a cuspidal end
(again, say on $A(C,10,90)$
for a boundary component $C \subset \D O$) 
admits a $<\bar{\theta}$-straight 1-strainer
of length $\hat{s}_{\bar{\mu},-1}(x)$ 
almost orthogonal to level sets of $d(C,.)$.
In particular,
we conclude $A(C,10,90) \subset S_{\bar\theta,\bar\mu,-1}$.

We fix this new value of $v$ and suppose from now on
that $(O,g)$ is compact with $(v,s_0)$-almost cuspidal ends
and $(v,-1)$-collapsed,
has $(v,s_0,K)$-curvature control
below scale $\rho_{-1}$
and contains no bad 2-suborbifolds.

We now define the {\em cusped necks} of $O$
to be the closed sets $\overline{N}_{25}(C)$
for all boundary components $C \subset \D O$.
On the neighbourhoods $N_{90}(C)$
of the cusped necks,
we have smooth gradient-like vector fields $V_C$
for the distance function
$d(C,.)$.
Cusped necks are homeomorphic to
$\Si^2 \times [0,1]$
for some toric 2-orbifold $\Si^2$.
Throughout the following discussion,
we are only interested in the ends of cusped necks
which are not boundary components of $O$.

By construction, cusped necks are disjoint from each other;
they are also disjoint from the humps in
$O \setminus \bigcup_C N_{10}(C)$
by \ref{lem:nhump} (i).

We will now show how cusped necks
can be integrated in our 
decomposition of $O$
according to its coarse stratification
much like humps.
As with humps, we call the end of a cusped neck {\em thin}
if the diameter of $\{d(C,.) = 25\}$
is not too large, say $\le \bar\theta^{\frac{401}{200}}s_1(\bar\mu,-1)$.
(Remember that we have seen
$\rho_{-1} \approx 4$
on a large neighbourhood of the end,
so the above condition means
that the diameter
is of order $\theta^{\frac{401}{200}}\hat{s}_{\bar\mu,-1}$
for all points in this neighbourhood.)
A thin end of a cusped neck corresponds to the {\em thin end}
of a neck (as defined in \ref{sec:necorbi})
and the interface can be matched up
using he flow of $V_C$.

Similarly, we say that the end of a cusped neck is {\em thick}
if the diameter of $\{d(C,.) = 25\}$
is sufficiently large, say $\ge \bar\theta^{\frac{49}{20}}$.
In this case, we can proceed as we did for humps in section \ref{sec:humpint}
and construct a 1-dimensional fibration 
on all or almost all of $\{d(C,.) = 25\}$,
with the possible exception of two tube cross sections.
We also can perturb $\{d(C,.) = 25\}$ such that it
intersects the tubes (if there are any)
in tube cross sections (again, using the flow of $V_C$).
If the end  of a cusped neck intersects two tubes,
it follows immediately
that both tubes have {\em annular} cross sections.

We now proceed to construct a decomposition 
of $O$ as we did in the closed case.
After adjusting the interfaces of the different
components of the decomposition
(using our Waldhausen-type arguments)
and performing a finite number of surgeries
we again obtain components
which are spherical
or admit a further decomposition along
(piecewise smooth)
toric suborbifolds into pieces which are
orbifold Seifert fibrations,
solid toric suborbifolds,
necks with toric cross section
or components of type $(\Si \times [-1,1])/\Z_2$ with toric $\Si$.
(The new components coming from cusped necks of $O$
are topologically of the same kind
as necks with toric cross sections.)
These decompositions are again graph,
which by virtue of Corollary \ref{cor:graphgeom}
completes the proof of the theorem.
\qed


\begin{thebibliography}{BLP05}



\bibitem[BLP05]{BLP}
M.\ Boileau, B.\ Leeb, J.\ Porti,
{\em Geometrization of 3--dimensional orbifolds},
Ann.\ Math.\ (2) 162, No.\ 1, 195--290 (2005).


\bibitem[BMP03]{BMP}
M.\ Boileau, S.\ Maillot, J.\ Porti,
{\em Three-dimensional orbifolds and their geometric structures},
Panoramas et Synth\`eses, vol.\ 15, SMF 2003.




\bibitem[BBBMP10]{BBBMP}
L.\ Bessi\`eres, G.\ Besson, M.\ Boileau, S.\ Maillot, J.\ Porti,
{\em Collapsing irreducible 3-manifolds with nontrivial fundamental group},
Invent.\ Math.\ 179, No.\ 2, 435--460 (2010).

\bibitem[BBI01]{BBI}
D.\ Burago, Y.\ Burago, S.\ Ivanov, 
{\em A course in metric geometry}, 
AMS 2001.


\bibitem[BGP92]{BGP}
Y.\ Burago, M.\ Gromov, G.\ Perelman,
{\em A. D. Alexandrov spaces with curvature bounded below}, 
Russ. Math. Surv. 47, No.\ 2, 1-58 (1992); 
translation from Usp. Mat. Nauk 47, No.\ 2(284), 3-51 (1992).

\bibitem[Bo02]{Bonahon_geomstrf}
F. Bonahon, 
{\em Geometric structures on 3-manifolds}, 
Handbook of geometric topology, 
93-164 (2002)

\bibitem[BS85]{BS_seif}
F. Bonahon, L.C. Siebenmann,
{\em The classification of Seifert fibered 3--orbifolds},
in {\em Low dimensional topology} (ed.\ Roger Fenn), Cambridge Univ. Press, Cambridge 1985.

\bibitem[BS87]{BS_split}
F. Bonahon, L. C. Siebenmann,
{\em The characteristic toric splitting of irreducible compact 3--orbifolds}
Math.\ Ann.\ 278, No. 1--4, 441--479 (1987).




\bibitem[DL09]{act3mfs}
J.\ Dinkelbach, B.\ Leeb,
{\em Equivariant Ricci flow with surgery
and applications to finite group actions on geometric 3-manifolds}, 
Geometry \& Topology 13, No. 2, 1129-1173 (2009).

\bibitem[Du88]{Dunbar}
W.D.\ Dunbar, 
{\em Geometric orbifolds}, 
Rev. Mat. Univ. Complutense Madr. 1, No.1-3, 67-99 (1988).

\bibitem[Ep66]{Epstein}
D.B.A.\ Epstein,
{\em Curves on 2-manifolds and isotopies},
Acta Math., 115, 83-107 (1966). 


\bibitem[FL]{locstr}
D.\ Faessler, B.\ Leeb,
{\em A compactness property for Riemannian orbifolds}, 
in preparation. 


\bibitem[FY92]{FukayaYamaguchi}
K.\ Fukaya, T.\ Yamaguchi, 
{\em The fundamental groups of almost nonnegatively curved manifolds}, 
Ann. Math. (2) 136, No.2, 253-333 (1992).




\bibitem[Ha82]{Hamilton1}
R.\ Hamilton,
{\em Three-manifolds with positive Ricci curvature},
J. Differential Geom. 17, No. 2, 255-306 (1982).

\bibitem[Ha95]{Hamilton}
R.\ Hamilton, 
{\em A compactness property for solutions of the Ricci flow},
Am. J. Math. 117, No.3, 545-572 (1995). 





\bibitem[KL08]{KL_notes}
B.\ Kleiner, J.\ Lott,
{\em Notes on Perelman's papers},
Geometry \& Topology 12, 2587--2855 (2008).

\bibitem[KL10]{KL_coll}
B.\ Kleiner, J.\ Lott, 
{\em Locally collapsed 3-manifolds}, 
arXiv:1005.5106 (2010).

\bibitem[KL11]{KL_preprint}
B.\ Kleiner, J.\ Lott,
{\em Geometrization of three-dimensional orbifolds via Ricci flow},
preprint, January 2011.

\bibitem[Kn29]{Kneser}
H.\ Kneser
{\em Geschlossene Fl\"achen in dreidimensionalen Mannigfaltigkeiten}, 
Jahresbericht D. M. V. 38, 248-260 (1929). 


\bibitem[Lu01]{Lu}
P.\ Lu, 
{\em A compactness property for solutions of the Ricci flow on orbifolds},
Amer.\ J.\ Math.\ 123 (2001), 1103-1134.



\bibitem[MY80]{MeeksYau}
W.H.\ Meeks, S.T.\ Yau, 
{\em Topology of three dimensional manifolds 
and the embedding problems in minimal surface theory}, 
Ann. Math. (2) 112, 441-484 (1980). 

\bibitem[MT07]{MT1}
J.\ Morgan, G.\ Tian,
{\em Ricci Flow and the Poincar\'e Conjecture},
Mathematics arXiv,
arXiv:math.DG/0607607v2 (2007).


\bibitem[MT08]{MT2}
J.\ Morgan, G.\ Tian,
{\em Completion of the proof of the Geometrization Conjecture},
Mathematics arXiv,
arXiv:0809.4040 (2008).


\bibitem[Mu60]{Munkres}
J.\ Munkres,
{\em Obstructions to the smoothing of piecewise-differentiable manifolds},
Ann.\ of Math.\ 72, No.\ 2, 521--554 (1960).


\bibitem[Pe03]{Perelman}
G.\ Perelman, 
{\em Ricci flow with surgery on three-manifolds}, 
Mathematics arXiv,
arXiv:math.DG/0303109 (2003).




\bibitem[SY00]{SY}
T.\ Shioya, T.\ Yamaguchi,
{\em Collapsing three-manifolds under a lower curvature bound},
J.\ Diff.\ Geom.\ 56, 1--66 (2000).

\bibitem[SY05]{SY05}
T.\ Shioya, T.\ Yamaguchi,
{\em Volume collapsed three-manifolds with a lower curvature bound},
Math.\ Ann.\ 333, No.\ 1, 131--155 (2005).

\bibitem[Th97]{Thurston_book}
W.\ Thurston,
{\em Three-Dimensional Geometry and Topology},
Princeton 1997. 


\bibitem[Wa67]{Waldhausen}
F.\ Waldhausen, 
{\em Eine Klasse von 3--dimensionalen Mannigfaltigkeiten: I, II}, 
Invent. Math. 3, 308--333; 4, 87--117 (1967).


\bibitem[Wh61]{Whitehead}
J.H.C.\ Whitehead,
{\em Manifolds with transverse fields in euclidean space},
Ann. of Math. 73, No.\ 2, 154--212 (1961).


\bibitem[Ya91]{Yamaguchi}
T.\ Yamaguchi, 
{\em Collapsing and pinching under a lower curvature bound},
Ann.\ Math.\ 133, No.\ 2, 317--357 (1991).


\end{thebibliography}
\end{document}